# SMOOTH ERGODIC THEORY OF $\mathbb{Z}^d$-ACTIONS PART 1: LYAPUNOV EXPONENTS, DYNAMICAL CHARTS, AND COARSE LYAPUNOV MANIFOLDS

AARON BROWN AND FEDERICO RODRIGUEZ HERTZ

ABSTRACT. In the first part of this paper, we formulate a general setting in which to study the ergodic theory of differentiable $\mathbb{Z}^d$-actions preserving a Borel probability measure. This framework includes actions by $C^{1+\text{Hölder}}$ diffeomorphisms of compact manifolds. We construct intermediate and coarse unstable manifolds for the action and establish controls on their local geometry.

In the second part we consider the relationship between entropy, Lyapunov exponents, and the geometry of conditional measures for rank-1 systems given by a number of generalizations of the Ledrappier–Young entropy formula. In the third part, for a smooth action of $\mathbb{Z}^d$ preserving a Borel probability measure, we show that entropy satisfies a certain "product structure" along coarse unstable manifolds. Moreover, given two smooth $\mathbb{Z}^d$-actions—one of which is a measurable factor of the other—we show that all coarse Lyapunov exponents contributing to the entropy of the factor system are coarse Lyapunov exponents of the total system. As a consequence, we derive an Abramov–Rohlin formula for entropy subordinated to coarse unstable manifolds.

## 1. INTRODUCTION

The primary motivation for this series of papers is to establish a "product structure of entropy" formula as well as a "coarse Abramov–Rohlin formula" for measure-preserving, non-uniformly hyperbolic $\mathbb{Z}^d$-actions. These formulas appear in Corollary 13.2 and Theorem 13.6 of Part 3. To establish these results, it is necessary to generalize the main result of seminal papers of Ledrappier and Young [LY1, LY2]. Parts 1 and 2 establish the necessary background and generalizations of the classical entropy formulas needed for Part 3.

The results of [LY1, LY2] establish a relationship between the metric entropy of a $C^2$ measure-preserving diffeomorphism, its Lyapunov exponents, and certain geometric data associated to the measure. In [LY1], a certain rigidity of measures is proven extending the main result of [Led]. Recall that for a $C^1$ diffeomorphism $f\colon M \to M$ of a compact manifold $M$ and an ergodic, $f$-invariant measure $\mu$ with Lyapunov exponents $\lambda_i$, the Margulis–Ruelle inequality gives a bound

$$h_\mu(f) \le \sum_{\lambda_i > 0} \lambda_i m_i \tag{1}$$

where $m_i$ is the multiplicity of the exponent. In [LY1] it is shown that if equality holds in (1), then the measure $\mu$ has the following geometric property: the conditional measures of $\mu$ along unstable manifolds are absolutely continuous (and in fact equivalent) to the Riemannian volume along unstable manifolds. When the unstable manifolds have an algebraic structure relative to which the dynamics acts as automorphisms, this rigidity can be used to obtain additional invariance of the measure. This idea is used, for example, in [MT, EM, BRHW, BFH].

In [LY2], the defect of equality in (1) is explicitly explained. For every positive Lyapunov exponent $\lambda_i > 0$ we define a corresponding contribution to the entropy $h_\mu(f)$. The maximal contribution of each $\lambda_i$ to the entropy $h_\mu(f)$ is $m_i\lambda_i$. The main result of [LY2] is that the entropy contribution of each $\lambda_i$ is given by $\gamma_i\lambda_i$ where $\gamma_i \leq m_i$ is an explicit geometric quantity.

The results of [LY1, LY2] have been generalized to many other settings including the case of i.i.d. random dynamics [LQ, LY3] on a compact manifold $M$, more general skew product systems [BL, QQX], infinite dimensional systems [BY], and endomorphisms [Shu, Liu]. In all these settings, the underlying dynamics occurs on compact subsets and the dynamics is assumed $C^2$. In the the case non-compact manifolds, finiteness of entropy and Margulis–Ruelle inequalities for diffeomorphisms have recently been studied in [Riq2]; see also [Riq1] for an example where the Margulis–Ruelle inequality fails for a $C^\infty$ diffeomorphism of a non-compact manifold.

In Part 1 of this paper we introduce the general setting in which we work. In Part 2 we reprove the main results of [LY1, LY2] under somewhat more general hypotheses and for more general notions of entropy. These generalizations are needed for the results in Part 3 which in turn are needed for [BRHW] and [BFH]. However, we expect the results here should be of interest and use in many other settings. In addition to [LY1, LY2], our results are highly derivative of [Led, KSLP, LS, Hu].

Some particular generalizations of existing results are the following:

(1) The original proof of the entropy formulas [LY1, LY2] requires the dynamics to be $C^2$. For hyperbolic measures (i.e. measures without zero exponents) the results still hold for $C^{1+\beta}$ diffeomorphisms of compact manifolds (see [Led], and [BP] and [BPS, Appendix]).[1] From the main result of [Bro], we have that the entropy formulas are still valid in the presence of zero Lyapunov exponents for $C^{1+\beta}$ diffeomorphisms of compact manifolds.

(2) Similarly, the entropy formulas appearing in [LY1, LY2] (and most extensions mentioned above) are proved assuming the manifold $M$ is compact. In particular, the derivative and its Lipschitz (or Hölder) norm are assumed to be bounded from above and the injectivity radius of $M$ is assumed to be bounded from below. In the case of random dynamics, on assumes the norms are slowly growing. We replace compactness (and consequent boundedness) with slow degeneration of the derivative, its Hölder norm, and the injectivity radii along orbits. This is accomplished through the introduction of dynamical charts $\{\phi_x\}$ in Section 3.1.

(3) We define two new related notions of metric entropy: metric entropy of a transformation *subordinate to a measurable partition* and metric entropy of a transformation *subordinate to a measurable foliation*. We prove a formula analogous to that of [LY2] for entropy subordinated to a measurable partition and, as in [LY1], show that the entropy subordinate to a measurable foliation attains its maximal theoretical value only when the measure is absolutely continuous along the unstable component of the foliation.

(4) As an application of the notion of metric entropy of a transformation *subordinate to a measurable partition* we obtain a Ledrappier–Young entropy formula for the fiber entropy of a diffeomorphism given a measurable factor of the system.

---

[1] Although [BP, Theorem 14.1.18] asserts the entropy formula holds for all (not just hyperbolic) measures, the proof outline of [BP, Lemma 14.1.19] fails because the Lipschitzness of unstable holonomies inside center-unstable sets is never established in [BP] and does not follow from the arguments in [BPS, Appendix].

(5) The main motivation for the generalizations outlined above was to establish the results in Part 3 which in turn are needed for [BRHW] and [BFH]. The main results are a *product structure of entropy* and a *coarse Abramov–Rohlin formula* for measures invariant under actions of higher-rank abelian groups.

(6) In Section 3.1, we introduce certain hypotheses to overcome non-compactness of the underlying space. These hypotheses only require controls on the local $C^{1+\beta}$ dynamics of the system. Thus, it is natural to work in systems where only the dynamics localized to an open set of full measures is assumed to be smooth (or even continuous). Such a setting was introduce in [KSLP] for diffeomorphisms with singularities and discontinuities. We introduce and prove our results in a setting similar to that introduced in [KSLP] for actions of $\mathbb{Z}^d$.

## CONTENTS

## 2. COCYCLES OVER $\mathbb{Z}^d$-ACTIONS

Let $(X, \mu)$ be a standard probability space. Consider $\alpha\colon \mathbb{Z}^d \times X \to X$ an action of $\mathbb{Z}^d$ by measurable, invertible, measure-preserving transformations of $(X, \mu)$. That is, for every $n, m \in \mathbb{Z}^d$ and $x \in X$

(1) $\alpha(n, \alpha(m, x)) = \alpha(n + m, x)$;

(2) $\alpha(0, x) = x$;

(3) $\alpha(n, \cdot)\colon X \to X$ is measurable;

(4) if $\phi \in L^1$ then $\int \phi(x)\, d\mu(x) = \int \phi \circ (\alpha(n, x))\, d\mu(x)$.

The action $\alpha$ of $\mathbb{Z}^d$ on $(X, \mu)$ is moreover said to be *ergodic* if

(5) given $\phi \in L^1(\mu)$, if $\phi(x) = \phi(\alpha(n, x))$ for all $n \in \mathbb{Z}^d$ then $\phi$ is constant on a set of full measure.

At times we write $\alpha(n) \colon X \to X$ for the map $\alpha(n)(x) = \alpha(n, x)$. Additional smoothness hypotheses on $\alpha$ and $X$ will be imposed starting in Section 3 below.

Fix a standard basis for $\mathbb{Z}^d$ and equip $\mathbb{Z}^d$ with the norm $|(n_1, \dots, n_d)| = \max |n_i|$. All definitions and facts stated below are independent of the choice of norm on $\mathbb{Z}^d$.

## 2.1. Slowly increasing functions.

**Definition 2.1.** Given a measurable set $X_0 \subset X$, a measurable function $L \colon X \to [1, \infty)$ is *slowly growing on $X_0$* (over the action $\alpha$) if

$$\lim_{\tau \to \infty} \frac{1}{\tau} \sup_{|n| \le \tau} \{\log(L(\alpha(n, x)))\} = 0 \tag{2}$$

for all $x \in X_0$. $L$ is *slowly growing* if $X_0$ can be taken with $\mu(X_0) = 1$. $L$ is said to be *$\varepsilon$-slowly growing on $X_0$* if for every $x \in X_0$ and every $n \in \mathbb{Z}^d$ we have

$$L(\alpha(n, x)) \le e^{|n|\varepsilon} L(x)$$

and *$\varepsilon$-slowly growing* if $X_0$ can be taken with $\mu(X_0) = 1$.

We have the following

**Claim 2.2.** *Let $X_0 \subset X$, and consider a measurable function $L \colon X \to [1, \infty)$ that is slowly growing on $X_0$. Let $Y_0$ be the $\alpha$-orbit of $X_0$. Then for any $\varepsilon > 0$ there is a measurable function $\hat{L} \colon Y_0 \to [1, \infty)$ that is $\varepsilon$-slowly growing on $Y_0$ with $L(x) \le \hat{L}(x)$.*

*Proof.* Given $x \in Y_0$ define

$$\hat{L}(x) := \sup_{n \in \mathbb{Z}^d} e^{-|n|\varepsilon} L(\alpha(n, x)). \tag{3}$$

Then $\hat{L}$ is defined for $x \in X_0$ by (2). Moreover, for $k \in \mathbb{Z}^d$

$$\hat{L}(\alpha(k, x)) := \sup_{n \in \mathbb{Z}^d} e^{-|n|\varepsilon} L(\alpha(n + k, x))$$

$$\le \sup_{n \in \mathbb{Z}^d} e^{|k|\varepsilon - |n + k|\varepsilon} L(\alpha(n + k, x))$$

$$= e^{|k|\varepsilon} \hat{L}(x)$$

whence (3) is defined for every $x \in Y_0$. Moreover $\hat{L}$ has the desired properties. $\square$

Applying either the higher-rank pointwise ergodic theorem (see [Kre]) or adapting [BP, Lemma 2.1.5] we have

**Claim 2.3.** *Let $\alpha$ be an action of $\mathbb{Z}^d$ on $(X, \mu)$ and let $\phi \colon X \to [1, \infty)$ be a measurable function with $\log(\phi) \in L^d(\mu)$. Then $\phi$ is slowly growing.*

Claim 2.3 fails for $\log(\phi) \in L^1(\mu)$ and $d > 1$. Indeed, see the example in Remark 2.3.

## 2.2. Multiplicative ergodic theorem and Lyapunov exponents.
Let $\alpha$ be an action of $\mathbb{Z}^d$ on $(X, \mu)$. A $k$-dimensional linear cocycle defined over $\alpha$ is measurable function

$$\mathcal{A} \colon \mathbb{Z}^d \times X_0 \to \mathrm{GL}(k, \mathbb{R})$$

satisfying the cocycle relation: $\mathcal{A}(0, x) = \mathrm{Id}$ and

$$\mathcal{A}(m, \alpha(n, x))\mathcal{A}(n, x) = \mathcal{A}(n + m, x) \tag{4}$$

for all $n, m \in \mathbb{Z}^d$ and $x \in X_0$ where $X_0 \subset X$ is an $\alpha$-invariant subset with $\mu(X_0) = 1$.

Write $L^p(\mu)$ for the standard $L^p$ spaces. $L^{p,q}(\mu)$ denotes the Lorentz space introduced, for instance, in [Lor]. For $1 \leq p, q \leq \infty$, $\varepsilon > 0$ and $q < q'$ we have that

(1) $L^{p,p}(\mu) = L^p(\mu)$;
(2) $L^{p+\varepsilon}(\mu) \subset L^{p,1}(\mu)$;
(3) $L^{p,q}(\mu) \subset L^{p,q'}(\mu)$.

Given a function $\phi \colon X \to \mathbb{R}$ write $\log^+ \phi := \max\{\log \phi, 0\}$. Write $|\cdot|$ for the standard norm on $\mathbb{R}^k$ and $\|\cdot\|$ for induced operator norm.

**Theorem 2.4** (Higher-rank Oseledec's Theorem). *Let $\alpha$ be an ergodic action of $\mathbb{Z}^d$ on $(X, \mu)$ and let $\mathcal{A} \colon \mathbb{Z}^d \times X \to \mathrm{GL}(k, \mathbb{R})$ be a measurable cocycle. Assume for every $m \in \mathbb{Z}^d$ that*

$$\big(x \mapsto \log^+ \|\mathcal{A}(m, x)\|\big) \in L^{d,1}(\mu). \tag{5}$$

*Then there are*

*(1) an $\alpha$-invariant subset $\Lambda_0 \subset X$ with $\mu(\Lambda_0) = 1$;*
*(2) linear functionals $\lambda_i \colon \mathbb{Z}^d \to \mathbb{R}$ for $1 \leq i \leq p$;*
*(3) and splittings $\mathbb{R}^k = \bigoplus_{i=1}^p E_{\lambda_i}(x)$ into families of mutually transverse, $\mu$-measurable subbundles $E_{\lambda_i}(x) \subset \mathbb{R}^k$ defined for $x \in \Lambda_0$*

*such that*

*(a) $\mathcal{A}(n, x) E_{\lambda_i}(x) = E_{\lambda_i}(\alpha(n, x))$ and*
*(b)* $$\lim_{|n| \to \infty} \frac{\log |\mathcal{A}(n, x)(v)| - \lambda_i(n)}{|n|} = 0$$

*for all $x \in \Lambda_0$ and all $v \in E_{\lambda_i}(p) \smallsetminus \{0\}$. Moreover, for $x \in \Lambda_0$ we have*

*(c)* $\displaystyle\lim_{n \to \infty} \frac{\log |\det \mathcal{A}(n, x)| - \sum m_i \lambda_i(n)}{|n|} = 0$ *where $m_i$ is the almost-surely constant value of $\dim(E_{\lambda_i}(x))$, and*
*(d) for every $\lambda_i$*

$$\lim_{n \to \infty} \frac{1}{|n|} \log \Big( \sin \angle \Big( E_{\lambda_i}(\alpha(n, x)), \bigoplus_{\lambda_j \neq \lambda_i} E_{\lambda_j}(\alpha(n, x)) \Big) \Big) = 0.$$

The modifications of Theorem 2.4 to the case of non-ergodic $\mu$ are standard. Here $\angle(V, U)$ denotes the smallest angle between subspaces $V, U$ in $\mathbb{R}^k$ relative to the standard inner product; that is $\angle(V, U) = \inf\{\angle(v, u) : v \in V \smallsetminus \{0\}, u \in U \smallsetminus \{0\}\}$.

Observe that the limits in (b)–(d) are taken along any sequence $n \to \infty$ in $\mathbb{Z}^d$. The limit in (b) implies the following weaker result: given $x \in \Lambda_0$, $v \in E_{\lambda_i} \smallsetminus \{0\}$, and $n \in \mathbb{Z}^d$

$$\lim_{k \to \infty} \frac{1}{k} \log |\mathcal{A}(kn, x)(v)| = \lambda_i(n). \tag{6}$$

However, convergence along rays in (6) does not imply the limit in (b) holds. It seems the $L^{d,1}(\mu)$ hypothesis is sharp as discussed in [BD].

To prove Theorem 2.4, first using only that $\big(x \mapsto \log^+ \|\mathcal{A}(m, x)\|\big) \in L^1(\mu)$, for every $m \in \mathbb{Z}^d$ one can produce the set $\Lambda_0$, the splitting $\mathbb{R}^k = \oplus_{i=1}^p E_{\lambda_i}(x)$ for $x \in \Lambda_0$ and the functionals $\lambda_i \colon \mathbb{R}^d \to \mathbb{R}$ such that convergence along (the countably many) rays in (6) holds in (b)–(d) (see [BP, Theorem 3.6.6]).

The convergence (b) and (c) can be derived from the maximal lemma in [Boi]. Alternatively, the convergence in (b) and (c) can be reinterpreted in terms of random semimetrics modeled on $\mathbb{Z}^d$ in which case the result follows from [Björ, Theorem 2.4]. The limit in (d) follows from (c) using arguments as in [BP, Theorem 1.3.11] and [BP, Section 1.3.2].

**Definition 2.5.** The linear functionals $\lambda_i$ in Theorem 2.4 are called *Lyapunov exponent functionals*, or simply *Lyapunov exponents*. The subspaces $E_{\lambda_i}(x)$ are the *Oseledec's subspaces*, and the set $\Lambda_0$ is the set of *regular points*. $m_i$ is the *multiplicity* of $\lambda_i$.

We write $\mathcal{L} = \{\lambda_i\}$ for the set of Lyapunov exponents functionals. Note that the exponents $\lambda_i$ are independent of the choice of norm on $\mathbb{R}^k$ and the choice of generating set and norm on $\mathbb{Z}^d$.

We have the following standard lemma.

**Proposition 2.6.** *For any $\varepsilon > 0$ there are $\varepsilon$-slowly growing measurable functions $A, K \colon \Lambda_0 \to [1, \infty)$ such that for every $x \in \Lambda_0$*

(a) *for all $n \in \mathbb{Z}^d$,*
$$A(x)^{-1} e^{\lambda_i(n) - |n|\varepsilon/2} \leq \|\mathcal{A}(n, x)\!\restriction_{E_{\lambda_i}(x)}\| \leq A(x) e^{\lambda_i(n) + |n|\varepsilon/2};$$

(b) *for all $\lambda_i$,*
$$\sin\left(\angle\left(E_{\lambda_i}(x), \bigoplus_{\lambda_j \neq \lambda_i} E_{\lambda_j}(x)\right)\right) \geq K(x)^{-1}.$$

### 2.3. **Failure of Theorem 2.4 in $L^1$.**

The integrability hypothesis $(x \mapsto \log^+ \|\mathcal{A}(m, x)\|) \in L^{d,1}(\mu)$ seems to be sharp in Theorem 2.4. In the literature, there are assertions that $L^1$ integrability is sufficient. For instance, both [BP, Theorem 3.6.7] and [KN, Theorem 1.7.1] are incorrect as stated. If we only have $(x \mapsto \log^+ \|\mathcal{A}(m, x)\|) \in L^1(\mu)$ then the limit in (1.7.4) (corresponding to (b) of our Theorem 2.4) of [KN, Theorem 1.7.1] need not converge. Similarly, in [BP, page 87] the sum defining the Lyapunov metric need not converge and the reduction theorem of [BP, Theorem 3.6.7] fails.

Both these defects can be seen in the following example. Take $\beta_1, \beta_2 \in [0, 1] \smallsetminus \mathbb{Q}$ to be badly approximable numbers and let $r_\beta \colon S^1 \to S^1$ denote the rotation on $S^1$ by $\beta$. Let $\alpha \colon \mathbb{Z}^2 \times \mathbb{T}^2 \to \mathbb{T}^2$ be the product action $\alpha((n, m), (x, y)) = (r_{n\beta_1}(x), r_{m\beta_2}(y))$. Given any $\gamma \in (1/2, 1)$ let $\phi \colon \mathbb{T}^2 \to \mathbb{R}$ be the function
$$\phi(x, y) = \frac{1}{d(x, 0)^\gamma} \frac{1}{d(y, 0)^\gamma}$$
where $d$ denotes the distance on $S^1$. Let $\mathcal{A} \colon \mathbb{Z}^2 \times \mathbb{T}^2 \to \mathrm{GL}(1, \mathbb{R})$ be the (abelian) cocycle
$$\mathcal{A}((n, m), (x, y)) = \exp\left[\phi(\alpha((n, m), (x, y))) - \phi(x, y)\right].$$

Let $\mu$ be the Lebesgue measure on $\mathbb{T}^2$. We have
$$(x, y) \mapsto \log^+ \|\mathcal{A}((n, m), (x, y))\| \in L^1(\mu)$$
but
$$(x, y) \mapsto \log^+ \|\mathcal{A}((n, m), (x, y))\| \notin L^2(\mu).$$
As $\beta_1, \beta_2$ are chosen to be badly approximable, there is a $C$ such that, for any $x, y \in S^1$ and any $N$, there are $n, m \leq N$ with
$$d((r_{n\beta_1}(x), 0) \leq \frac{C}{N}, \quad d((r_{m\beta_2}(y), 0) \leq \frac{C}{N}.$$
Then for any $(x, y) \in \mathbb{T}^2$ one can find a sequence $(n_i, m_i) \to \infty$ with
$$\frac{\log^+ \|\mathcal{A}((n_i, m_i), (x, y))\|}{|(n_i, m_i)|} \to \infty.$$

### 2.4. **Restricted cocycles.**

Let $\alpha$ be an ergodic action of $\mathbb{Z}^d$ on $(X, \mu)$. Let $\mathcal{A} \colon \mathbb{Z}^d \times X \to \mathrm{GL}(k, \mathbb{R})$ be a linear cocycle satisfying the hypotheses of Theorem 2.4 and let $\mathcal{L} = \{\lambda_j\}$ be the Lyapunov exponents of $\alpha$. Let $H \subset \mathbb{Z}^d$ be a subgroup and let $\tilde{\alpha}$ denote the restriction of $\alpha$ to $H$. It may be that $\tilde{\alpha}$ is no longer ergodic. However, for almost every ergodic component $\tilde{\mu}^e_x$ of the action of $\tilde{\alpha}$ on $(X, \mu)$, the restriction of the cocycle $\mathcal{A}$ to the action $\tilde{\alpha}$ satisfies the hypotheses of Theorem 2.4. For such an ergodic component $\tilde{\mu}^e_x$, let $\tilde{\mathcal{L}}_x = \{\tilde{\lambda}_{i,x}\}$ be the Lyapunov exponents of the restricted cocycle. Clearly, for every $i$ we have $\tilde{\lambda}_{i,x} = \lambda_j {\restriction}_H$ for some $j$. In particular, the functionals $\{\tilde{\lambda}_{i,x}\}$ are independent of the ergodic component $\tilde{\mu}^e_x$.

### 2.5. Lyapunov metric.
A standard technique which simplifies certain dynamical arguments is to specify a family of inner products and related norms on $\mathbb{R}^k$ relative to which the dynamics of the cocycle $\mathcal{A}$ is uniformly partially hyperbolic. Let $\mathcal{A} \colon \mathbb{Z}^d \times M \to \mathrm{GL}(d, \mathbb{R})$ be a measurable cocycle satisfying (5). Given any $\varepsilon > 0$, $x \in \Lambda_0$, $\lambda_i \in \mathcal{L}$, and $v, w \in E_{\lambda_i}(x)$ define a measurable family of inner products, called the $\varepsilon$-*Lyapunov metric*, by

$$\langle\!\langle\!\langle v, w \rangle\!\rangle\!\rangle_{x,\varepsilon} := \sum_{n \in \mathbb{Z}^d} e^{-2\lambda_i(n) - 2\varepsilon|n|} \langle \mathcal{A}(n, x)v, \mathcal{A}(n, x)w \rangle. \tag{7}$$

The expression in (7) converges for $x \in \Lambda_0$ by Proposition 2.6(a). We extend $\langle\!\langle\!\langle \cdot, \cdot \rangle\!\rangle\!\rangle_{x,\varepsilon}$ to $\mathbb{R}^k$ be declaring

$$\langle\!\langle\!\langle v, w \rangle\!\rangle\!\rangle_{x,\varepsilon} = 0$$

for $v \in E_{\lambda_i}(x)$ and $w \in E_{\lambda_j}(x)$ with $\lambda_i \neq \lambda_j$.

For $v_i \in E_{\lambda_i}(x)$ let $\|\!\|v_i\|\!\|_{x,\varepsilon}$ be the norm on $E_{\lambda_i}(x)$ induced by $\langle\!\langle\!\langle \cdot, \cdot \rangle\!\rangle\!\rangle_{x,\varepsilon}$. Given $v \in \mathbb{R}^d$, we decompose $v = \sum v_i$ where $v_i \in E_{\lambda_i}(x)$ and define a measurable family of norms $\|\!\|\cdot\|\!\|_{x,\varepsilon}$, called the $\varepsilon$-*Lyapunov norm*, on $\mathbb{R}^k$ by

$$\|\!\|v\|\!\|_{x,\varepsilon} = \max\{\|\!\|v_i\|\!\|_{x,\varepsilon}\}. \tag{8}$$

We have the following two facts about the family of norms $\|\!\|\cdot\|\!\|_{x,\varepsilon}$.

**Proposition 2.7.** *For $x \in \Lambda_0$, $v \in E_{\lambda_i}(x)$, and all $k \in \mathbb{Z}^d$ we have*

$$e^{\lambda_i(k) - \varepsilon|k|} \|\!\|v\|\!\|_{x,\varepsilon} \leq \|\!\|\mathcal{A}(k, x)v\|\!\|_{\alpha(k,x),\varepsilon} \leq e^{\lambda_i(k) + \varepsilon|k|} \|\!\|v\|\!\|_{x,\varepsilon}. \tag{9}$$

**Lemma 2.8.** *There is uniform constant $k_0 > 1$ and, for every $\varepsilon > 0$, an $\varepsilon$-slowly growing function $L \colon \Lambda_0 \to [1, \infty)$ such that for $x \in \Lambda_0$ and $v \in \mathbb{R}^k$*

$$k_0^{-1}|v| \leq \|\!\|v\|\!\|_{x,\varepsilon} \leq L(x)|v|. \tag{10}$$

## 3. Smooth actions of $\mathbb{Z}^d$ with singularities and discontinuities

In this section we establish the notational conventions as well as the standing hypotheses for the remainder of the paper. In particular, we present hypotheses under which generalizations of the entropy formulas of [LY1, LY2] will hold.

### 3.0.1. *Notational conventions.*
Given a map $f \colon X \to Y$ between metric spaces let $\mathrm{Lip}(f)$ be the Lipschitz constant of $f$ and $\mathrm{H\ddot{o}l}^\beta(f)$ the $\beta$-Hölder constant of $f$. Given a norm $\|\cdot\|$ on $\mathbb{R}^m$ and subspace $V \subset (\mathbb{R}^m, \|\cdot\|)$ we write $V(r) = V(r, \|\cdot\|) := \{v \in V : \|v\| < r\}$.

Consider $\mathbb{R}^k$ and $\mathbb{R}^j$ equipped, respectively, with norms norms $\|\cdot\|_1$ and $\|\cdot\|_2$, an open set $U \subset \mathbb{R}^k$, and a differentiable map $g \colon U \to \mathbb{R}^j$. Let

$$\|Dg\| = \sup_{u \in U} \|D_u g\|$$

and let $\|g\|_{C^1}$ be the usual $C^1$ norm of $g$. Note we often consider the case $0 \in U \subset \mathbb{R}^k(1, \|\cdot\|_1)$ and $g(0) = 0$ and ignore the $C^0$ part of $g$. Considering $Dg$ as map from $U$ to a space of linear maps we have

$$\mathrm{H\ddot{o}l}^\beta(Dg) = \sup_{u \neq v \in U} \frac{\|D_u g - D_v g\|}{\|u - v\|_1^\beta}$$

for the $\beta$-Hölder constant of $Dg$. Set $\|g\|_{C^{1+\beta}} = \max\{\|g\|_{C^1}, \mathrm{H\ddot{o}l}^\beta(Dg)\}$. If $\|g\|_{C^{1+\beta}} < \infty$ we say $g$ is uniformly $C^{1+\beta}$.

### 3.1. Standing hypotheses.

We work in a setting similar to that introduced in [KSLP]. Let $M$ be a Hausdorff, second countable, $k$-dimensional, $C^\infty$ manifold without boundary. Let $\mu$ be a Borel probability measure on $M$. Unlike in [KSLP] we don't explicitly assume any properties of any metric on $M$ or the metric completion of $M$. Let $\alpha \colon \mathbb{Z}^d \times M \to M$ be an action by measurable, invertible, $\mu$-preserving transformations. We do not assume $\mu$ to be ergodic.

With $|\cdot|$ the standard norm on $\mathbb{R}^k$, we assume the following hypotheses for the system $(M, \alpha, \mu)$.

**Standing hypotheses.** We assume there are

- a measurable, $\alpha$-invariant subset $\Lambda \subset M$ with $\mu(\Lambda) = 1$;
- an open set $U_0 \supset \Lambda$ equipped with a continuous Riemannian metric and an induced distance $d \colon U_0 \times U_0 \to [0, \infty)$;
- measurable functions $\rho, D \colon \Lambda \to [1, \infty)$ that are slowly increasing (over $\alpha$) on $\Lambda$ (see Definition 2.1);
- a measurable family of $C^1$ embeddings

$$\{\phi_x : x \in \Lambda\}, \quad \phi_x \colon \mathbb{R}^k(\rho(x)^{-1}) \to U_0$$

with the following properties:

(H1) $\phi_x \colon \mathbb{R}^k(\rho(x)^{-1}) \to U_0$ is a $C^1$ diffeomorphism onto its image with $\phi_x(0) = x$;

(H2) $\|D\phi_x\| \leq D(x)$ and $\|D\phi_x^{-1}\| \leq D(x)$; in particular $\phi_x \colon \mathbb{R}^k(\rho(x)^{-1}) \to (U_0, d)$ is a bi-Lipschitz embedding with $D(x)^{-1} \leq \mathrm{Lip}(\phi_x) \leq D(x)$.

Moreover, given any finite, symmetric subset $F \subset \mathbb{Z}^d$ that generates $\mathbb{Z}^d$, we assume there are

- an open subset $\Lambda \subset U \subset U_0$ such that for every $n \in F$ the restriction $\alpha(n){\restriction}_U \colon U \to U_0$ is a diffeomorphism of $U$ onto its range;
- measurable functions $r, C \colon \Lambda \to [1, \infty)$ that are slowly increasing (over $\alpha$) on $\Lambda$ with $\rho(x) \leq r(x)$;

such that

(H3) $\phi_x(\mathbb{R}^k(r(x)^{-1})) \subset U$ and for every $m \in F$

$$\alpha(m)(\phi_x(\mathbb{R}^k(r(x)^{-1})) \subset \phi_{\alpha(m,x)}(\mathbb{R}^k(\rho(\alpha(m,x)^{-1}))).$$

Moreover, for each $m \in F$, setting $f(\cdot) = \alpha(m, \cdot)$ and defining

$$\hat{f}_x := \phi_{f(x)}^{-1} \circ f \circ \phi_x \tag{11}$$

for $x \in \Lambda$, we assume

(H4) $\hat{f}_x \colon \mathbb{R}^k(r(x)^{-1}) \to \mathbb{R}^k(\rho(f(x))^{-1})$ is uniformly $C^{1+\beta}$ with

$$\|\hat{f}_x\|_{1+\beta} \leq C(x).$$

As $F$ is a generating set for $\mathbb{Z}^d$, it follows that given $x \in \Lambda$ and $n \in \mathbb{Z}^d$, the map

$$\hat{\alpha}(n, x) := \phi_{\alpha(n,x)}^{-1} \circ \alpha(n, x) \circ \phi_x$$

is defined and is a uniformly $C^{1+\beta}$ diffeomorphism on some neighborhood of 0. This induces a measurable cocycle $\mathcal{A}\colon \mathbb{Z}^d \times \Lambda \to \mathrm{GL}(d, \mathbb{R})$ given by the derivative

$$\mathcal{A}(n, x) = D_0 \hat{\alpha}(n, x). \tag{12}$$

We moreover assume that

(H5) $\left(x \mapsto \log^+ \|\mathcal{A}(n, x)\|\right) \in L^{d,1}(\mu)$ for every $n \in \mathbb{Z}^d$.

Applying Theorem 2.4, we write $\Lambda_0 \subset \Lambda$ in the remainder for the set of *regular points* of the cocycle $\mathcal{A}$ over the action of $\alpha$ on $(M, \mu)$.

**Remark 3.1.** Given local diffeomorphisms $g_1, g_2 \colon \mathbb{R}^k(1) \to \mathbb{R}^k$ preserving 0, defining $h = g_2 \circ g_1$ where defined we have

(1) $\|Dh\| \le \|Dg_2\| \|Dg_1\|$;
(2) $\mathrm{H\ddot{o}l}^\beta(Dh) \le \|Dh\| + \|Dg_2\| \, \mathrm{H\ddot{o}l}^\beta(Dg_1) + \|Dg_1\| \, \mathrm{H\ddot{o}l}^\beta(Dg_2) \|Dg_1\|^\beta$.

Suppose that a family of charts $\{\phi_x\}$, open set $U_0$, and functions $\rho$ and $D$ satisfying (H1)–(H2) exist such that for some finite symmetric generating set $F \subset \mathbb{Z}^d$ there are $U, r$ and $C$ such that (H3)–(H4) hold. Then, given any other finite symmetric generating set $F'$, we may modify the functions $r$ and $C$ and the open set $U$ so that (H3)–(H4) hold for $F'$.

Condition (H5) is independent of $F, U, r$, and $C$ hence remains valid passing to $F'$. Moreover, by the cocycle property (4), it is enough to verify $(x \mapsto \log^+ \|\mathcal{A}(n, x)\|) \in L^{d,1}(\mu)$ only for $n$ in some finite symmetric generating subset.

**Remark 3.2.** Given an invertible, measurable, measure-preserving transformation $f\colon (M, \mu) \to (M, \mu)$ we say the $f$ *satisfies the standing hypotheses* if the $\mathbb{Z}$-action generated by $f$ does. In this case, we take $F = \{1, -1\}$ so that (H3) and (H4) hold for both $f$ and $f^{-1}$.

**Remark 3.3.** As $L^{d,1}(\mu) \subset L^{d-1,1}(\mu)$ it follows that if the action $\alpha\colon \mathbb{Z}^d \times (M, \mu) \to (M, \mu)$ satisfies our standing hypotheses, then for any subgroup $H \subset \mathbb{Z}^d$, setting $\overline{\alpha}\colon H \times M \to M$ to be the restriction of $\alpha$ to $H$, we have that $(M, \overline{\alpha}, \mu)$ also satisfy our standing hypotheses. Moreover, for almost every $\alpha{\restriction}_H$-ergodic component $\mu_x^e$ of $\mu$, $(M, \alpha{\restriction}_H, \mu_x^e)$ satisfies the standing hypotheses.

## 4. Unstable manifolds and $C^{1+\beta}$-tame foliations

We introduce one of our primary dynamical objects of study. Recall $U_0 \subset M$ is equipped with a distance $d$. $B(x, r) \subset U_0$ denotes the open ball centered at $x$ of radius $r$.

### 4.1. $C^{1+\beta}$-tame foliations. Let $\mathcal{F}$ be a partition of $(M, \mu)$. We do not assume $\mathcal{F}$ to be measurable. Let $\mathcal{F}(x)$ denote the atom of $\mathcal{F}$ containing $x$.

**Definition 4.1.** A *measurable foliation* is a partition $\mathcal{F}$, a set $B = B(\mathcal{F}) \subset M$ with $\mu(B) = 0$, and a measurable function $r\colon M \to (0, \infty)$ such that

(1) for almost every $x \in M$, $\mathcal{F}(x) \smallsetminus B$ is a $C^1$ injectively immersed manifold in $M$ of constant dimension

and writing $\mathcal{F}(x, r(x))$ for the path connected component of $(\mathcal{F}(x) \smallsetminus B) \cap B(x, r(x))$ (relative to the immersed topology in $\mathcal{F}(x) \smallsetminus B$) containing $x$

(2) the family $\{\mathcal{F}(x, r(x))\}$ is a measurable family of $C^1$ embedded discs.

That is, by removing the set $B$ and sets of arbitrary small measure, given almost every $x$ there is neighborhood on which $\mathcal{F}$ restricts to a continuous lamination with uniformly $C^1$ leaves. We write $\mathcal{F}'(x) := \mathcal{F}(x) \smallsetminus B$ and refer to $\mathcal{F}'(x)$ as the *leaf of $\mathcal{F}$ through $x$*. Note $\mathcal{F}'(x)$ need not be connected.

Consider a $\mu$-preserving action $\alpha\colon \mathbb{Z}^d \times M \to M$ satisfying the hypotheses of Section 3.1.

**Definition 4.2.** A measurable foliation $\mathcal{F}$ is $\alpha$-*invariant* if, for every finite symmetric generating set $F \subset \mathbb{Z}^d$, the set $U \subset U_0$ in Section 3.1 and the null-set $B$ in Definition 4.1 can be taken so that for all $m \in F$,

$$\alpha(m)(\mathcal{F}'(x) \cap U) \subset \mathcal{F}'(\alpha(m, x)).$$

Recall that $\alpha(m)\restriction_U\colon U \to U_0$ is a diffeomorphism onto its image for $m \in F$. It follows by recurrence that the dimension of $\mathcal{F}'(x)$ is constant along orbits whence for $m \in F$, $\alpha(m)$ is a diffeomorphism between $\mathcal{F}'(x) \cap U$ and an open subset of $\mathcal{F}'(\alpha(n, x))$. In particular, if $\mu$ is $\alpha$-ergodic and $\mathcal{F}$ is $\alpha$-invariant then the leaves of $\mathcal{F}$ have constant dimension a.s.

Note that the geometry of the leaves of $\mathcal{F}$ as embedded in $M$ may degrade along orbits of $\alpha$ arbitrarily fast. We restrict ourselves to foliations for which this degradation is subexponential. We also impose additional regularity on the local geometry of leaves of the foliation.

Let $\alpha\colon \mathbb{Z}^d \times (M, \mu) \to (M, \mu)$ be an action satisfying the hypotheses of Section 3.1. Recall the family of charts $\phi_x$ introduced therein.

**Definition 4.3.** A measurable foliation $\mathcal{F}$ is $C^{1+\beta}$-*tame* (for the action $\alpha$ and relative to the charts $\phi_x$) if there is as set $\Lambda_{\mathcal{F}} \subset \Lambda$ with $\mu(\Lambda_{\mathcal{F}}) = 1$ such that for every $\varepsilon > 0$ there is a measurable function $\ell_{\mathcal{F}}\colon \Lambda_{\mathcal{F}} \to [1, \infty)$ that is $\varepsilon$-*slowly growing* (over $\alpha$) on $\Lambda_F$ and such that, writing $\hat{\mathcal{F}}(x)$ for the the path component (relative to the immersed topology) of

$$\phi_x^{-1}(\mathcal{F}'(x)) \cap \mathbb{R}^k(\ell_{\mathcal{F}}^{-1}(x))$$

containing 0, for $x \in \Lambda_{\mathcal{F}}$, $\hat{\mathcal{F}}(x)$ is the graph of a $C^{1+\beta}$ function

$$h_x^{\mathcal{F}}\colon U_x \subset T_x\hat{\mathcal{F}}(x) \to T_x\hat{\mathcal{F}}(x)^{\perp};$$

with

(1) $h_x^{\mathcal{F}}(0) = 0$ and $D_0 h_x^{\mathcal{F}} = 0$;
(2) $\mathrm{H\ddot{o}l}^{\beta}(Dh_x^{\mathcal{F}}) \le (\ell_{\mathcal{F}}(x))^{\beta}$ whence $\|Dh_x^{\mathcal{F}}\| \le 1$.

Note that the family of functions $\{h_x^{\mathcal{F}} : x \in \Lambda_{\mathcal{F}}\}$ depend measurably on $x$.

The primary examples of $\alpha$-invariant, $C^{1+\beta}$-tame, measurable foliations are the partitions into strong unstable manifolds and coarse and intermediate Lyapunov foliations arising in higher-rank, non-uniformly hyperbolic dynamics. (See Proposition 4.5.)

4.2. **Global unstable manifolds.** Let $\alpha\colon \mathbb{Z}^d \times (M, \mu) \to (M, \mu)$ be an action satisfying the hypotheses of Section 3.1. To simplify notation, we moreover assume the action $\alpha$ is ergodic with respect to $\mu$. All definitions here may be generalized to the non-ergodic case by considering ergodic components.

Let $\mathcal{L} = \{\lambda_i : 1 \le i \le p\}$ denote the Lyapunov exponent functionals of the derivative cocycle (12) with some fixed enumeration. Given $n \in \mathbb{Z}^d$ we choose a permutation $\sigma(n)$ of $\{1, 2, \dots, p\}$ and $u(n) \in \mathbb{N}$ so that

$$\lambda_{\sigma(n)(1)}(n) \ge \lambda_{\sigma(n)(2)}(n) \ge \cdots \ge \lambda_{\sigma(n)(u(n))}(n) > 0 \ge \cdots \ge \lambda_{\sigma(n)(p)}(n).$$

Note that by choosing $n$ in general position we can ensure all inequalities above are strict if necessary.

**Definition 4.4.** Given $x \in M$, $n \in \mathbb{Z}^d$, and $1 \leq i \leq u(n)$ define the *$i$th unstable manifold through $x$ for $\alpha(n)$* to be the set

$$W_n^{u,i}(x) = \left\{ y \in M \mid \limsup_{k \to -\infty} \frac{1}{k} \log d(\alpha(kn, x), \alpha(kn, y)) \leq -\lambda_{\sigma(n)(i)}(n) \right\}.$$

The *unstable manifold through $x$ for $\alpha(n)$* is the set

$$W_n^u(x) := \left\{ y \in M \mid \limsup_{k \to -\infty} \frac{1}{k} \log d(\alpha(kn, x), \alpha(kn, y)) < 0 \right\}.$$

It follows from Proposition 4.5 below that $W_n^u(x) = W_n^{u,u(n)}(x)$ for almost every $x$. Note that if $n$ is not in general position so that $\lambda_{\sigma(n)(i)}(n) = \lambda_{\sigma(n)(i+1)}(n)$ for some $i$ with $\lambda_{\sigma(n)(i+1)}(n) > 0$ then we declare

$$W_n^{u,i}(x) = W_n^{u,i+1}(x).$$

Implicit in the above definition is that $d(\alpha(kn, x), \alpha(kn, y))$ is defined for all but finitely many $k < 0$; in particular, we have $\alpha(kn, x) \in U_0$ and $\alpha(kn, y) \in U_0$ for all but finitely many $k \leq 0$. If $\alpha(kn, x) \in M \smallsetminus U_0$ for infinitely many $k < 0$, declare $W_n^{u,i}(x) = \{x\}$. From the above definition, it follows for each $n \in \mathbb{Z}^d$ that the collection of $i$th unstable manifolds forms a partition $\mathscr{W}_n^{u,i}$ of $M$. We also write $\mathscr{W}_n^u$ of $M$ for the partition into unstable manifolds for $\alpha(n)$. It follows form Proposition 4.5 below that $\mathscr{W}_n^u$ and $\mathscr{W}_n^{u,u(n)}$ coincide off a null set.

As we do not assume the action $\alpha$ to be by diffeomorphisms, $W_n^{u,i}(x)$ is in general not a manifold. However, from Proposition 4.5 below, under our standing hypotheses, the partition $\mathscr{W}_n^{u,i}$ has the structure of a $C^{1+\beta}$-tame, measurable foliation tangent a.e. to $D_0\phi_x(\bigoplus_{1 \leq j \leq i} E_{\lambda_{\sigma(n)(j)}}(x))$. Then with the notation introduced above, $\left(\mathscr{W}_n^{u,i}\right)'(x)$ is an injectively immersed manifold of dimension $\dim\left(\bigoplus_{1 \leq j \leq i} E_{\lambda_{\sigma(n)(j)}}(x)\right)$.

Given an $\alpha$-invariant, measurable foliation $\mathcal{F}$ and $x \in M$, let $d_\mathcal{F}$ denote the distance on $\mathcal{F}'(x) \cap U_0$ induced by restriction of the Riemannian structure on $U_0$. In particular, $d_\mathcal{F}(x, y)$ is defined if and only if $y \in \mathcal{F}'(x)$ and there is a $C^1$ path in $\mathcal{F}'(x) \cap U_0$ from $x$ to $y$. Given a measurable foliation $\mathcal{F}$ and $n \in \mathbb{Z}^d$ define

$$\left(\mathcal{F} \vee \mathscr{W}_n^{u,i}\right)(x) := \left\{ y \in M \mid \limsup_{k \to -\infty} \frac{1}{k} \log d_\mathcal{F}(\alpha(kn, x), \alpha(kn, y)) \leq -\lambda_{\sigma(n)(i)}(n) \right\}.$$

This defines a partition $\mathcal{F} \vee \mathscr{W}_n^{u,i}$ of $(M, \mu)$. Similarly define $\mathcal{F} \vee \mathscr{W}_n^u$.

The following proposition will be shown in the next section.

**Proposition 4.5.** *Given any $n \in \mathbb{Z}^d$*

(a) *$\mathscr{W}_n^{u,i}$ is a $C^{1+\beta}$-tame, $\alpha$-invariant, measurable foliation with*

$$T_x\left(\mathscr{W}_n^{u,i}\right)'(x) = D_0\phi_x\left(\bigoplus_{1 \leq j \leq i} E_{\lambda_{\sigma(n)(j)}}(x)\right)$$

*for almost every $x$;*

(b) *given any $C^{1+\beta}$-tame, $\alpha$-invariant, measurable foliation $\mathcal{F}$ the partition $\mathcal{F} \vee \mathscr{W}_n^i$ is a $C^{1+\beta}$-tame, $\alpha$-invariant, measurable foliation with*

$$T_x\left(\left(\mathcal{F} \vee \mathscr{W}_n^{u,i}\right)'(x)\right) = T_x\left(\mathscr{W}_n^{u,i}\right)'(x) \cap T_x\mathcal{F}'(x)$$

*for almost every $x$.*

**Remark 4.6.** Note that the term *measurable foliation* indicates that the transverse structure of the foliation is measurable. In general, a measurable foliation is *not* a measurable partition. In particular, the partition of $(M, \mu)$ into global unstable manifolds $\mathscr{W}_n^u$ is never a measurable partition if $\alpha(n)$ has positive entropy.

As a corollary of Proposition 4.5 we have the following uniqueness property.

**Lemma 4.7.** *Let* $n_1, n_2 \in \mathbb{Z}^d$ *have the following property: for some* $1 \le i \le \min\{u(n_1), u(n_2)\}$

    (1) $\{\sigma(n_1)(j) : 1 \le j \le i\} = \{\sigma(n_2)(j) : 1 \le j \le i\}$,

    (2) $\lambda_{\sigma(n_1)(i)}(n_1) > \lambda_{\sigma(n_1)(i+1)}(n_1)$, *and*

    (3) $\lambda_{\sigma(n_2)(i)}(n_2) > \lambda_{\sigma(n_2)(i+1)}(n_2)$.

*Then* $\mathscr{W}_{n_1}^{u,i} \overset{\circ}{=} \mathscr{W}_{n_2}^{u,i}$. *In particular, if* $\operatorname{sgn}(\lambda_i(n_1)) = \operatorname{sgn}(\lambda_i(n_2))$ *for all* $\lambda_i \in \mathcal{L}$ *then* $\mathscr{W}_{n_1}^u \overset{\circ}{=} \mathscr{W}_{n_2}^u$.

Indeed, the lemma follows from Proposition 4.5 as for almost every $x$

$$T_x \left( \mathscr{W}_{n_1}^{u,i} \right)'(x) = T_x \left( \mathscr{W}_{n_2}^{u,i} \right)'(x)$$

whence $\mathscr{W}_{n_1}^{u,i} \vee \mathscr{W}_{n_2}^{u,i} \overset{\circ}{=} \mathscr{W}_{n_1}^{u,i}$.

4.3. **Coarse Lyapunov exponents, manifolds, and foliations.** We continue to assume $\alpha \colon \mathbb{Z}^d \times (M, \mu) \to (M, \mu)$ satisfies our standing hypotheses with $\mu$ ergodic. Let $\mathcal{L} = \{\lambda_i\}$ denote the Lyapunov exponent functionals of the derivative cocycle (12).

**Definition 4.8.** Two Lyapunov exponents $\lambda_i$ and $\lambda_j \in \mathcal{L}$ are *coarsely equivalent* if there is a $c > 0$ such that $\lambda_i = c\lambda_j$. A *coarse Lyapunov exponent* is an equivalence class of $\mathcal{L}$ under the coarse equivalence relation

We write $\hat{\mathcal{L}}$ for the set of coarse Lyapunov exponents. Given $\chi \in \hat{\mathcal{L}}$, the sign of $\chi(n)$ is well-defined.

**Definition 4.9.** Given $\chi \in \hat{\mathcal{L}}$ with $\chi \neq 0$, the *coarse Lyapunov foliation* corresponding to $\chi$ is

$$\mathscr{W}^\chi := \bigvee_{\{n \in \mathbb{Z}^d : \chi(n) > 0\}} \mathscr{W}_n^u. \tag{13}$$

The *coarse Lyapunov manifold* corresponding to $\chi$ through $x$ is the corresponding leaf $W^\chi(x) := (\mathscr{W}^\chi)'(x)$.

From Lemma 4.7 the intersection in (13) above is equivalent to an intersection taken over a finite subset of $\mathbb{Z}^d$. It then follows from Proposition 4.5 that $\mathscr{W}^\chi$ is a $C^{1+\beta}$-tame, $\alpha$-invariant, measurable foliation.

## 5. LYAPUNOV CHARTS; PROPERTIES OF TAME AND UNSTABLE FOLIATIONS

Let $(M, \mu)$ and $\alpha \colon \mathbb{Z}^d \times M \to M$ satisfy the hypotheses of Section 3.1. To simplify notation, we assume that $\alpha$ acts ergodically on $(M, \mu)$. We present here a standard construction which, via a local change of coordinates, converts the local non-uniformly partially hyperbolic dynamics into uniformly partially hyperbolic dynamics.

Recall the Lyapunov exponent functionals $\mathcal{L} = \{\lambda_i\}$. For the remainder of this section, fix $F \subset \mathbb{Z}^d$ to be a finite, symmetric generating set as in Section 3.1. We recall the sets $\Lambda$, $U_0$ and $U$, functions $r, \rho$ and $C$, and all other notation from Section 3.1. All constructions below are relative to this choice of $F$ and corresponding $U, r, \rho$ and $C$. Write

$$\lambda_0 = \max\{|\lambda_i(n)| : n \in F\}$$

and

$$\varepsilon_0 := \min\{1, |\lambda_i(n)|, |\lambda_i(n) - \lambda_j(n)| : n \in F, \lambda_i(n) \neq \lambda_j(n), \lambda_i(n) \neq 0\}/100.$$

Note if $\lambda_i(n) > 0$ then for any $0 \leq \varepsilon \leq \varepsilon_0$, then

$$e^{\lambda_i(n) - 2\varepsilon} \leq e^{\lambda_i(n) - \varepsilon} - \varepsilon \leq e^{\lambda_i(n) + \varepsilon} + \varepsilon \leq e^{\lambda_i(n) + 2\varepsilon};$$

if $\lambda_i(n) = 0$ then $e^{\lambda_i(n) + \varepsilon} + \varepsilon \leq e^{\lambda_i(n) + 2\varepsilon}$. Recall that $\Lambda_0 \subset \Lambda$ denotes the set of regular points in Theorem 2.4.

## 5.1. Lyapunov charts.

We specify an alternative norm on $\mathbb{R}^k$. Fix an orthogonal (with respect to the standard inner product) decomposition $\mathbb{R}^k = \oplus \mathbb{R}^i$ where $\dim \mathbb{R}^i = m_i$ is the dimension of $E_{\lambda_i}(x)$ for $x \in \Lambda_0$. Define the norm $\|\cdot\|$ on $\mathbb{R}^k$ as follows: writing $v = \sum v_i$ for $v_i \in \mathbb{R}^i$ set $\|v\| = \max\{|v_i|\}$ where $|v_i|$ restricts to the norm induced by the standard inner product on each $\mathbb{R}^i$.

**Proposition 5.1** (Lyapunov charts). *For every $0 < \varepsilon < \varepsilon_0$ there is an $\varepsilon$-slowly growing function $\ell : \Lambda_0 \to [1, \infty)$ and a measurable family of invertible linear maps $\{L_x : (\mathbb{R}^k, \|\cdot\|) \to (\mathbb{R}^k, |\cdot|) : x \in \Lambda\}$ such that, defining a measurable family of $C^1$ embeddings*

$$\Phi_x = \phi_x \circ L_x : \mathbb{R}^k(1, \|\cdot\|) \to M,$$

*we have*

(a) $\Phi_x(0) = x$;

(b) $L_x \mathbb{R}^i = E_{\lambda_i}(x)$;

(c) $L_x : \mathbb{R}^k(e^{-\lambda_0 - 2\varepsilon}, \|\cdot\|) \subset \mathbb{R}^k(r(x)^{-1}, |\cdot|)$ *and* $L_x : \mathbb{R}^k(1, \|\cdot\|) \subset \mathbb{R}^k(\rho(x)^{-1}, |\cdot|)$.

*Furthermore, for every $n \in F$, writing $f = \alpha(n)$ we have that the map*

(d) $\tilde{f}_x : \mathbb{R}^k(e^{-\lambda_0 - 2\varepsilon}, \|\cdot\|) \to \mathbb{R}^k(1, \|\cdot\|)$ *given by*

$$\tilde{f}_x(v) := \Phi_{f(x)}^{-1} \circ f \circ \Phi_x(v) = L_{f(x)}^{-1} \circ \hat{f}_x \circ L_x \tag{14}$$

*is well-defined (where $\hat{f}_x$ is as in (11));*

(e) $D_0 \tilde{f}_x \mathbb{R}^i = \mathbb{R}^i$ *for every $i$ and for $v \in \mathbb{R}^i$*

$$e^{\lambda_i(n) - \varepsilon} \|v\| \leq \|D_0 \tilde{f}_x v\| \leq e^{\lambda_i(n) + \varepsilon} \|v\|;$$

(f) $\mathrm{H\ddot{o}l}^\beta(D\tilde{f}_x) \leq \varepsilon$ *hence* $\mathrm{Lip}(\tilde{f}_x - D_0 \tilde{f}_x) \leq \varepsilon$;

(g) $\ell(x)^{-1} \leq \|L_x\| \leq 1$ *and* $\ell(x)^{-1} \leq \mathrm{Lip}(\Phi_x) \leq 1$.

We call a family of embeddings $\{\Phi_x\}$ satisfying the above properties a family of $\varepsilon$-*Lyapunov charts*. As our construction is somewhat different than others in the literature we include a proof.

*Proof.* Take $0 < \varepsilon' < \varepsilon/4$ sufficiently small. Recall the $\varepsilon'$-Lyapunov metric $\langle\!\langle\!\langle \cdot, \cdot \rangle\!\rangle\!\rangle_{x, \varepsilon'}$ defined by (7) and the corresponding family of norms $\|\!|\cdot|\!\|_{x, \varepsilon'}$. Relative to the inner products $\langle\!\langle\!\langle \cdot, \cdot \rangle\!\rangle\!\rangle_{x, \varepsilon'}$ on $\mathbb{R}^k$, choose a measurable orthonormal basis for each $E_{\lambda_i}(x)$ which, in turn, defines a measurable family of linear isometries $\hat{\tau}_x : (\mathbb{R}^k, \|\!|\cdot|\!\|_{x, \varepsilon'}) \to (\mathbb{R}^k, \|\cdot\|)$ with $\hat{\tau}_x E_{\lambda_i}(x) = \mathbb{R}^i$ for every $i$.

Let $k_0$ be the constant and $L(x)$ the $\varepsilon'$-slowly increasing function in Lemma 2.8. Let $\hat{C}(x)$ be an $\varepsilon'$-slowly increasing function such that the functions $D(x), C(x), r(x), \rho(x)$ appearing in our standing hypotheses of Section 3.1 are bounded above by $\hat{C}(x)$ for all $x \in \Lambda_0$. Take $\tau_x : (\mathbb{R}^k, \|\cdot\|) \to (\mathbb{R}^k, |\cdot|)$ defined by $\tau_x(v) = \hat{\tau}_x^{-1}(v)$. We have

$$\|\tau_x\| \leq k_0, \qquad \|\tau_x^{-1}\| \leq L(x).$$

Take
$$\hat{\ell}(x) = \left( (\varepsilon/2)^{-1} k_0^{1+\beta} L(x) \hat{C}(x) \right)^{\frac{1}{\beta}}$$
and take $L_x \colon (\mathbb{R}^k, \|\cdot\|) \to (\mathbb{R}^k, |\cdot|)$ to be the linear map defined by
$$L_x(v) = \tau_x(\hat{\ell}^{-1}(x)v).$$
Having taken $\varepsilon'$ sufficiently small, $\hat{\ell}$ is $(\varepsilon/2)$-slowly increasing and we verify properties (a), (b), and (c). Also, the upper bounds in (g) follow.

Given $n \in F$ write $f = \alpha(n)$. Given $x \in \Lambda_0$, write $\hat{f}_x$ as in (11). Consider first $\overline{f}_x \colon \mathbb{R}^k(r(x) k_0^{-1}, \|\cdot\|) \to (\mathbb{R}^k, \|\cdot\|)$ given by $\overline{f}_x = \tau_{f(x)}^{-1} \circ \hat{f}_x \circ \tau_x$. We have
$$D_v \overline{f}_x(\xi) = \tau_{f(x)}^{-1} \left( D_{\tau_x(v)} \hat{f}_x \left( \tau_x(\xi) \right) \right)$$
hence
$$\left| D_v \overline{f}_x(\xi) - D_u \overline{f}_x(\xi) \right| \le k_0 L(x) \operatorname{Höl}^\beta(Df_x) \left| \tau_x(v) - \tau_x(u) \right|^\beta$$
$$\le k_0 L(x) \hat{C}(x) k_0^\beta \left| u - v \right|^\beta$$
and
$$\operatorname{Höl}^\beta(D\overline{f}_x) \le k_0 L(x) \hat{C}(x) k_0^\beta \le \hat{\ell}(x)^\beta \varepsilon/2.$$
Moreover for $v \in \mathbb{R}^i$ we have
$$e^{\lambda_i(n) - \varepsilon'} \|v\| \le \|D_0 \overline{f}_x v\| \le e^{\lambda_i(n) + \varepsilon'} \|v\|;$$
By (c), $\tilde{f}_x \colon \mathbb{R}^k(e^{-\lambda_0 - 2\varepsilon}, \|\cdot\|) \to \mathbb{R}^k$ defined by
$$\tilde{f}_x(v) = \hat{\ell}(f(x)) \overline{f}_x(\hat{\ell}(x)^{-1} v)$$
is well-defined. We have
$$D_v \tilde{f}_x(\xi) = \hat{\ell}(f(x)) \hat{\ell}(x)^{-1} \left( D_{\hat{\ell}(x)^{-1} v} \overline{f}_x(\xi) \right)$$
hence
$$\left| D_v \tilde{f}_x(\xi) - D_u \tilde{f}_x(\xi) \right| \le \hat{\ell}(f(x)) \hat{\ell}(x)^{-1} \operatorname{Höl}^\beta(D\overline{f}_x) \left| \hat{\ell}(x)^{-1} v - \hat{\ell}(x)^{-1} u \right|^\beta$$
$$\le \hat{\ell}(f(x)) \hat{\ell}(x)^{-1} \varepsilon/2 \left| u - v \right|^\beta \le \varepsilon \left| u - v \right|^\beta$$
whence (f) follows. Also,
$$e^{-\varepsilon/2} \|D_0 \overline{f}_x(v)\| \le \|D_0 \tilde{f}_x(v)\| = \hat{\ell}(f(x)) \hat{\ell}(x)^{-1} \|D_0 \overline{f}_x(v)\| \le e^{\varepsilon/2} \|D_0 \overline{f}_x(v)\|$$
and (e) follows.

Finally, with $\ell(x) = k_0 L(x) \hat{C}(x) \hat{\ell}(x)$, we have that $\ell(x)$ is $\varepsilon$-slowly growing and satisfies the lower bounds in (g). Also, as $e^{\lambda_0 + \varepsilon} + \varepsilon \le e^{\lambda_0 + 2\varepsilon}$, the Lipschitz constant of $\tilde{f}_x$ is bounded above by $e^{\lambda_0 + 2\varepsilon}$ and (d) follows. $\qquad\square$

## 5.2. Lyapunov charts adapted to $C^{1+\beta}$-tame foliations.

Let $\mathcal{F}$ be an $\alpha$-invariant $C^{1+\beta}$-tame measurable foliation. The following proposition guarantees that Lyapunov charts above may be chosen so that, relative to the charts $\Phi_x$, the leaves $\mathcal{F}(x)$ are uniformly $C^{1+\beta}$ embedded.

**Proposition 5.2.** *An $\alpha$-invariant measurable foliation $\mathcal{F}$ of $(M, \mu)$ is $C^{1+\beta}$-tame if and only if for every $0 < \varepsilon < \varepsilon_0$, the charts $\Phi_x$ in Proposition 5.1 can be chosen so that, in addition to the properties in Proposition 5.1, there are*

- *a set $\Lambda' \subset \Lambda_0$ with $\mu(\Lambda') = 1$;*

- a subspace $V \subset \mathbb{R}^k$ with orthogonal complement[2] $W$;
- a measurable family of $C^{1+\beta}$ functions

$$\tilde{h}_x^{\mathcal{F}} : V(e^{-\varepsilon}, \|\cdot\|) \to W(1, \|\cdot\|)$$

  defined for $x \in \Lambda'$

such that, writing $\tilde{\mathcal{F}}_x$ for the path component (relative to the immersed topology) of $\Phi_x^{-1}(\mathcal{F}(x)') \cap \mathbb{R}^k(1)$ containing 0, $\tilde{\mathcal{F}}_x$ contains as an open set the graph of $\tilde{h}_x^{\mathcal{F}} : V(e^{-\varepsilon}, \|\cdot\|) \to W(1, \|\cdot\|)$ and

(i) $\tilde{h}_x^{\mathcal{F}}(0) = 0$; $D_0 \tilde{h}_x^{\mathcal{F}} = 0$;
(j) $\mathrm{H\ddot{o}l}^\beta(D\tilde{h}_x^{\mathcal{F}}) \le \epsilon$ and hence $\|D\tilde{h}_x^{\mathcal{F}}\| \le \epsilon$;
(k) for all $n \in F$, writing $\tilde{f}_x$ as in (14) we have $\tilde{f}_x$ is a diffeomorphism between $\mathrm{graph}(\tilde{h}_x^{\mathcal{F}}) \cap \mathbb{R}^n(e^{-\lambda_0 - 2\varepsilon})$ and an open subset of $\mathrm{graph}(\tilde{h}_{\alpha(n)(x)}^{\mathcal{F}})$.

Note that for each $n \in F$ writing $f = \alpha(n)$ and $\tilde{f}_x$ as in (14), $V$ is $D_0 \tilde{f}_x$-invariant.

The proof of Proposition 5.2 relies on the following lemma which, in turn, follows from the Implicit Function Theorem with Hölder estimates [Pes, Lemma 2.1.1]. Consider $\mathbb{R}^k$ equipped with 2 norms $\|\cdot\|_1$ and $\|\cdot\|_2$. Let $V^1, V^2 \subset \mathbb{R}^k$ be a subspaces and for $j \in \{1, 2\}$ let $W^i$ be a subspace of complementary dimension transverse to $V^j$. We assume the decompositions $V^j \oplus W^j$ and norms $\|\cdot\|_j$ have the property that, for any vector $u \in \mathbb{R}^k$, writing $u = v + w$ for $v \in V^j$ and $w \in W^j$ we have

$$\|u\|_j \ge \max\{\|v\|_j, \|w\|_j\}. \tag{15}$$

Our applications below of the above setup are: $\|\cdot\|_j$ is the Euclidean norm and $W^i$ is the orthogonal complement of $V^j$; or $\|\cdot\|_j$ is the norm on $\mathbb{R}^k$ specified in Section 5.1, $V^j$ satisfies $V^j = \oplus_i (V \cap \mathbb{R}^i)$, and $W^2$ is the orthogonal complement of $V^2$.

With the above setup we have the following.

**Lemma 5.3.** *Let* $h\colon V^1(r, \|\cdot\|_1) \to W^1$ *be a* $C^{1+\beta}$ *function such that (relative to* $\|\cdot\|_1$*) we have*

- $D_0 h = 0$
- $\mathrm{H\ddot{o}l}^\beta(Dh) \le a$
- $\|Dh\| \le 1$.

*Let* $L\colon \mathbb{R}^n \to \mathbb{R}^n$ *be an invertible linear map with* $L(V^1) = V^2$. *Take*

- $a_0 = a\|L^{-1}\|^{1+\beta}$
- $b_0 = 2\|L^{-1}\|$
- $c_0 = \|L\|$.

*Then with*

$$r_0 = \min\left\{ \frac{r/2}{\|L^{-1}\|}, \frac{r/2}{2b_0 c_0 \|L^{-1}\|}, \frac{1}{(1 + 2b_0 c_0)(2a_0 c_0)^{1/\beta}}, \frac{1}{((1 + b_0 c_0)^2 (8a_0 c_0))^{1/\beta}} \right\}$$

*we have that* $L(\mathrm{graph}(h))$ *contains as an open set the graph of a* $C^{1+\beta}$*-function* $\hat{h}\colon V^2(r_0) \to W^2$ *such that (relative to* $\|\cdot\|_2$*) we have*

(1) $D_0 \hat{h} = 0$;
(2) $\mathrm{H\ddot{o}l}^\beta(D\hat{h}) \le 8 a_0 c_0 (1 + b_0 c_0)^2$;

*whence*

(3) $\|D\hat{h}\| \le 1$.

---
[2]relative to the standard inner product on $\mathbb{R}^k$

*Proof.* We have that $\mathrm{graph}(h)$ in $(\mathbb{R}^k, \|\cdot\|_1)$ is the solution set to $\psi(v, w_1) \equiv 0$ where

$$\psi(v_1, w_1) = w_1 - h(v_1)$$

is defined for all $(v_1, w_1)$ with $\|v_1\|_1 < r$. Then $L(\mathrm{graph}(h))$ contains as an open set the solution set $\overline{\psi}(v_2, w_2) \equiv 0$ of the function

$$\overline{\psi} \colon V^2\left(\frac{r}{2\|L^{-1}\|}, \|\cdot\|_2\right) \times W^2\left(\frac{r}{2\|L^{-1}\|}, \|\cdot\|_2\right) \to W^1$$

given by

$$\overline{\psi}(v_2, w_2) = \psi \circ L^{-1}(v_2, w_2).$$

Decomposing the domains $\mathbb{R}^k = V^2 \oplus W^2$ and $\mathbb{R}^k = V^1 \oplus W^1$ of $L^{-1}$ and $\psi$ we have

$$D_{(v,w_2)}\overline{\psi}(\xi) = \left[\begin{array}{c|c} -D_{L^{-1}v}h & I \end{array}\right] L^{-1}(\xi).$$

We check the following.

(1) The partial derivative $D_{2,(0,0)}\overline{\psi} \colon W^2 \to W^1$ is of the form

$$D_{2,(0,0)}\overline{\psi}(\xi) = \left[\begin{array}{c|c} 0 & I \end{array}\right] L^{-1}(\xi)$$

and is a linear isomorphism with

$$\left\|\left(D_{2,(0,0)}\overline{\psi}\right)^{-1}\right\| \leq \|L\|.$$

Indeed, given $w_2 \in W^2$ write $L^{-1}(w_2) = v_1 + w_1$. Then $D_{2,(0,0)}\overline{\psi}(w_2) = w_1$, and writing $L(w_1) = w_2 - L(v_1)$ we have from (15) that

$$\|w_2\|_2 \leq \|Lw_1\|_2 \leq \|L\|\|w_1\|_1.$$

(2) The partial derivative $D_{1,(v,0)}\overline{\psi} \colon V^2 \to W^1$ is of the form

$$D_{1,(v,0)}\overline{\psi}(\xi) = \left[\begin{array}{c|c} -D_{L^{-1}v}h & I \end{array}\right] L^{-1}(\xi)$$

hence

$$\max \|D_{1,(v,0)}\overline{\psi}\| \leq 2\|L^{-1}\|.$$

(3) $D_{(v,w_2)}\overline{\psi}(\xi) - D_{(\hat{v},\hat{w}_2)}\overline{\psi}(\xi) = \left[\begin{array}{c|c} -D_{L^{-1}v}h + D_{L^{-1}\hat{v}}h & 0 \end{array}\right] L^{-1}(\xi)$, so

$$\mathrm{H\ddot{o}l}^\beta(D\overline{\psi}) \leq \|L^{-1}\|^{1+\beta}a.$$

The conclusion of the proposition then follows from [Pes, Lemma 2.1.1] and the fact that $L(\mathrm{graph}(h))$ is tangent to $V^2$. □

*Proof of Proposition 5.2.* The sufficiency follows from Lemma 5.3 with $L = L_x \colon (\mathbb{R}^k, \|\cdot\|) \to (\mathbb{R}^k, |\cdot|)$ the maps from Proposition 5.1.

For the necessity, let $\mathcal{F}$ be an $\alpha$-invariant $C^{1+\beta}$-tame, measurable foliation with $C^{1+\beta}$ leaves. We retain all notation from Definitions 4.2 and 4.3. In particular, we assume for $n \in F$ that $\alpha(n)(\mathcal{F}'(x) \cap U) \subset \mathcal{F}'(\alpha(n, x))$ for almost all $x$. Take $\Lambda' = \Lambda_0 \cap \Lambda_{\mathcal{F}}$. Note that the $\alpha$-invariance of $\mathcal{F}$ implies that $T_0\hat{\mathcal{F}}(x) = \oplus(E_{\lambda_i} \cap T_0\hat{\mathcal{F}}(x))$ for almost every $x$ (where $\hat{\mathcal{F}}(x)$ is as in Definition 4.3.) In the proof of Proposition 5.1 we may select $V \subset \mathbb{R}^k$ and construct the maps $\tau_x \colon (\mathbb{R}^k, \|\cdot\|) \to (\mathbb{R}^k, |\cdot|)$ so that $\tau_x(V) = T_0\hat{\mathcal{F}}(x)$ for all $x \in \Lambda'$.

Taking $\varepsilon' > 0$ sufficiently small in the proof of Proposition 5.1 and an $\varepsilon'$-slowly growing function $\ell_{\mathcal{F}}$ as in Definition 4.3, applying Lemma 5.3 to the maps $\tau_x^{-1} \colon (\mathbb{R}^k, |\cdot|) \to (\mathbb{R}^k, \|\cdot\|)$ appearing the proof of Proposition 5.1 we may find an $(\varepsilon/2)$-slowly increasing function $\overline{\ell}(x)$ so that

$$\tau_x^{-1}(\hat{\mathcal{F}}(x)) \cap \mathbb{R}^k(\overline{\ell}(x)^{-1}, \|\cdot\|)$$

is the graph of a $C^{1+\beta}$ function

$$\hat{h}_x \colon V(\overline{\ell}(x)^{-1}, \|\cdot\|) \to V^\perp(\overline{\ell}(x)^{-1}, \|\cdot\|)$$

with

(1) $\hat{h}_x(0) = 0$ and $D_0\hat{h}_x(0) = 0$;
(2) $\operatorname{H\ddot{o}l}^\beta(D\hat{h}_x) \le \varepsilon/2\overline{\ell}(x)^\beta$.

Then taking

$$\hat{\ell}(x) = \max\left\{\overline{\ell}(x), \left((\varepsilon/2)^{-1} k_0^{1+\beta} L(x)^p \hat{C}(x)^2\right)^{\frac{1}{\beta}}\right\}$$

in the proof of Proposition 5.1, the results of Proposition 5.1 remain valid.

Moreover from similar computations as in the proof of Proposition 5.1

$$L_x^{-1}(\hat{\mathcal{F}}(x)) \cap \mathbb{R}^k(1, \|\cdot\|)$$

is the graph of a $C^{1+\beta}$ function

$$\tilde{h}_x \colon V(1, \|\cdot\|) \to V^\perp(1, \|\cdot\|)$$

with

(1) $\tilde{h}_x(0) = 0$ and $D_0\tilde{h}_x(0) = 0$;
(2) $\operatorname{H\ddot{o}l}^\beta(D\tilde{h}_x) \le \varepsilon$ whence $\|D\tilde{h}_x\| \le \varepsilon$.         $\square$

The invariance in (k) follows from construction.

5.2.1. *Lyapunov charts for dynamics restricted to leaves.* Note that if $\mathcal{F}$ is an $\alpha$-invariant, $C^{1+\beta}$-tame, measurable foliation then with $V$ as in the notation of Proposition 5.2 we have

$$V = \bigoplus(V \cap \mathbb{R}^i).$$

As in Propositions 5.2 and 5.1, given an $\alpha$-invariant, $C^{1+\beta}$-tame, measurable foliation $\mathcal{F}$ we construct a family of charts for the restriction of the dynamics to leaves of $\mathcal{F}$. Given $0 < \varepsilon' < \varepsilon_0$ sufficiently small, let $\{\Phi_x\}$ a family of $\varepsilon'$-Lyapunov charts satisfying Proposition 5.2. Take $\Lambda'$ and $\tilde{h}_x^{\mathcal{F}}(v)$ as in Proposition 5.2. Given $x \in \Lambda'$ define

$$\tilde{H}_x^{\mathcal{F}} \colon \mathbb{R}^k(1, \|\cdot\|) \to \mathbb{R}^k(\|\cdot\|)$$

relative to the decomposition $\mathbb{R}^k = V \oplus W$ by

$$\tilde{H}_x^{\mathcal{F}}(v, w) = (v, w + \tilde{h}_x^{\mathcal{F}}(v)).$$

We then define a measurable family of embeddings $\Psi_x^{\mathcal{F}} \colon V(1, \|\cdot\|) \to \mathcal{F}(x)$ by

$$\Psi_x^{\mathcal{F}}(v) = \Phi_x \circ \tilde{H}_x^{\mathcal{F}}(v).$$

**Proposition 5.4.** *For every $0 < \varepsilon < \varepsilon_0$ there is a family of $\varepsilon$-Lyapunov charts $\{\Phi_x : x \in \Lambda_0\}$, a $\varepsilon$-slowly growing function $\overline{\ell} \colon \Lambda' \to [1, \infty)$ and a measurable family of $C^1$ embeddings $\{\Psi_x^{\mathcal{F}} : x \in \Lambda'\}$ defined as above with*

(a) $\Psi_x^{\mathcal{F}}(0) = x$ and $\Psi_x^{\mathcal{F}}(V(e^{-\varepsilon}, \|\cdot\|)) \subset \mathcal{F}'(x)$ for almost every $x$.

*Furthermore, for every $n \in F$, writing $f = \alpha(n)$ we have*

(b) *the function $\overline{f}_x \colon V(e^{-\lambda_0 - 2\varepsilon}, \|\cdot\|) \to V(e^{-\varepsilon}, \|\cdot\|)$ given by*

$$\overline{f}_x(v) := (\Psi_{f(x)}^{\mathcal{F}})^{-1} \circ f \circ \Psi_x^{\mathcal{F}}(v)$$

*is well-defined;*

(c) $D_0\overline{f}_x(\mathbb{R}^i \cap V) = (\mathbb{R}^i \cap V)$ *for every $i$ and for $v \in (\mathbb{R}^i \cap V)$*

$$e^{\lambda_i(n) - \varepsilon}\|v\| \le \|D_0\overline{f}_x v\| \le e^{\lambda_i(n) + \varepsilon}\|v\|;$$

  (d) $\mathrm{Höl}^\beta(D\overline{f}_x) \le \varepsilon$ hence $\mathrm{Lip}(\overline{f}_x - D_0\overline{f}_x) \le \varepsilon$;
  (e) $\overline{\ell}(x)^{-1} \le \mathrm{Lip}(\Psi_x^{\mathcal{F}}) \le 1$.

*Proof.* For $0 < \varepsilon' < \varepsilon_0$ sufficiently small, let $\{\Phi_x\}$ a family of $\varepsilon'$-Lyapunov charts satisfying Proposition 5.2. Then with $\tilde{H}_x^{\mathcal{F}}$ as above

  (1) $\|D\tilde{H}_x^{\mathcal{F}}\| \le 1 + \varepsilon'$;
  (2) $\|D(\tilde{H}_x^{\mathcal{F}})^{-1}\| \le 1 + \varepsilon'$;
  (3) $\mathrm{Höl}^\beta(D\tilde{H}_x^{\mathcal{F}}) \le \varepsilon'$.

Given $n \in F$, write $f = \alpha(n)$ and define $\tilde{f}_x$ as in (14). Define $\overline{F}_x \colon \mathbb{R}^k\left(\frac{1}{(1+\varepsilon')^2}e^{-\lambda_0 - 2\varepsilon'}, \|\cdot\|\right) \to \mathbb{R}^k(1, \|\cdot\|)$ by

$$\overline{F}_x = (\tilde{H}_{f(x)}^{\mathcal{F}})^{-1} \circ \tilde{f}_x \circ \tilde{H}_x^{\mathcal{F}}.$$

  We have

  (1) $(1+\varepsilon')^2\|D_v\tilde{f}_x(\xi)\| \le \|D_v\overline{F}_x(\xi)\| \le (1+\varepsilon')^2\|D_v\tilde{f}_x(\xi)\|$
  (2) $\mathrm{Höl}^\beta(D\overline{F}_x) \le (1+\varepsilon')^{2+\beta}(\varepsilon')^\beta\left(e^{(1+\beta)(\lambda_0+\varepsilon')} + \varepsilon'\right) + e^{(1+\beta)(\lambda_0 2\varepsilon')}(1+\varepsilon')^{1+\beta}\varepsilon'.$

From (1), $\overline{F}_x$ is well-defined. Taking $\varepsilon' > 0$ sufficiently small $\Phi_x$ is a family of $\varepsilon$-charts, and we can ensure $\overline{f}_x = \overline{F}_x\restriction_{V(e^{-\lambda_0 - 2\varepsilon})}$ is well defined, $\Psi_x^{\mathcal{F}}$ is a family of $\varepsilon$-charts, and find a function $\overline{\ell}$ satisfying the desired properties. $\qquad\square$

### 5.3. Local unstable manifolds and Proof of Proposition 4.5.
Relative to either the charts $\Phi_x$ in Proposition 5.1 or $\Psi_x^{\mathcal{F}}$ in Proposition 5.4 we may perform either the Perron–Irwin method or Hadamard graph transform method to construct (un)stable manifolds. See [FHY] or [HPS] for more details.

Fix $n \in \mathbb{Z}^d$. Let $F \subset \mathbb{Z}^d$ be a finite, symmetric, generating set containing $n$. Let $\mathcal{F}$ be an $\alpha$-invariant, $C^{1+\beta}$-tame measurable foliation. Let $U$ be as in Section 3.1 and Definition 4.2 and take $\varepsilon_0$ and $\lambda_0$ as above. Let $f = \alpha(n)$ and fix $0 < \varepsilon < \varepsilon_0$. Take $\Lambda = \Lambda_0$ or $\Lambda'$, $E = \mathbb{R}^k$ or $V$, and $f_x = \tilde{f}_x$ or $f_x = \overline{f}_x$, respectively, with the notation of either Proposition 5.1 or Proposition 5.4. write $E^i = \mathbb{R}^i \cap E$ and $E_{\le i} := \bigoplus_{1 \le j \le i} E^{\sigma(n)(j)}$ $E_{>i} := \bigoplus_{j > i} E^{\sigma(n)(j)}$.

**Lemma 5.5.** *If $\lambda_{\sigma(n)(i)}(n) > 0$ and $\lambda_{\sigma(n)(i)}(n) > \lambda_{\sigma(n)(i+1)}(n)$ then for every $x \in \Lambda$ there is a $C^{1+\beta}$ function*

$$h_x^i \colon E_{\le i}(e^{-\lambda_0 - 2\varepsilon}, \|\cdot\|) \to E_{>i}(1, \|\cdot\|)$$

*with*

  (a) *$h_x^i(0) = 0$, and $D_0 h_x^i(0) = 0$;*
  (b) *$\mathrm{Höl}^\beta(Dh_x^i) < c$ for some uniform constant $c$;*
  (c) *$\|Dh_x^i\| \le 1/3$;*
  (d) *for $\delta \le e^{-\lambda_0 - 2\varepsilon}$, writing $W_{x,\delta} := \mathrm{graph}\left(h_x^i\restriction_{E_{\le i}(\delta, \|\cdot\|)}\right)$ we have*

$$W_{f(x),\delta} \subset f_x(W_{x,\delta});$$

  (e) *if $u, v \in W_{x,\delta}$ then*

$$e^{\lambda_{\sigma(n)(i)}(n) - 2\varepsilon}\|v - u\| \le \|f_x(v) - f_x(w)\| \le e^{\lambda_{\sigma(n)(1)}(n) + 2\varepsilon}\|v - u\|.$$

*Moreover, the family $\{h_x^i\}$ depends measurably on $x \in \Lambda$. Write $f_x^{-k} = f_{f^{-k}(x)}^{-1} \circ \cdots \circ f_{f^{-1}(x)}^{-1}$ where defined. Then*

  (f) *for $u \in W_{x,\delta}$ and $n \ge 0$, $f_x^{-n}(u)$ is defined and*

$$W_{x,\delta} := \{u \in E(\delta, \|\cdot\|) : \limsup \tfrac{1}{k}\log\|f_x^{-k}(u)\| \le -\lambda_{\sigma(n)(i)}(n) + 10\varepsilon\}. \tag{16}$$

Continue to write $f = \alpha(n)$. Let $V_{\mathrm{loc},x,\varepsilon}$ be the image of $W_{x,e^{-2\lambda_0-4\varepsilon}}$ in Lemma 5.5 under either $\Phi_x$ or $\Psi_x^{\mathcal{F}}$. With $f = \alpha(n)$ we still have $f(V_{\mathrm{loc},x,\varepsilon}) \supset V_{\mathrm{loc},f(x),\varepsilon}$. Moreover, for $m \in F$ and a.e. $x \in F$ we have

(1) $\alpha(m)(V_{\mathrm{loc},x,\varepsilon})$ is contained in the image of $W_{f(x),e^{-\lambda_0-2\varepsilon}}$ under either $\Phi_{\alpha(m,x)}$ or $\Psi_{\alpha(m,x)}^{\mathcal{F}}$;

(2) $f^{-k}(\alpha(m)(V_{\mathrm{loc},x,\varepsilon})) \subset U$ for all $k \geq 0$;

(3) for $y \in V_{\mathrm{loc},x,\varepsilon}$
$$d(f^{-k}(\alpha(m,y), f^{-k}(\alpha(m,x)))) \leq \ell(x) e^{|m|\varepsilon} e^{-\lambda_{\sigma(n)(i)}(n)+2\varepsilon}.$$

*Proof of Proposition 4.5.* Let $\mathcal{F}$ be an $\alpha$-invariant, $C^{1+\beta}$-tame measurable foliation. In the case of (a) of Proposition 4.5, take $\mathcal{F} = \{M\}$. Recall we write $f = \alpha(n)$ for our distinguished $n$. Also recall we fix a finite symmetric generating set $F \ni n$ and set $U \subset U_0$ with $\alpha(m) \colon U \to U_0$ a diffeomorphism onto its image for each $m \in F$. Take $U \subset \hat{U} \subset U_0$ open such that $f\!\restriction_{\hat{U}} \to U_0$ is a diffeomorphism onto its image. We may then replace $U$ in the hypotheses of Section 3.1 with a smaller open set so that
$$\alpha(m)(U) \subset \hat{U}$$
for all $m \in F$.

With $\Lambda$ as in Lemma 5.5 and $x \in \Lambda$ take
$$V_x := \bigcup_{k \in \mathbb{Z}} (f\!\restriction_{\hat{U}})^k (V_{\mathrm{loc},f^{-k}(x),\varepsilon})$$
$$= \{y : f^{-k}(y) \in \hat{U} \text{ for all } k \geq 1 \text{ and } f^{-k}(y) \in V_{\mathrm{loc},f^{-k}(x),\varepsilon} \text{ for some } k \geq 0\}.$$

We have

**Claim 5.6.** *For $x, y \in \Lambda$ and $m \in F$*

(1) *if $y \in V_x \cap U$ then $\alpha(m)(y) \in V_{\alpha(m)(x)}$;*

(2) *$y \in V_x$ then $\limsup_{k \to -\infty} d(f^{-k}(y), f^{-k}(x)) \leq -\lambda_{\sigma(n)(i)}(n)$;*

(3) *if $y \notin V_x$ then $V_x \cap V_y = \varnothing$;*

(4) *if $y \in V_x$ then $V_x = V_y$.*

Take $B(n) = M \smallsetminus \left( \bigcup_{x \in \Lambda} V_x \right)$. From the above discussion, for $x \in \Lambda$ we have

(1) $V_x$ is a $C^1$ injectively immersed manifold and is defined independently of $\varepsilon$;

(2) $\left( \mathcal{F} \vee \mathscr{W}_n^{u,i} \right)(x) \smallsetminus B(n) = V_x$.

In particular,

(3) $\mathcal{F} \vee \mathscr{W}_n^{u,i}$ is an $\alpha$-invariant, $C^{1+\beta}$-tame, measurable foliation.

The proposition follows.                                                      $\square$

**Remark 5.7.** Given an $\alpha$-invariant, $C^{1+\beta}$-tame, measurable foliation $\mathcal{F}$ and sufficiently small $\varepsilon > 0$, let $\Phi_x$ be a family of $\varepsilon$-charts such the family $\Psi_x^{\mathcal{F}}$ constructed from $\Phi_x$ in Proposition 5.4 is a family of $\varepsilon$-charts.

Let $h_x^u$ be as in Lemma 5.5 for the charts $\Phi_x$ and $i = r$. Let $\tilde{h}_x^{\mathcal{F}}$ be as in Proposition 5.2 and let $h_x^{\mathcal{F},u}$ be as in Lemma 5.5 for the charts $\Psi_x$ with $i = r$. Let $W_{x,\delta}^u :=$ graph $\left( h_x^u\!\restriction_{\mathbb{R}_{\leq r}(\delta,\|\cdot\|)} \right)$ and let $\tilde{\mathcal{F}}_x$ be the graph of $\tilde{h}_x^{\mathcal{F}}$.

Let $W_x^{\mathcal{F},u}$ be the graph of $\tilde{h}_x^{\mathcal{F}} \circ h_x^{\mathcal{F},u}$ and let $W_x^{\mathcal{F},u}(\delta) = W_x^{\mathcal{F},u} \cap \mathbb{R}^k(\delta)$.

It follows from the characterization (16), the Lipschitness of $H_x$ in the proof of Proposition 5.4, and the local dynamics and invariance of manifolds in charts that for all $\delta < e^{-\lambda_0-2\varepsilon}$, $W_x^{\mathcal{F},u}(\delta) = \tilde{\mathcal{F}}_x \cap W_{x,\delta}^u$.

UNIVERSITY OF CHICAGO, CHICAGO, IL 60637, USA
*E-mail address*: `awb@uchicago.edu`

PENNSYLVANIA STATE UNIVERSITY, STATE COLLEGE, PA 16802, USA
*E-mail address*: `hertz@math.psu.edu`

# SMOOTH ERGODIC THEORY OF $\mathbb{Z}^d$-ACTIONS PART 2: ENTROPY FORMULAS FOR RANK-1 SYSTEMS

## AARON BROWN

ABSTRACT. We show the main result of [LY1] and the entropy formulas from [LY2] continue to hold for diffeomorphisms satisfying our standing hypotheses. Moreover, we establish similar results for notions of entropy subordinate to a foliation or a measurable partition. As a corollary of the proof, we obtain the finiteness of entropy for systems satisfying our standing hypotheses of Section 3.1.

In this part, we extend and generalize the main result of [LY1] and the entropy formulas from [LY2] to the setting of diffeomorphisms satisfying our standing hypotheses. As a corollary of the proof we obtain the finiteness of entropy for systems satisfying our standing hypotheses of Section 3.1. Although we provide most details here, the arguments in this section are adapted from the original papers [LY1, LY2]

## 6. DEFINITIONS AND FACTS ABOUT METRIC ENTROPY

We recall some standard facts about the metric entropy of measure-preserving transformations as well as the main definitions used in this part. Let $(X, \mu)$ be a standard probability space. Given a measurable partition $\xi$ of $(X, \mu)$, for the remainder we indicate by $\{\mu_x^\xi\}$ a family of conditional probability measures relative to the partition $\xi$. That is, the assignment $x \mapsto \mu_x^\xi$ is $\sigma(\xi)$-measurable (where $\sigma(\xi)$ denotes the $\sigma$-algebra of $\xi$-saturated sets) and given a measurable $A \subset X$,

$$\mu(A) = \int \mu_x^\xi(A) \ d\mu(x).$$

Given measurable partitions $\eta, \xi$ of $(X, \mu)$, the *mean conditional information* of $\eta$ relative to $\xi$ is $I_\mu(\eta \mid \xi)(x) = -\log(\mu_x^\xi(\eta(x)))$ and the *mean conditional entropy* of $\eta$ relative to $\xi$ is

$$H_\mu(\eta \mid \xi) = \int I_\mu(\eta \mid \xi)(x) \ d\mu(x).$$

The *entropy of* $\eta$ is $H_\mu(\eta) = H_\mu(\eta \mid \{\varnothing, X\})$. Note that if $H_\mu(\eta) < \infty$ then $\eta$ is necessarily countable. Let $f \colon (X, \mu) \to (X, \mu)$ be an invertible, measurable, measure-preserving transformation. Let $\eta$ be an arbitrary measurable partition of $(X, \mu)$. We define

$$\eta^+ := \bigvee_{i=0}^{\infty} f^i \eta, \qquad \eta^f := \bigvee_{i \in \mathbb{Z}}^{\infty} f^i \eta.$$

We define the ("unstable" or "future[1]") *entropy of* $f$ *given the partition* $\eta$ to be

$$h_\mu(f, \eta) := H_\mu(\eta \mid f\eta^+) = H_\mu(\eta^+ \mid f\eta^+) = H_\mu(f^{-1}\eta^+ \mid \eta^+).$$

---

[1] It is more standard to define the entropy $h_\mu(f, \eta)$ as $H_\mu(\eta \mid f^{-1}\eta^-)$. Note that we typically expect asymmetry of these definitions: $H_\mu(\eta \mid f\eta^+) \neq H_\mu(\eta \mid f^{-1}\eta^-)$. However, if $\eta$ satisfies $H_\mu(\eta) < \infty$ then the symmetry of the two definitions holds. We choose to define $h_\mu(f, \eta) = H_\mu(\eta \mid f\eta^+)$ to have results most consistent with statements in [LY1, LY2] and related work.

We define the *$\mu$-entropy of $f$* to be $h_\mu(f) = \sup\{h_\mu(f, \eta)\}$ where the supremum is taken over all measurable partitions of $(X, \mu)$.

### 6.0.1. *Entropy subordinate to a partition.*

**Definition 6.1.** Given a measurable partition $\eta$ of $(M, \mu)$ we define the *entropy of $f$ subordinate to $\eta$* to be the quantity

$$h_\mu(f \mid \eta) := \sup\{h_\mu(f, \xi) : \xi \geq \eta\}.$$

In general, $h_\mu(f \mid \eta) \neq h_\mu(f^{-1} \mid \eta)$. However, if $\eta$ is an $f$-invariant partition then

$$h_\mu(f \mid \eta) = h_\mu(f^{-1} \mid \eta).$$

### 6.0.2. *Entropy subordinate to a measurable foliation.*

We now take $X = M$ to be a $C^\infty$ manifold equipped with a Borel probability measure $\mu$. Consider a measurable foliation $\mathcal{F}$ of $M$. Note that the partition into leaves of $\mathcal{F}$ is generally not a measurable partition of $(M, \mu)$.

**Definition 6.2.** A measurable partition $\eta$ of $(M, \mu)$ is said to be *subordinate to $\mathcal{F}$* if for a.e. $x \in M$

(1) $\eta(x) \subset \mathcal{F}'(x)$ (where $\mathcal{F}'(x)$ is as in Section 4.1);
(2) $\eta(x)$ contains an open neighborhood (in the immersed topology) of $x$ in $\mathcal{F}'(x)$;
(3) $\eta(x)$ is precompact in the internal topology of $\mathcal{F}'(x)$ for a positive measure set of $x$.

**Definition 6.3.** Suppose $\mu$ is $f$-ergodic. Let $\mathcal{F}$ be an $f$-invariant, measurable foliation of $(M, \mu)$. We define the *entropy of $f$ subordinate to $\mathcal{F}$* to be

$$h_\mu(f \mid \mathcal{F}) := \sup\{h_\mu(f \mid \xi) : \xi \text{ is subordinate to } \mathcal{F}\}.$$

More generally, if $\eta$ is any measurable partition of $(M, \mu)$ we define

$$h_\mu(f \mid \eta \vee \mathcal{F}) := \sup\{h_\mu(f \mid \xi \vee \eta) : \xi \text{ is subordinate to } \mathcal{F}\}.$$

### 6.1. **Properties of metric entropy.**

We recall some properties of the above definitions. A primary reference for proofs is [Rok]. Consider a standard probability space $(X, \mu)$, an invertible, measurable, measure-preserving transformation $f: (X, \mu) \to (X, \mu)$, and measurable partitions $\eta, \xi$, and $\zeta$ of $(X, \mu)$.

(1) $h_\mu(f, \xi) \leq H_\mu(\xi)$.
(2) $I_\mu(\eta \vee \zeta \mid \xi)(x) = I_\mu(\eta \mid \xi)(x) + I_\mu(\zeta \mid \eta \vee \xi)(x)$ a.e. whence $H_\mu(\eta \vee \zeta \mid \xi) = H_\mu(\eta \mid \xi) + H_\mu(\zeta \mid \eta \vee \xi)$.
(3) $h_\mu(f, \eta \vee \zeta) \leq h_\mu(f, \eta) + h_\mu(f, \zeta)$.
(4) If $\eta \geq \xi$ then $H_\mu(\eta \mid \zeta) \geq H_\mu(\xi \mid \zeta)$ and $H_\mu(\zeta \mid \eta) \leq H_\mu(\zeta \mid \xi)$.
(5) If $\zeta_n \nearrow \zeta$ and if $H_\mu(\eta \mid \zeta_1) < \infty$ then $I_\mu(\eta \mid \zeta_n)(x) \searrow I_\mu(\eta \mid \zeta)$ and $H_\mu(\eta \mid \zeta_n) \searrow H_\mu(\eta \mid \zeta)$.
(6) If $\mu_x^\eta(\xi(x)) > 0$ for almost every $x$ then $h_\mu(f, \eta \vee \xi) \geq h_\mu(f, \eta)$; in particular, if $H_\mu(\mathcal{P}) < \infty$ then $h_\mu(f, \eta \vee \mathcal{P}) \geq h_\mu(f, \eta)$.
(7) $h_\mu(f) = \sup\{h_\mu(f, \mathcal{P}) : H_\mu(\mathcal{P}) < \infty\}$ and $h_\mu(f \mid \eta) = \sup\{h_\mu(f, \eta \vee \mathcal{P}) : H_\mu(\mathcal{P}) < \infty\}$.
(8) $h_\mu(f, \eta \vee \xi) = h_\mu(f, \eta \vee f^k(\xi))$ for $k \in \mathbb{Z}$.
(9) If $\xi \geq \eta$ and $H_\mu(\xi \mid f\eta^+) < \infty$ then $\frac{1}{n} H_\mu\left(\bigvee_0^{n-1} f^n(\xi) \mid f^n\eta^+\right) \searrow h_\mu(f, \xi)$.
(10) If either $h_\mu(f, \xi) < \infty$ or $H_\mu(\xi \mid \eta) < \infty$ then

$$h_\mu(f, \xi \vee \eta) \leq h_\mu(f, \eta) + h_\mu(f, \xi \vee \eta^f).$$

We note that all inequalities hold for $\infty$-valued quantities.

(9) follows from [Rok, 7.3]. (10) holds as

$$
\begin{aligned}
h_\mu(f, \xi \vee \eta) &= H_\mu(\eta \vee f^n(\xi) \mid f(\eta^+) \vee f^{n+1}(\xi^+)) \\
&= H_\mu(\eta \mid f(\eta^+) \vee f^{n+1}(\xi^+)) + H_\mu(f^n(\xi) \mid \eta^+ \vee f^{n+1}(\xi^+)) \\
&= H_\mu(\eta \mid f(\eta^+) \vee f^{n+1}(\xi^+)) + H_\mu(\xi \mid f^{-n}(\eta^+) \vee f(\xi^+)).
\end{aligned}
$$

We have $H_\mu(\eta \mid f(\eta^+) \vee f^{n+1}(\xi^+)) \leq h_\mu(f, \eta)$ and, assuming either $h_\mu(f, \xi) < \infty$ or $H_\mu(\xi \mid \eta) < \infty$, $H_\mu(\xi \mid f^{-n}(\eta^+) \vee f(\xi^+))$ decreases to $H_\mu(\xi \mid \eta^f \vee f(\xi^+))$ as $n \to \infty$. As $\eta^f$ is $f$-invariant, we have $H_\mu(\xi \mid \eta^f \vee f(\xi^+)) = h_\mu(f, \xi \vee \eta^f)$.

For (6) first note that there is a countable partition $\mathcal{P}$ with $\eta \vee \xi = \eta \vee \mathcal{P}$. Also note that the result holds whenever $\eta$ and $\xi$ are finite partitions. Take $\eta_n$ and $\mathcal{P}_n$ to be finite partitions with $\eta_n \nearrow \eta$ and $\mathcal{P}_n \nearrow \mathcal{P}$. If the inequality fails take $M < \infty$ with

$$
h_\mu(f, \eta \vee \xi) = H_\mu(\eta \vee \mathcal{P} \mid f(\eta^+ \vee \mathcal{P}^+)) < M < H_\mu(\eta \mid f(\eta^+)) = h_\mu(f, \eta).
$$

We have $H_\mu(\eta_n \mid f(\eta^+)) \nearrow H_\mu(\eta \mid f(\eta^+)) = h_\mu(f, \eta)$. Since $H_\mu(\eta_n \mid f(\eta_n^+)) \geq H_\mu(\eta_n \mid f(\eta^+))$ take $n_1$ with $M < H_\mu(\eta_n \mid f(\eta_n^+))$ for all $n \geq n_1$. As

$$
\begin{aligned}
h_\mu(f, \eta_n \vee \mathcal{P}_n) &= H_\mu(\eta_n \vee \mathcal{P}_n \mid f(\eta_n^+ \vee \mathcal{P}_n^+)) \searrow H_\mu(\eta_n \vee \mathcal{P}_n \mid f(\eta^+ \vee \mathcal{P}^+)) \\
&\leq H_\mu(\eta \vee \mathcal{P} \mid f(\eta^+ \vee \mathcal{P}^+)) < M
\end{aligned}
$$

take $n_2 \geq n_1$ so that $H_\mu(\eta_n \vee \mathcal{P}_n \mid f(\eta_n^+ \vee \mathcal{P}^+)) < M$ for all $n \geq n_2$. Then $h_\mu(f, \eta_n \vee \mathcal{P}_n) < h_\mu(f, \eta_n)$ for $n \geq n_2$, a contradiction.

### 6.2. Abramov–Rohlin-type formulas.

We say a partition $\eta$ of $(X, \mu)$ is $f$-increasing if $f\eta \leq \eta$. Note in this case $\eta^+ = \eta$ and $h_\mu(f, \eta) = H_\mu(\eta \mid f\eta)$. We say $\eta$ generates if $\eta^f$ is the point partition.

From (9) and (10) of Section 6.1 we obtain the following

**Proposition 6.4.** *Let $\eta$ be an $f$-increasing, measurable partition of $(X, \mu)$. If $H_\mu(\xi \mid \eta) < \infty$ then*

$$
h_\mu(f, \xi \vee \eta) = h_\mu(f, \eta) + h_\mu(f, \xi \vee \eta^f). \tag{17}
$$

Let $g \colon (Y, \nu) \to (Y, \nu)$ be second invertible, measure-preserving transformation of a standard probability space. Suppose there is a measurable $\psi \colon X \to Y$ with $\psi_* \mu = \nu$ and $\psi \circ f = g \circ \psi$. We say that $g$ is a *factor of $f$ induced by $\psi$*. Write $\mathcal{A}^\psi$ for the partition of $(X, \mu)$ into preimages of $\psi$. Note $\mathcal{A}^\psi$ is $f$-invariant.

Taking the supremum over all partitions $\mathcal{P}$ of $(X, \mu)$ with $H_\mu(\mathcal{P}) < \infty$ and over $\eta = \mathcal{Q}^+$ for all partitions $\mathcal{Q}$ of $(Y, \nu)$ with $H_\nu(\mathcal{Q}) < \infty$, from Proposition 6.4 we obtain the following equalities including the classical Abramov–Rohlin formula (see [LW, BC]).

**Corollary 6.5.** *Let $g \colon (Y, \nu) \to (Y, \nu)$ be a measurable factor of $f \colon (X, \mu) \to (X, \mu)$ induced by $\psi$. Let $\hat\eta$ be a measurable partition of $(Y, \nu)$ that is increasing for $g$. Then*

$$
h_\mu(f \mid \psi^{-1}(\hat\eta)) = h_\nu(g, \hat\eta) + h_\mu(f \mid \psi^{-1}(\hat\eta^g))
$$

*and if $\hat\eta$ generates for $g$ then*

$$
h_\mu(f \mid \psi^{-1}(\hat\eta)) = h_\nu(g, \hat\eta) + h_\mu(f \mid \mathcal{A}^\psi). \tag{18}
$$

*In particular,*

$$
h_\mu(f) = h_\nu(g) + h_\mu(f \mid \mathcal{A}^\psi). \tag{19}
$$

*Proof of Proposition 6.4.* Note that if $h_\mu(f, \eta) = \infty$ we have equality from (6) of Section 6.1. We thus assume $h_\mu(f, \eta) < \infty$. Take $\xi$ with $H_\mu(\xi \mid \eta) < \infty$. By (10),

$$h_\mu(f, \xi \vee \eta) \leq h_\mu(f, \eta) + h_\mu(f, \xi \vee \eta^f).$$

For the reverse inequality, we have $H_\mu(\xi \vee \eta \mid f\eta) = H_\mu(\xi \mid f\eta) + H_\mu(\eta \mid f\eta) < \infty$ whence from (9) we have $\frac{1}{n} H_\mu(\bigvee_0^{n-1} f^n(\xi \vee \eta) \mid f^n \eta) \searrow h(f, \xi \vee \eta)$. But,

$$H_\mu \left( \bigvee_0^{n-1} f^n(\xi \vee \eta) \mid f^n \eta \right)$$

$$= H_\mu \left( \bigvee_0^{n-1} f^n \eta \mid f^n \eta \right) + H_\mu \left( \bigvee_0^{n-1} f^n \xi \mid f^n \eta \vee \bigvee_0^{n-1} f^n \eta \right)$$

$$= H_\mu \left( \eta \mid f^n \eta \right) + H_\mu \left( \bigvee_0^{n-1} f^n \xi \mid \eta \right)$$

$$\geq n h_\mu(f, \eta) + H_\mu \left( \bigvee_0^{n-1} f^n \xi \mid \eta^f \right).$$

Dividing by $n$ and taking $n \to \infty$ we have $h_\mu(f, \eta) + h_\mu(f, \xi \vee \eta^f) \leq h(f, \xi \vee \eta)$. $\square$

## 7. Statement of results

We begin with an extension of the main results of [LY1, LY2] to the setting of rank-1 systems satisfying the hypotheses introduced in Section 3.1. Let $M$ be a $k$-dimensional, $C^\infty$ manifold equipped with a Borel probability measure $\mu$. Let $f \colon M \to M$ be an invertible, measurable, $\mu$-preserving transformation. Then $f$ generates an action of $\mathbb{Z}$ on $M$. We assume the induced $\mathbb{Z}$-action satisfies the hypotheses of Section 3.1. In particular, we fix the generating set $F = \{-1, 1\}$ and consider all constructions from Part 1 including the set $U \subset U_0$ from Section 3.1 and the Lyapunov charts from Section 5 to be relative to $F$.

For the remainder of Part 2 we further assume that $\mu$ is $f$-ergodic as the generalizations of all results stated here to the non-ergodic case are standard.

As $f$ induces an action of $\mathbb{Z}$, we may identify the Lyapunov exponent functionals $\lambda_i \colon \mathbb{Z} \to \mathbb{R}$ for the derivative cocycle (12) with their coefficients $\lambda_i := \lambda_i(1) \in \mathbb{R}$. Thus, for the remainder of Part 2, Lyapunov exponents are numbers and are listed without multiplicity and ordered so that

$$\lambda_1 > \lambda_2 > \cdots > \lambda_r > 0 \geq \lambda_{r+1} > \cdots > \lambda_p.$$

We write $E^i(x) = E_{\lambda_i}(x)$ and write $m_i$ for the almost-surely constant value of $\dim E^i(x)$. Also write $E^u = \oplus_{\lambda_i > 0} E^i$, $E^0 = E_0$, and $E^s = \oplus_{\lambda_i < 0} E_{\lambda_i}$ Given $1 \leq i \leq r$ we write $W^i(x) = W^i_1(x)$ and $\mathscr{W}^i = \mathscr{W}^i_1$ for the $i$th unstable manifold and $i$th unstable foliation corresponding to the generator $f$. Similarly write $\mathscr{W}^u = \mathscr{W}^u_1$. Given any $C^{1+\beta}$-tame, $f$-invariant, measurable foliation $\mathcal{F}$ we write $\mathcal{F}^u := \mathcal{F} \vee \mathscr{W}^u$. We similarly define stable manifolds $W^s(x)$ and the stable foliation $\mathscr{W}^s$ to be the unstable manifold and foliation relative to $f^{-1}$.

### 7.1. Finiteness of entropy.
Our first result is the following version of the Margulis–Ruelle inequality for systems satisfying our standing hypotheses.

**Proposition 7.1.** *Let $(M, \mu)$ and $f$ be as above. Then*

$$h_\mu(f) \leq \sum_{\lambda_i > 0} \lambda_i m_i.$$

*In particular, $h_\mu(f) < \infty$.*

Note, in particular, that we do not assume any control on the capacity of the manifold $M$ as in [KSLP] to obtain finiteness of entropy. On the other hand, the finiteness of entropy in [KSLP, Theorem IV.1.1] uses only $C^1$ estimates on the dynamics.

### 7.2. Geometric rigidity of measures satisfying the entropy formula.
Let $\mathcal{F}$ be an $f$-invariant, measurable foliation. Define the multiplicity of $\lambda_i$ relative to $\mathcal{F}$ to be (the almost-surely constant value of)

$$m_i(\mathcal{F}) := \dim \left( E^i(x) \cap T_x \mathcal{F}'(x) \right).$$

We say a measure $\mu$ is *absolutely continuous along* $\mathcal{F}$ if $\mu_x^\eta$ is absolutely continuous with respect to the Riemannian volume on $\eta_x \subset \mathcal{F}(x)'$ for any measurable partition $\eta$ subordinated to $\mathcal{F}$ and almost every $x$.

**Theorem 7.2.** *Let $\mathcal{F}$ be an $f$-invariant, $C^{1+\beta}$-tame, measurable foliation. Then*

$$h_\mu(f \mid \mathcal{F}) \leq \sum_{1 \leq i \leq r} \lambda_i m_i(\mathcal{F}). \tag{20}$$

*Moreover, equality holds if and only if $\mu$ is absolutely continuous along $\mathcal{F}^u$.*

Furthermore, in the case of equality in (20), the measures $\mu_x^\eta$ are *equivalent* to the Riemannian volume on $\eta_x \subset \mathcal{F}^u(x)'$ for any measurable partition $\eta$ subordinated to $\mathcal{F}^u$. (See [LY1, Corollary 6.1.4]).

Note that Proposition 7.1 follows from the statement of Theorem 7.2 by taking $\mathcal{F} = M$. However, we establish Proposition 7.1 separately from Theorem 7.2.

As in the main result of [Led], Theorem 7.2 follows from Jensen's inequality after we establish that all entropy of the system is carried by unstable manifolds. That is,

**Proposition 7.3.** *For an $f$-invariant, $C^{1+\beta}$-tame, measurable foliation $\mathcal{F}$ we have*

$$h_\mu(f \mid \mathcal{F}) = h_\mu(f \mid \mathcal{F}^u).$$

### 7.3. Pointwise transverse dimensions.
Let $\eta$ be an arbitrary measurable partition of $(M, \mu)$. Given $1 \leq i \leq r$, let $\xi^i$ be a measurable partition subordinate to $\mathscr{W}^i$. Let $\{\mu_x^{\xi^i \vee \eta}\}$ be a family of conditional measures relative to the measurable partition $\xi^i \vee \eta$. We define the $i$th upper and lower pointwise dimensions of $\mu$ relative to $\eta$ at $x$ to be

$$\overline{\dim}^i(\mu, x|\eta) := \limsup_{r \to 0} \frac{\log \left( \mu_x^{\xi^i \vee \eta}(B(x, r)) \right)}{\log r}$$

$$\underline{\dim}^i(\mu, x|\eta) := \liminf_{r \to 0} \frac{\log \left( \mu_x^{\xi^i \vee \eta}(B(x, r)) \right)}{\log r}$$

and the $i$th pointwise dimension of $\mu$ relative to $\eta$ at $x$ to be

$$\dim^i(\mu, x|\eta) := \lim_{r \to 0} \frac{\log \left( \mu_x^{\xi^i \vee \eta}(B(x, r)) \right)}{\log r}$$

when the limit exists. One verifies that the functions $\overline{\dim}^i(\mu, x|\eta)$ and $\underline{\dim}^i(\mu, x|\eta)$ are measurable and independent of the choice of $\xi^i$. Moreover, if $f\eta \leq \eta$ and $h_\mu(f, \eta) < \infty$,

it follows that $\overline{\dim}^i(\mu, x|\eta)$ and $\underline{\dim}^i(\mu, x|\eta)$ are constant along orbits of $f$ and hence, by ergodicity of $\mu$, constant a.s. Let $\overline{\dim}^i(\mu|\eta)$, $\underline{\dim}^i(\mu|\eta)$ denote the a.s. constant values.

We have

**Proposition 7.4.** *Let $\eta$ be a measurable partition of $(M, \mu)$. Then*

$$\overline{\dim}^i(\mu|\eta^+) = \underline{\dim}^i(\mu|\eta^+)$$

Set $\dim^0(\mu|\eta^+) = 0$. For $1 \leq i \leq r$ define the *$i$th transverse dimension of $\mu$ relative to $\eta^+$* to be

$$\gamma^i(\mu|\eta^+) = \dim^i(\mu|\eta^+) - \dim^{i-1}(\mu|\eta^+).$$

We claim

**Claim 7.5.** *For $1 \leq i \leq r$, $\gamma^i(\mu|\eta^+) \leq m_i := \dim E^i$.*

Note that if $\hat{\eta}$ is a measurable partition of $(M, \mu)$ with $\hat{\eta} \geq \eta$ then it will follow from Proposition 8.4 and Claim 8.6 below that

$$h_\mu(f \mid \eta) \geq h_\mu(f \mid \hat{\eta}).$$

We have a similar result for the transverse dimensions above.

**Proposition 7.6.** *If $\hat{\eta}^+ \geq \eta^+$ then for every $1 \leq i \leq r$*

$$\gamma^i(\mu|\hat{\eta}^+) \leq \gamma^i(\mu|\eta^+).$$

Proposition 7.6 follows from discussion at the end of Section 11.3.

### 7.4. Geometric characterization of the defect in the entropy formula.

As in [LY2] we have an explicit geometric description of the defect of equality in (20).

**Theorem 7.7.** *Let $\eta$ be any measurable partition of $(M, \mu)$. Then*

$$h_\mu(f \mid \eta) = \sum_{1 \leq i \leq r} \lambda_i \gamma_i(\mu|\eta^+).$$

Taking $\eta$ to be the trivial partition $\{X, \varnothing\}$ we obtain an extension of the entropy formula from [LY2] to $C^{1+\beta}$ diffeomorphisms of noncompact manifolds satisfying our standing hypotheses.

Consider $g\colon (Y, \nu) \to (Y, \nu)$ a measurable factor of $f$ induced by $\psi\colon (M, \mu) \to (Y, \nu)$. As a primary application of Theorem 7.7, we obtain a Ledrappier–Young entropy formula for the fiber entropy $h_\mu(f \mid \mathcal{A}^\psi)$ of smooth systems when the elements of the fiber partition $\mathcal{A}^\psi$ are only measurable. In particular, from (19) the entropy of the factor system $g\colon (Y, \nu) \to (Y, \nu)$ can be computed in terms of the Lyapunov exponents of the total system $f$ and the geometry of the conditional measures of $\mu$ and $\{\mu_x^{\mathcal{A}^\psi}\}$ along unstable manifolds. In the case that the fibers are smooth manifolds, a Ledrappier–Young formula for fiber entropy follows from [QQX].

### 7.5. A characterization of the Pinsker Partition.

Let $\pi$ denote the Pinsker partition for the action of $f$ on $(M, \mu)$. Also, let $\mathscr{B}^u$ and $\mathscr{B}^s$ be the measurable hulls of the (typically non-measurable) partitions $\mathscr{W}^u$ and $\mathscr{W}^s$.

**Theorem 7.8.** *Under the above assumptions*

$$\pi \stackrel{\circ}{=} \mathscr{B}^u \stackrel{\circ}{=} \mathscr{B}^s.$$

Theorem 7.8 follows from the discussion in [Rok, Section 12] exactly as in [LY1] from Proposition 8.4 and Remark 8.1 below.

## 8. Preparations: special partitions and their entropy properties

### 8.1. Standing notation.
We write $\Lambda_0$ for the set of regular points of $\mu$. Let

$$\varepsilon_0 = \min\{1, |\lambda_i|, |\lambda_i - \lambda_j| : i \neq j, \lambda_i \neq 0\}/100.$$

For Sections 8–11, fix $0 < \varepsilon < \varepsilon_0$. Let $\{\Phi_x\}$ be a fixed family of $\varepsilon$-Lyapunov charts as in Proposition 5.1 with corresponding function $\ell$. Given $0 \leq i \leq r$ and $x \in \Lambda_0$, let $h_x^i$ be as in Lemma 5.5 (relative to the charts $\Phi_x$). Let $\lambda_0 = \max\{|\lambda_i|\}$. For $0 < \delta \leq e^{-\lambda_0 - 2\varepsilon}$ let

$$W_{x,\delta}^i = \operatorname{graph}(h_x^i \!\restriction_{\bigoplus_{j \leq i} \mathbb{R}^j (\delta)})$$

and write

$$V_{\mathrm{loc},x,\varepsilon}^i := \Phi_x(\operatorname{graph}(h_x^i \!\restriction_{\bigoplus_{j \leq i} \mathbb{R}^j (-2\lambda_0 - 4\varepsilon)})) = \Phi_x(W_{x,e^{-2\lambda_0 - 4\varepsilon}}^i)$$

for the *local $i$th unstable manifold* relative to the charts $\{\Phi_x\}$.

Consider $\mathcal{F}$ a $C^{1+\beta}$-tame, $f$-invariant, measurable foliation. We may choose $\{\Phi_x\}$ above so that the charts $\Psi_x^{\mathcal{F}}$ built from $\Phi_x$ as in Proposition 5.4 as in Proposition 5.4 is a family of $\varepsilon$-charts. Let $h_x^{\mathcal{F},u} = h_x^r$ be as in Lemma 5.5 (relative to the charts $\{\Psi_x^{\mathcal{F}}\}$) and write

$$V_{\mathrm{loc},x,\varepsilon}^{\mathcal{F},u} := \Psi_x^{\mathcal{F}}(\operatorname{graph}(h_x^{\mathcal{F},u} \!\restriction_{\bigoplus_{j \leq r} (V \cap \mathbb{R}^j)(e^{-2\lambda_0 - 4\varepsilon})}))$$

for the local manifold of $\mathcal{F} \vee \mathscr{W}^u$ through $x$ relative to the charts $\{\Psi_x^{\mathcal{F}}\}$. Also write

$$V_{\mathrm{loc},x,\varepsilon}^{\mathcal{F}} := \Psi_x^{\mathcal{F}}(V(e^{-\varepsilon})) = \Phi_x(\operatorname{graph}(\tilde{h}_x^{\mathcal{F}} \!\restriction_{V(e^{-\varepsilon})}))$$

where $\tilde{h}_x^{\mathcal{F}}$ is as in Proposition 5.2.

Recall the ambient metric $d$ on $U_0$. Write $d_{V_{\mathrm{loc},x,\varepsilon}^i}$ for the restriction of $d$ to the embedded manifold $V_{\mathrm{loc},x,\varepsilon}^i$ where $d_{V_{\mathrm{loc},x,\varepsilon}^i}(x,y) = \infty$ if $y \notin V_{\mathrm{loc},x,\varepsilon}^i$. Denote by $B(x,\delta) := \{y \in U_0 : d(x,y) < \delta\}$ and

$$B^i(x,\delta) := \{y \in V_{\mathrm{loc},x,\varepsilon}^i : d_{V_{\mathrm{loc},x,\varepsilon}^i}(x,y) < \delta\}.$$

### 8.2. Expanding partitions subordinate to unstable manifolds.
It follows from the same constructions in [LY2, Lemma 9.1.1] (see also [LS] for aditional details) that given $1 \leq i \leq r$ there exists a measurable partition $\xi^i$ with

(1) $\xi^i$ subordinate to $\mathscr{W}^i$ (see Section 6.0.2); moreover $\xi^i(x) \subset V_{\mathrm{loc},x,\varepsilon}^i$ for a positive measure set of $x$;
(2) $\xi^i \leq f^{-1}\xi^i$;
(3) $\bigvee_{n=0}^{\infty} f^{-n}\xi^i$ is the point partition;
(4) for $1 \leq i \leq r-1$, $\xi^{i+1} \leq \xi^i$.

A partition satisfying 1–3 is said to be an *increasing generator subordinate to $\mathscr{W}^i$*.

**Remark 8.1.** Note that with $i = r$, for any partition $\xi^u = \xi^r$ satisfying 1 and 2 we have

$$\bigwedge_{n \geq 0} f^n \xi^u = \mathscr{B}^u.$$

In the sequel, we will need an additional property of the partitions $\xi^i$ given by Remark 8.2 below as well as certain sets used in their construction. Additionally, we will need partitions $\xi^{\mathcal{F}^u}$ subordinate to $\mathcal{F}^u$ with similar properties as above. For this reason we briefly review their construction.

Recall our fixed $0 < \varepsilon < \varepsilon_0$ and family of $\varepsilon$-Lyapunov charts $\Phi_x$ with corresponding function $\ell$. Also recall the charts $\Psi_x^{\mathcal{F}}$ built from $\Phi_x$. Let $\Lambda_1 \subset \Lambda_0$ be a set of positive measure such that

- $\ell(x) \leq \ell_0$ for $x \in \Lambda_1$;
- the charts $x \mapsto \Phi_x$ and $x \mapsto \Psi_x^{\mathcal{F}}$ vary continuously in the $C^1$ topology on $\Lambda_1$;
- for $1 \leq i \leq r$ the family of functions $h_x^{u,i}$ vary continuously in the $C^1$ topology;
- the family of functions $h_x^{\mathcal{F},u}$ vary continuously in the $C^1$ topology.

It follows that the family of embedded manifolds $x \mapsto V_{\mathrm{loc},x,\varepsilon}^i$, $x \mapsto V_{\mathrm{loc},x,\varepsilon}^{\mathcal{F}}$, and $x \mapsto V_{\mathrm{loc},x,\varepsilon}^{\mathcal{F},u}$ vary continuously in the $C^1$ topology on $\Lambda_1$. Let $x_0$ be a density point of $\Lambda_1$. Given $\rho > 0$ sufficiently small and $y \in B(x_0, \rho) \cap \Lambda_1$, let $V^i(y, \rho)$ be the connected component of $V_{\mathrm{loc},y,\varepsilon}^i \cap B(x_0, \rho)$ containing $y$. Write $V^u(y, \rho) = V^r(y, \rho)$. Similarly define $V^{\mathcal{F},u}(y, \rho)$ for the connected component of $V_{\mathrm{loc},y,\varepsilon}^{\mathcal{F},u} \cap B(x_0, \rho)$ containing $y$. Set $S_\rho = B(x_0, \rho) \cap \Lambda_1$.

Taking $0 < \rho_0$ sufficiently small, we may assume that $N = B(x_0, \rho_0)$ is an embedded manifold with $\overline{N} \subset U$ an embedded manifold with boundary and that for all $0 < \rho < \rho_0$ and $y \in S_\rho$

(a) $V^{u,i}(y, \rho)$ is a connected open neighborhood of $y$ in $V_{\mathrm{loc},y,\varepsilon}^i$;
(b) the elements $\{V^{u,i}(y, \rho) \cap S_\rho : y \in S_\rho\}$ form a partition of $S_\rho$;
(c) if $y_1, y_2 \in S_\rho$ and $y_1 \notin V^{u,i}(y_2, \rho)$ then $y_1 \notin V_{\mathrm{loc},y_2,\varepsilon}^i$;
(d) $V^{\mathcal{F},u}(y, \rho)$ is a connected open neighborhood of $y$ in $V_{\mathrm{loc},y,\varepsilon}^{\mathcal{F},u}$ and

$$V^{\mathcal{F},u}(y, \rho) = V^u(y, \rho) \cap V_{\mathrm{loc},y,\varepsilon}^{\mathcal{F}}.$$

(e) the elements $\{V^{\mathcal{F},u}(y, \rho) \cap S_\rho : y \in S_\rho\}$ form a partition of $S_\rho$;
(f) if $y_1, y_2 \in S_\rho$ and $y_1 \notin V^{\mathcal{F},u}(y_2, \rho)$ then $y_1 \notin V_{\mathrm{loc},y_2,\varepsilon}^{\mathcal{F},u}$.

Let $\hat{\xi}^i$ be the partition

$$\hat{\xi}^i(y) = \begin{cases} V^{u,i}(y, \rho) & y \in S_\rho \\ M \smallsetminus S_\rho & y \notin S_\rho. \end{cases}$$

Take $\xi^i := (\hat{\xi}^i)^+$. Properties 2 and 3 follow immediately. For 4 note that $\hat{\xi}^{i+1} \leq \hat{\xi}^i$ for $1 \leq i \leq r - 1$ whence the conclusion follows by construction. Finally, property 1 follows for Lebesgue almost-every choice of $r < r_0$ by the same arguments as in [LS]. We fix such an $\rho$ and set $S = S_\rho$ for the remainder.

We write $\xi^u = \xi^r$.

**Remark 8.2.** $\xi^i$ has the following characterization: $y \in \xi^i(x)$ if and only if for every $m \leq 0$ such that $f^m(x) \in S$, we have $f^m(y) \in S$ and $f^m(y) \in V_{\mathrm{loc},f^m(x),\varepsilon}^i$.

**Remark 8.3.** Let $\hat{\xi}$ be the partition

$$\hat{\xi}(y) = \begin{cases} V^{\mathcal{F},u}(y, \rho) & y \in S_\rho \\ M \smallsetminus S_\rho & y \notin S_\rho \end{cases}$$

and $\xi^{\mathcal{F}^u} := \hat{\xi}^+$. Choosing from a full measure set of $\rho$, we similarly obtain an increasing generator $\xi^{\mathcal{F}^u}$ subordinate to $\mathcal{F}^u$ satisfying the analogues of 1–4. Moreover, we have $y \in \xi^{\mathcal{F}^u}(x)$ if and only if for every $m \leq 0$ such that $f^m(x) \in S$, we have $f^m(y) \in S$ and $f^m(y) \in V_{\mathrm{loc},f^m(x),\varepsilon}^{\mathcal{F},u}$. It follows for $x \in S$ that

$$\xi^u(x) \cap V_{\mathrm{loc},x,\varepsilon}^{\mathcal{F},u} = \xi^{\mathcal{F}^u}(x). \tag{21}$$

We have the following version of [LY1, Corollary 5.3] which will follow from Proposition 9.1 and the inequalities in 10.1 below.

**Proposition 8.4.** *Let $\xi^u = \xi^r$ be a partition as above. Then for any measurable partition $\eta$*

$$h_\mu(f \mid \eta) = h_\mu(f, \eta \vee \xi^u) = H_\mu(\eta^+ \vee \xi^u \mid f(\eta^+ \vee \xi^u)).$$

*In particular, $h_\mu(f) = h_\mu(f, \xi^u)$.*

### 8.3. Entropy properties of the partitions $\xi^i$.

From the proof of [LY2, (9.2)] we have

**Claim 8.5.** *Let $\xi^i$ be any partition as in Section 8.2. Then $h_\mu(f, \xi^i) < \infty$.*

*Proof.* Given $\delta > 0$ let

$$A_\delta := \{x : B^i(x, \delta) \subset \xi^i(x)\}.$$

Then $\mu(A_\delta) \to 1$ as $\delta \to 0$. Let $g(x) = -\log \mu_x^{\xi^i}(f^{-1}\xi^i(x))$. Given any $M < h_\mu(f, \xi^i) = H_\mu(f^{-1}\xi^i \mid \xi^i)$ take $\delta > 0$ so that

$$\int_{A_\delta} g \, d\mu \geq M.$$

Write

$$U^i(x, n, \delta) := \bigcap_{\{j : 0 \leq j \leq n, f^j(x) \in A_\delta\}} (f^{-j}\xi^i)(x).$$

Relative to the family of $\varepsilon$-Lyapunov charts $\{\Phi_x\}$, we observe for every $n \geq 0$ that

$$B^i(x, \delta\ell(x)^{-1}e^{(-\lambda_1 - 2\varepsilon)n}) \subset U^i(x, n, \delta).$$

Moreover

$$-\log \mu_x^{\xi^i}\left(B^i(x, \delta\ell(x)^{-1}e^{(-\lambda_1 - 2\varepsilon)n})\right) \geq -\log \mu_x^{\xi^i}\left(U^i(x, n, \delta)\right) \geq \sum_{j=0}^{n}(\mathbb{1}_{A_\delta} \cdot g)(f^j(x)).$$

Then, by the pointwise ergodic theorem, for $\mu$-a.e. $x$ we have

$$\liminf_{n \to \infty} -\frac{1}{n}\log\mu_x^{\xi^i}(B^i(x, \delta\ell(x)^{-1}e^{(-\lambda_1 - 2\varepsilon)n}) \geq \int_{A_\delta} g \geq M.$$

On the other hand, (see [LY1, Lemma 4.1.4]) fixing a bi-Lipschitz identification of the embedded manifold $V_{\text{loc}, x, \varepsilon}^i$ with $\mathbb{R}^{m_i + \cdots + m_1}$ we have

$$\limsup_{n \to \infty} -\frac{1}{n}\log\mu^i(B^i(x, \delta\ell(x)^{-1}e^{(-\lambda_1 - 2\varepsilon)n})) \leq (m_i + m_{i-1} + \cdots + m_1)(\lambda_1 + 2\varepsilon)$$

whence $M \leq (m_i + m_{i-1} + \cdots + m_1)(\lambda_1 + 2\varepsilon)$. $\qquad\square$

We have the following observation of which we make repeated use.

**Claim 8.6.** *Let $\xi$ be any measurable partition of $(M, \mu)$ satisfying properties 2 and 3 of Section 8.2. Let $\eta$ and $\zeta$ be measurable partitions with $h_\mu(f, \zeta) < \infty$. Then*

$$h_\mu(f, \xi \vee \eta \vee \zeta) \leq h_\mu(f, \xi \vee \eta).$$

Indeed, this follows from (10) of Section 6.1 as $\bigvee_{n \in \mathbb{Z}} f^n(\xi \vee \eta)$ is the point partition. Combined with (6) of Section 6.1, we have (compare to [LY1, Lemma 3.2.1])

**Corollary 8.7.** *Let $\xi$ be any measurable partition of $(M, \mu)$ satisfying 2 and 3 above. Let $\eta$ and $\mathcal{P}$ be measurable partitions with $H_\mu(\mathcal{P}) < \infty$. Then $h_\mu(f, \xi \vee \eta \vee \mathcal{P}) = h_\mu(f, \xi \vee \eta)$.*

As in [LY1, Lemma 3.1.2],

**Lemma 8.8.** *For each $1 \leq i \leq r$ let $\xi_1^i$ and $\xi_2^i$ be two partitions as in 8.2 and let $\eta$ be an arbitrary measurable partition. Then*

*(1) $h_\mu(f, \xi_1^i) = h_\mu(f, \xi_2^i)$.*
*(2) $h_\mu(f, \xi_1^i \vee \eta) = h_\mu(f, \xi_2^i \vee \eta)$.*

8.4. **Finite entropy partitions adapted to Lyapunov charts.** Recall our family of $\varepsilon$-Lyapunov charts $\{\Phi_x\}$. Let $0 < \delta < 1$ be a reduction factor. For $x \in \Lambda$, define the corresponding *center-unstable sets*:

$$S^{cu}_{\delta,x} = \{y \in \mathbb{R}^k : \|\Phi^{-1}_{f^{-m}(x)} \circ f^{-m} \circ \Phi_x(y)\| < \delta : \text{ for all } m \geq 0\}.$$

**Definition 8.9.** We say a measurable partition $\mathcal{P}$ of $(M, \mu)$ is adapted to $(\{\Phi_x\}, \delta)$ if, for almost every $x$,

$$\mathcal{P}^+(x) \subset \Phi_x(S^{cu}_{\delta,x}).$$

**Lemma 8.10.** *For every $0 < \delta < 1$ there is a measurable partition $\mathcal{P}$ subordinate to $(\{\Phi_x\}, \delta)$ with $H_\mu(\mathcal{P}) < \infty$.*

*Proof.* Let $\ell \colon \Lambda \to [1, \infty)$ be the function associated with the charts $\{\Phi_x\}$. Recall in Section 8.2 we fixed and open set $N$ with $\overline{N} \subset U$ to be a compact, $k$-dimensional submanifold with boundary and $S \subset N$ with $\mu(S) > 0$ and $\ell(x) \leq \ell_0$ for $x \in S$.

Let $n \colon S \to \mathbb{N}$ be the first return function: $n(x) = \inf\{j \geq 1 : f^j(x) \in S\}$. We have $\int n(x) \, d\mu\!\upharpoonright_S = 1$.

Let $\rho \colon S \to (0,1)$ be

$$\rho(x) = \delta \ell_0^{-1} e^{-n(x)(\lambda_0 + 2\varepsilon)}.$$

Let $\hat\mu = \frac{1}{\mu(S)} \mu\!\upharpoonright_S$. We have $\int -\log(\rho) \, d\hat\mu < \infty$ hence, adapting [Mañ, Lemma 2] to $\overline{N}$, there is a partition $\hat{\mathcal{P}}$ of $(S, \hat\mu)$ with $H_{\hat\mu}(\hat{\mathcal{P}}) < \infty$ and $\mathrm{diam}(\hat{\mathcal{P}}(x)) < \rho(x)$ for almost every $x \in S$. Let $\mathcal{P} = \hat{\mathcal{P}} \cup \{M \smallsetminus S\}$. Then $H_\mu(\mathcal{P}) < \infty$. Moreover, from the choice of $\rho$ and the properties of the charts $\{\Phi_x, x \in S\}$, the same computations as in [LY1, Lemma 2.4.2] show that $\mathcal{P}$ is adapted to $(\{\Phi_x\}, \delta)$. □

**Claim 8.11.** *Let $\delta < e^{-2\lambda_0 - 4\varepsilon}$ and let $\mathcal{P}$ be adapted to $(\{\Phi_x\}, \delta)$. Then for almost every $x$,*

$$(\mathcal{P}^+ \vee \xi^u)(x) \subset V^u_{\mathrm{loc},x,\varepsilon}.$$

*Proof.* Let $\eta = \mathcal{P}^+ \vee \xi^u$. Note that $f^{-k}(\eta(x)) \subset \eta(f^{-k}(x))$. For $x \in S$ we have

$$\eta(x) \subset \xi^u(x) \subset V^u_{\mathrm{loc},x,\varepsilon}$$

whence $\Phi^{-1}_x(\eta(x)) \subset W^u_{x,\delta}$. By the choice of $\delta$, it follows that if $\Phi^{-1}_x(y) \in W^u_{x,\delta}$ and $f(y) \in \eta(f(x))$ then

$$\Phi^{-1}_{f(x)}(f(y)) \in W^u_{f(x), e^{\lambda_0 + 2\varepsilon}\delta} \cap \mathbb{R}^k(\delta) = W^u_{f(x),\delta}$$

whence $f(y) \in V^u_{\mathrm{loc},f(x),\varepsilon}$. For $x \notin S$ let $n(x) = \max\{k \geq 0 \, f^{-k}(x) \in S\}$. Since $f^{-k}(\eta(x)) \subset \eta(f^{-k}(x))$, proceeding inductively, the claim follows for all $x$. □

## 9. LOCAL ENTROPIES ALONG UNSTABLE FOLIATIONS

Let $\eta$ be a measurable partition of $(M, \mu)$. We do not assume $h_\mu(f \mid \eta) < \infty$. Given $1 \leq i \leq r$, let $\xi^i$ be any measurable partition subordinated to $\mathscr{W}^i$.

For $1 \leq i \leq r$ define the *$i$th unstable Bowen ball*

$$V^i(x, n, \delta) := \{y : d_{V^i_{\mathrm{loc},f^k(x),\varepsilon}}(f^k(x), f^k(y)) < \delta \text{ for } 0 \leq k \leq n\}.$$

Define for $1 \leq i \leq r$

- $\overline{h}_i(x, \delta, \eta) := \limsup_{n\to\infty} -\frac{1}{n}\log(\mu^{\eta^+ \vee \xi^i}_x V^i(x, n, \delta));$
- $\underline{h}_i(x, \delta, \eta) := \liminf_{n\to\infty} -\frac{1}{n}\log(\mu^{\eta^+ \vee \xi^i}_x V^i(x, n, \delta));$
- $\overline{h}_i(x, \eta) := \lim_{\delta\to 0} \overline{h}_i(x, \delta, \eta);$

- $\underline{h}_i(x,\eta) := \lim_{\delta\to 0}\underline{h}_i(x,\delta,\eta).$

The last two limits exist by monotonicity. It is clear that the definitions are independent of the choice of partitions $\xi^i$ subordinate to $\mathscr{W}^i$.

To unify formulas later, let $\xi^0$ denote the point partition on $(M,\mu)$. Then $\mu_x^{\xi^0\vee\eta^+} = \delta_x$ is atomic and, with the same notations as above, we have

- $\underline{h}_0(x,\eta) = \overline{h}_0(x,\eta) = 0.$

**Proposition 9.1.** *For $1 \le i \le r$, let $\xi_i$ be a partition as in Section 8.2. Then for $\mu$-a.e. $x$,*

$$\underline{h}_i(x,\eta) = \overline{h}_i(x,\eta) =: h_i(x,\eta) = h_\mu(f,\xi^i\vee\eta) = H_\mu(f^{-1}(\eta^+\vee\xi^i)\mid\eta^+\vee\xi^i).$$

*Moreover, $h_\mu(f,\xi^i\vee\eta) < \infty$.*

Note that the finiteness of $h_\mu(f,\xi^i\vee\eta)$ follows from 9.1.1 below and the inequality (see the proof of Claim 8.5 and [LY1, Lemma 4.1.4])

$$\overline{h}_i(x,\eta) < \infty.$$

Write $h_i(\eta) = h_\mu(f,\xi^i\vee\eta)$ for the almost-surely constant value of $h_i(x,\eta)$. Note that $h_i(\eta)$ is independent of the choice of $\varepsilon < \varepsilon_0$ and the family of Lyapunov charts $\{\Phi_x\}$.

### 9.1. **Proof of Proposition 9.1.** We prove Proposition 9.1 in two steps.

9.1.1. *Proof that $\underline{h}_i(x,\eta) \ge H_\mu(f^{-1}(\eta^+\vee\xi^i)\mid\eta^+\vee\xi^i)$.* For a fixed $1 \le i \le r$ and $k > 0$ write

$$\xi^{i,k}(x) := (f^{-k}\xi^i)(x).$$

As $\xi^{i,k}(x)$ contains an open ball in $\xi^i(x)$ for almost every $x$, we have for almost every $x$ that

$$\mu_x^{\xi^{i,k}\vee\eta^+} = \frac{1}{\mu_x^{\xi^i\vee\eta^+}((\xi^{i,k}\vee\eta^+)(x))}\mu_x^{\xi^i\vee\eta^+}\!\upharpoonright_{(\xi^{i,k}\vee\eta^+)(x)}.$$

Then,

$$\frac{\mu_x^{\xi^{i,k}\vee\eta^+}(f^{-1}(\xi^{i,k}\vee\eta^+)(x))}{\mu_x^{\xi^{i,k}\vee\eta^+}(f^{-1}\xi^{i,k}(x))} = \frac{\mu_x^{\xi^i\vee\eta^+}(f^{-1}(\xi^{i,k}\vee\eta^+)(x))}{\mu_x^{\xi^i\vee\eta^+}(f^{-1}\xi^{i,k}(x))}$$

$$= \mathbb{E}_{\mu_x^{\xi^i\vee\eta^+}}\left(\mathbb{1}_{f^{-1}\eta^+(x)}\mid\sigma(f^{-1}\xi^{i,k})\right)(x).$$

Given $b > 0$ set

$$A_{b,k} := \left\{x : \frac{\mu_x^{\xi^{i,n}\vee\eta^+}(f^{-1}(\xi^{i,n}\vee\eta^+)(x))}{\mu_x^{\xi^{i,n}\vee\eta^+}(f^{-1}\xi^{i,n}(x))} \ge e^{-b} \text{ for all } n \ge k\right\}.$$

By pointwise convergence for martingales, we have $\mu(A_{b,k}) \to 1$ as $k \to \infty$.

Given $b > 0, \delta > 0$ and $k \in \mathbb{N}$ set

$$A_{b,k,\delta} := \{x \in A_{b,k} : B^i(x,\delta) \subset \xi^{i,k}(x)\}.$$

We have $\mu(A_{b,k,\delta}) \to \mu(A_{b,k})$ as $\delta \to 0$. Given any

$$M < H_\mu(f^{-1}(\eta^+\vee\xi^{i,k})\mid\eta^+\vee\xi^{i,k}) = H_\mu(f^{-1}(\eta^+\vee\xi^i)\mid\eta^+\vee\xi^i)$$

and

$$0 < b < H_\mu(f^{-1}(\eta^+\vee\xi^{i,k})\mid\eta^+\vee\xi^{i,k}) - M$$

choose first $k$ then $\delta_0$ so that for any $0 < \delta \le \delta_0$

$$M + b \le \int_{A_{b,k,\delta}} -\log(\mu_x^{\eta^+\vee\xi^{i,k}}(f^{-1}(\eta^+\vee\xi^{i,k})(x)))\,d\mu(x).$$

Then
$$M \leq \int_{A_{b,k,\delta}} -\log(\mu_x^{\eta^+ \vee \xi^{i,k}}((f^{-1}\xi^{i,k})(x))) \; d\mu(x).$$

Define
$$U_i(x,n,\delta) := \bigcap_{0 \leq j \leq n : f^j(x) \in A_{b,k,\delta}} (f^{-j}\xi^{i,k})(x).$$

Then $V^i(x,n,\delta) \subset U_i(x,n,\delta)$ and with
$$g(x) := -\log \mu^{\eta^+ \vee \xi^i}(f^{-1}\xi^{i,k})(x)$$

we have
$$-\log \mu^{\eta^+ \vee \xi^{i,k}}(U_i(x,n,\delta)) \geq \sum_{j=0}^{n}(\mathbb{1}_{A_{b,k,\delta}} \cdot g)(f^j(x)).$$

Then for a.e. $x$,
$$\begin{aligned}
\underline{h}_i(x,\delta,\eta) &:= \liminf_{n \to \infty} -\frac{1}{n}\log \mu_x^{\eta^+ \vee \xi^i} V^i(x,n,\delta) \\
&\geq \liminf_{n \to \infty} -\frac{1}{n}\log \mu^{\eta^+ \vee \xi^i}(U_i(x,n,\delta)) \\
&\geq \liminf_{n \to \infty} -\frac{1}{n}\sum_{j=0}^{n}(\mathbb{1}_{A_{b,k,\delta}} \cdot g)(f^j(x)) \\
&= \int_{A_{b,k,\delta}} g \; d\mu \geq M
\end{aligned}$$

and the inequality follows.                                                    $\square$

9.1.2. *Proof that $\overline{h}_i(x,\eta) \leq H_\mu(f^{-1}(\eta^+ \vee \xi^i) \mid \eta^+ \vee \xi^i)$.* Given a partition $\eta$ of $(M, \mu)$ and $\ell < k$ define
$$\eta_\ell^k := \bigvee_{j=\ell}^{k-1} f^{-j}\eta \qquad \eta_0^0 = \{\varnothing, M\}.$$

We have exactly as in [LY2, Lemma 9.3.1]

**Lemma 9.2.** *Let $\zeta$ be a partition of $M$ with $H_\mu(\zeta \mid \xi^i \vee \eta^+) < \infty$. Then for $\mu$-a.e. $x$*
$$\lim_{n \to \infty} -\frac{1}{n}\log \mu_x^{\xi^i \vee \eta^+}\left((\zeta \vee \xi^i \vee \eta)_0^n(x)\right) = H_\mu(\xi^i \vee \eta \mid f(\xi_i \vee \eta^+)).$$

Recall the unstable manifolds $W_{x,\delta}^i = \operatorname{graph}(h_x^i|_{\mathbb{R}_{j \leq i}(\delta)})$ defined inside charts the $\Phi_x$ in Section 8.1.

**Lemma 9.3.** *Given $0 < \varepsilon'$ sufficiently small, there exists a partition $\mathcal{P}$ of $(M, \mu)$ with $H_\mu(\mathcal{P}) < \infty$ and a measurable $n_0 \colon M \to \mathbb{N}$ such that for almost every $x \in M$*
$$\Phi_x^{-1}\left((\mathcal{P}_0^n \vee \xi^i)(x)\right) \subset W_{x,\varepsilon'e^{-(\lambda_i - 2\varepsilon)n}}^i.$$

*In particular, for $n \geq n_0(x)$,*
  (a) $(\mathcal{P}_0^n \vee \xi^i)(x) \subset V^i(x,n,\varepsilon')$;
  (b) $(\mathcal{P}_0^n \vee \xi^i)(x) \subset B^i(x, e^{-(\lambda_i - 2\varepsilon)n})$.

Note the final two conclusions follow as $\Phi_x$ is 1-Lipschitz for almost every $x$. The proof is nearly identical to that of [LY2, Lemmas 9.3.3, 9.3.2]. We included it for completeness.

*Proof.* Let $S$ be as in Section 8.2. Recall $\ell(x) \leq \ell_0$ for $x \in S$, $\mu(S) > 0$, and there is an open manifold $N$ with $\overline{N}$ a compact manifold with boundary and $S \subset N$.

Given $x \in S$ define
$$n_+(x) = \min\{n > 0 : f^n(x) \in S\}, \qquad n_-(x) = \min\{n > 0 : f^{-n}(x) \in S\}.$$
Let $\delta = e^{-2\lambda_0 - 4\varepsilon}$ and consider any $0 < \varepsilon' < \delta$. Let $\psi, \psi_+ : S \to (0, 1)$ be
$$\psi(x) = \varepsilon' \ell_0^{-1} e^{-(\lambda_0 + 2\varepsilon)\max\{n_+(x), n_-(x)\}}, \quad \psi_+(x) = \varepsilon' \ell_0^{-1} e^{-(\lambda_0 + 2\varepsilon)n_+(x)}.$$
Let $\mu_S = \frac{1}{\mu(S)}\mu{\restriction}_S$. Then $\int \log(\psi)\, d\mu_S < \infty$. Adapting [Mañ, Lemma 2] (to $(\overline{N}, \mu_S)$), we may find a measurable partition $\mathcal{P}'$ of $S$ with

(1) $H_{\mu_S}(\mathcal{P}') < \infty$;
(2) $\mathcal{P}'(x) \subset B(x, \psi(x))$ for all $x \in S$.

Let $\mathcal{P} = \{S, M \smallsetminus S\} \vee \mathcal{P}'$. Then still $H_\mu(\mathcal{P}) < \infty$ and $\mathcal{P}(x) = \mathcal{P}'(x) \subset B(x, \psi(x))$ for all $x \in S$.

From the dynamics inside Lyapunov charts we have, exactly as in [LY2, Lemma 9.3.2(1)], that for $x \in S$, if $\Phi_x^{-1}(y) \in W_{x,\delta}^i$ and $d(x, y) \le \psi_+(x)$ then
$$\|\Phi_{f^j(x)}^{-1}(f^j(y))\| < \varepsilon', \quad \text{for } 0 \le j \le n_+(x)$$
and
$$\Phi_{f^{n_+(x)}(x)}^{-1}(f^{n_+(x)}(y)) \in W_{f^{n_+(x)}(x),\delta}^i.$$
For $x \in M$ define
$$n_0(x) := \min\{n \ge 0 : f^n(x) \in S\}, \quad r_0(x) := \max\{n < 0 : f^n(x) \in S\}.$$
Consider $x$ with $f^n(x) \in S$ infinitely often as $n \to \pm\infty$ and $y \in \mathcal{P}_0^{n_0(x)}(x) \cap \xi^i(x)$. As $f^{n_0(x)}(y) \in \mathcal{P}(f^{n_0(x)}(x))$ it follows from the choice of $\psi(x)$ that for $0 \le j \le n_0(x) - r_0(x)$,
$$\Phi_{f^{n_0(x)-j}(x)}^{-1}(f^{n_0(x)-j}(\mathcal{P}(f^{n_0(x)}(x)))) \subset \mathbb{R}^k(\delta). \tag{22}$$
Since $\xi^i$ is increasing we have $f^{r_0(x)}(y) \in \xi^i(f^{r_0(x)}(x))$ and since $f^{r_0(x)}(x) \in S$ we have
$$\Phi_{f^{r_0(x)}(x)}^{-1}(f^{r_0(x)}(y)) \in W_{f^{r_0(x)}(x),\delta}^i.$$
It follows from the choice of $\delta$, (22), and the dynamics inside Lyapunov charts that for $0 \le j \le n_0(x) - r_0(x)$,
$$\Phi_{f^{r_0(x)+j}(x)}^{-1}(f^{r_0(x)+j}(y)) \in W_{f^{r_0(x)+j}(x),\delta}^i.$$
In particular,
$$\Phi_{f^{n_0(x)}(x)}^{-1}(f^{n_0(x)}(y)) \in W_{f^{n_0(x)}(x),\delta}^i$$
and
$$\|\Phi_{f^{n_0(x)}(x)}^{-1}(f^{n_0(x)}(y))\| < \varepsilon' e^{-(\lambda_0 + 2\varepsilon)\max\{n_+(x), n_-(x)\}}.$$
Now, let $n_k(x)$ be the subsequent returns of $x$ to $S$. For $n \ge n_0(x)$ take $k$ with $n_k(x) \le n < n_{k+1}(x)$. If $y \in (\mathcal{P}_0^n \vee \xi^i)(x)$ then $f^{n_i(x)}(y) \in \mathcal{P}(f^{n_i(x)}(x))$. For all $0 \le i \le k$ we recursively verify as in the above that
$$\Phi_{f^{n_k(x)}(x)}^{-1} f^{n_k(x)}(y) \in W_{f^{n_k(x)}(x),\delta}^i$$
and
$$\|\Phi_{f^{n_k(x)}(x)}^{-1} f^{n_k(x)}(y)\| \le \varepsilon' e^{-(\lambda_0 + 2\varepsilon)\max\{n_+(f^{n_k(x)}(x)), n_-(f^{n_k(x)}(x))\}}.$$

Applying the dynamics inside Lyapunov charts we obtain

$$\|\Phi_x^{-1}(y)\| \leq e^{-(\lambda_i - 2\varepsilon)n_k(x)}\varepsilon' e^{-(\lambda_0 + 2\varepsilon)\max\{n_+(f^{n_k(x)}), n_-(f^{n_k(x)})\}}$$
$$\leq e^{-(\lambda_i - 2\varepsilon)n}\varepsilon'. \qquad \square$$

*Proof that* $\overline{h}_i(x, \eta) \leq H_\mu(f^{-1}(\eta^+ \vee \xi^i) \mid \eta^+ \vee \xi^i)$. Given $\varepsilon' > 0$ sufficiently small, from Lemmas 9.2 and 9.3 we have

$$\overline{h}_i(x, \varepsilon', \eta) := \limsup_{n \to \infty} -\frac{1}{n} \log(\mu_x^{\eta^+ \vee \xi^i} V^i(x, n, \varepsilon'))$$
$$\leq \limsup_{n \to \infty} -\frac{1}{n} \log(\mu_x^{\eta^+ \vee \xi^i} \mathcal{P}_0^n(x))$$
$$\leq \limsup_{n \to \infty} -\frac{1}{n} \log(\mu_x^{\eta^+ \vee \xi^i}(\eta^+ \vee \xi^i \vee \mathcal{P})_0^n(x))$$
$$= H_\mu(f^{-1}(\eta^+ \vee \xi^i) \mid \eta^+ \vee \xi^i). \qquad \square$$

## 10. Bounds on local entropies

### 10.1. Proof of Propositions 7.1 and 7.4, Claim 7.5, and Theorem 7.7 and Proposition 8.4.
Consider $\eta$ a measurable partition of $(M, \mu)$. As we have not yet proven Proposition 7.1, assume moreover that $h_\mu(f, \eta) < \infty$. Propositions 7.1 and 7.4, Claim 7.5, Theorem 7.7 , and Proposition 8.4 then follow directly from the following four inequalities whose proofs occupy the rest of this and the following section.

Recall the local entropies $h_i(\eta)$ defined in Section 9. We claim that for $1 \leq i \leq r$ the following inequalities hold:

(I) $h_i(\eta) - h_{i-1}(\eta) \geq \lambda_i \left(\overline{\dim}^i(\mu|\eta^+) - \overline{\dim}^{i-1}(\mu|\eta^+)\right)$

(II) $h_i(\eta) - h_{i-1}(\eta) \leq \lambda_i \left(\underline{\dim}^i(\mu|\eta^+) - \underline{\dim}^{i-1}(\mu|\eta^+)\right)$

(III) $h_i(\eta) - h_{i-1}(\eta) \leq \lambda_i m_i$

(IV) $h_r(\eta) = h_\mu(f \mid \eta)$.

### 10.2. Proof of (I).
This is identical to [LY2, (10.2)]; we include it for completeness. Recall our fixed $0 < \varepsilon < \varepsilon_0$ and family of $\varepsilon$-Lyapunov charts $\{\Phi_x\}$ with corresponding function $\ell$. For $0 \leq i \leq r$, let $\xi^i$ be a measurable partitions as in Section 8.2.

**Lemma 10.1.** *For each $1 \leq i \leq r$ there exists a partition $\mathcal{P}$ with $H(\mathcal{P}) < \infty$ and a measurable function $n_0 \colon M \to \mathbb{N}$ such that for $\mu$-a.e. $x$ the following hold for all $n \geq n_0(x)$:*

(a) $\dfrac{\log \mu_x^{\xi^{i-1} \vee \eta^+} B^{i-1}\left(x, e^{-n(\lambda_i - 2\varepsilon)}\right)}{-n(\lambda_i - 2\varepsilon)} \leq \overline{\dim}^{i-1}(\mu \mid \eta^+) + \varepsilon;$

(b) $-\dfrac{1}{n} \log \mu_x^{\xi^{i-1} \vee \eta^+} \mathcal{P}_0^n(x) \geq h_{i-1}(\eta) - \varepsilon;$

(c) $\xi^i(x) \cap \mathcal{P}_0^n(x) \subset B^i\left(x, e^{-n(\lambda_i - 2\varepsilon)}\right);$

(d) $-\dfrac{1}{n} \log \mu_x^{\xi^i \vee \eta^+} \mathcal{P}_0^n(x) \leq h_i(\eta) + \varepsilon;$

(e) $B^{i-1}\left(x, e^{-n(\lambda_i - 2\varepsilon)}\right) \subset \xi^{i-1}(x).$

*Moreover, for infinitely many $n \geq n_0(x)$*

(f) $\dfrac{\log \mu_x^{\xi^i \vee \eta^+} B^i\left(x, 2e^{-n(\lambda_i - 2\varepsilon)}\right)}{-n(\lambda_i - 2\varepsilon)} \geq \overline{\dim}^i(\mu \mid \eta^+) - \varepsilon.$

*Proof.* (a), (e) and (f) follow from definition. (b) and (c) follows from Lemma 9.3 taking sufficiently small $\varepsilon' > 0$. (d) follows from Lemma 9.2. Note also that (a), (b), (d) and (e) hold trivially when $i - 1 = 0$. □

We now prove the first inequality of Section 10.1.

*Proof of (I).* We retain all notation from Lemma 10.1.

With $\Gamma := \{x : n_0(x) \le N_1\}$ select $N_1$ large enough so that for some $x \in \Gamma$ satisfying Lemma 10.1 we have for all $n \ge N_1$

$$\mu_x^{\xi^{i-1} \vee \eta^+}\left(\Gamma \cap B^{i-1}(x, e^{-n(\lambda_i - 2\varepsilon)})\right) \ge \frac{1}{2} \mu_x^{\xi^{i-1} \vee \eta^+}\left(B^{i-1}(x, e^{-n(\lambda_i - 2\varepsilon)})\right)$$
$$\ge \frac{1}{2} e^{(-n(\lambda_i - 2\varepsilon))(\overline{\dim}^{i-1}(\mu|\eta^+) + \varepsilon)}.$$

Fix such an $x$.

For $n \ge N_1$, take $L_n := B^{i-1}(x, e^{-n(\lambda_i - 2\varepsilon)})$. By Lemma 10.1(e), for $y \in \Gamma \cap L_n \cap \eta^+(x)$ and $n \ge N_1$ we have $B^{i-1}\left(y, e^{-n(\lambda_i - 2\varepsilon)}\right) \subset \xi^{i-1}(x)$ and hence from Lemma 10.1(b),

$$\mu_x^{\xi^{i-1} \vee \eta^+} \mathcal{P}_0^n(y) = \mu_y^{\xi^{i-1} \vee \eta^+} \mathcal{P}_0^n(y) \le e^{-n(h_{i-1}(\eta) - \varepsilon)}.$$

For $n \ge N_1$ we the obtain lower a bound on the cardinality of the number of distinct $\mathcal{P}_0^n \cap \xi^i$-atoms meeting $\Gamma \cap L_n \cap \eta^+(x)$ by

$$\#\{(\mathcal{P}_0^n \cap \xi^i)(y) : y \in \Gamma \cap L_n \cap \eta^+(x)\} \ge \mu_x^{\xi^{i-1} \vee \eta^+}\left(\Gamma \cap L_n\right)/e^{-n(h_{i-1}(\eta) - \varepsilon)}$$
$$\ge \frac{1}{2} e^{(-n(\lambda_i - 2\varepsilon))(\overline{\dim}^{i-1}(\mu|\eta^+) + \varepsilon)} e^{n(h_{i-1}(\eta) - \varepsilon)}.$$

For $y \in \Gamma \cap L_n \cap \eta^+(x)$, we have by Lemma 10.1(c) that $(\mathcal{P}_0^n \cap \xi^i)(y) \subset B^i\left(y, e^{-n(\lambda_i - 2\varepsilon)}\right)$ whence $(\mathcal{P}_0^n \cap \xi^i)(y) \subset B^i\left(x, 2e^{-n(\lambda_i - 2\varepsilon)}\right)$. From Lemma 10.1(d), we have

$$\mu_y^{\xi^i \vee \eta^+} \mathcal{P}_0^n(y) \ge e^{-n(h_i(\eta) + \varepsilon)}$$

and obtain inequalities

$$\mu_x^{\xi^i \vee \eta^+} B^i\left(x, 2e^{-n(\lambda_i - 2\varepsilon)}\right)$$
$$\ge \#\{(\mathcal{P}_0^n \cap \xi^i)(y) : y \in \Gamma \cap L_n \cap \eta^+(x)\} \cdot e^{-n(h_i(\eta) + \varepsilon)}$$
$$\ge \frac{1}{2} e^{(-n(\lambda_i - 2\varepsilon))(\overline{\dim}^{i-1}(\mu|\eta^+) + \varepsilon)} e^{n(h_{i-1}(\eta) - \varepsilon)} e^{-n(h_i(\eta) + \varepsilon)}.$$

Comparing to Lemma 10.1(f) we have for infinitely many $n$ that

$$(\lambda_i - 2\varepsilon)(\overline{\dim}^i(\mu \mid \eta^+) - \varepsilon) \le \frac{\log 2}{n} + (\lambda_i - 2\varepsilon)(\overline{\dim}^{i-1}(\mu \mid \eta^+) + \varepsilon) + h_i(\eta) - h_{i-1}(\eta) + 2\varepsilon.$$

Choosing $n$ sufficiently large we have

$$h_i(\eta) - h_{i-1}(\eta) \ge (\lambda_i - 2\varepsilon)\left(\overline{\dim}^i(\mu \mid \eta^+) - \overline{\dim}^{i-1}(\mu \mid \eta^+) - 2\varepsilon\right) - 3\varepsilon.$$

Inequality (I) follows by the arbitrariness of $\varepsilon < \varepsilon_0$. □

## 11. Proofs of (II), (III), and (IV)

11.1. **Properties of fake unstable manifolds inside charts.** Recall we fix $0 < \varepsilon \le \varepsilon_0$, a family of $\varepsilon$-Lyapunov charts $\{\Phi_x\}$, and for $x \in \Lambda_0$ define $\tilde{f}_x \colon \mathbb{R}^k(e^{-\lambda_0 - 2\varepsilon}, \|\cdot\|) \to \mathbb{R}^k(1, \|\cdot\|)$ as in Proposition 5.1. It is convenient to extend each $\tilde{f}_x$ to a function $F_x \colon \mathbb{R}^k \to \mathbb{R}^k$. Taking a $C^\infty$ function $\Theta \colon \mathbb{R}^k \to [0, 1]$ with $\Theta(v) = 0$ for $\|v\| \ge 1$ and $\Theta(v) = 1$ for

$\|v\| \le \frac{1}{2}$ and $0 < \hat{\delta} < 1$ define $F_x : \mathbb{R}^k \to \mathbb{R}^k$ by

$$F_x(v) = D_0 \tilde{f}_x(v) + \Theta(\hat{\delta}^{-1}v)(\tilde{f}_x(v) - D_0\tilde{f}_x(v)). \tag{23}$$

If $\hat{\delta} < e^{-\lambda_0 - 2\varepsilon}$ then $F_x$ is well-defined. Moreover, if $\hat{\delta}$ is sufficiently small, then there is a $C > 1$ so that $\|F_x - D_0\tilde{f}_x\|_{C^1} < \varepsilon$ and $\mathrm{Höl}^\beta(DF_x) \le C$. Fix such sufficiently small $\hat{\delta}$.

It follows for each $1 \le i \le r+1$ and $z \in \mathbb{R}^k$ that there is a $C^1$ function $g^i_{x,z} \colon \bigoplus_{j \le i} \mathbb{R}^j \to \bigoplus_{j > i} \mathbb{R}^j$ with $\|Dg^i_{x,z}\| \le \frac{1}{3}$ such that, writing $\widetilde{W}^i_x(z)$ for the graph of $g^i_{x,z}$, we have $z \in \widetilde{W}^i_x(z)$, $F_x(\widetilde{W}^i_x(z)) = \widetilde{W}^i_{f(x)}(F_x(z))$, $\widetilde{W}^{i-1}_x(z) \subset \widetilde{W}^i_x(z)$, and if $\hat{z} \in \widetilde{W}^i_x(z)$ then $\widetilde{W}^i_x(z) = \widetilde{W}^i_x(\hat{z})$ and

$$e^{\lambda_i - 2\varepsilon}\|z - \hat{z}\| \le \|F_x(z) - F_x(\hat{z})\| \le e^{\lambda_1 + 2\varepsilon}\|z - \hat{z}\|.$$

Write $\widetilde{W}^i_{x,r}(z) := \{\hat{z} \in W^i_x(z) : \|\hat{z}\| < r\}$. (Note we use $\widetilde{W}$ to denote the "fake" unstable manifolds. These depend on the choice of globalized dynamics $F_x$. In particular, when $i = r + 1$ one only expects "fake" center-unstable manifold through $0$.)

Recall the center-unstable sets $S^{cu}_{x,\delta}$ defined in Section 8.4.

**Claim 11.1.** *For* $1 \le i \le r$ *and* $\delta < \hat{\delta}/4$

(a) $S^{cu}_{2\delta,x} = \{z \in \mathbb{R}^k(2\delta) : F^{-1}_{f^{-n}(x)} \circ \cdots \circ F^{-1}_{f^{-1}(x)}(z) \in \mathbb{R}^k(2\delta) \text{ for all } n \ge 0\};$

(b) $\widetilde{W}^i_{x,2\delta}(0) = W^i_{x,\delta};$

(c) $S^{cu}_{2\delta,x} \subset \widetilde{W}^r_{x,2\delta}(0)$ *if* $\lambda_{r+1} < 0$ *and* $S^{cu}_{2\delta,x} \subset \widetilde{W}^{r+1}_{x,2\delta}(0)$ *if* $\lambda_{r+1} = 0$.

*Moreover, for* $z \in S^{cu}_{\delta,x}$

(d) $\widetilde{W}^i_{x,2\delta}(z) \subset S^{cu}_{2\delta,x};$

(e) $\widetilde{W}^i_{x,2\delta}(z) = \tilde{f}_{f^{-1}(x)}(\widetilde{W}^i_{f^{-1}(x),2\delta}(\tilde{f}^{-1}_x(z))) \cap \mathbb{R}^k(2\delta);$

(f) *if* $\Phi_x(z) \in \Lambda_0$ *then*

$$\Phi^{-1}_x(V^i_{\mathrm{loc},\Phi_x(z),\varepsilon}) \cap S^{cu}_{x,\delta} \subset \widetilde{W}^i_{x,2\delta}(z);$$

(g) *if* $\Phi_x(z) \in \Lambda_0$ *and* $\ell(\Phi_x(z)) \le e^{-2\lambda_0 - 4\varepsilon}(4\delta)^{-1}$ *then*

$$\Phi_x(\widetilde{W}^i_{x,2\delta}(z)) \subset V^i_{\mathrm{loc},\Phi_x(z),\varepsilon}.$$

## 11.2. Construction of auxiliary partitions.

Fix a family of measurable partitions $\{\xi^i\}_{0 \le i \le r}$ as in Section 8.2. Recall in the construction of $\xi^i$ in Section 8.2 we fix $x_0 \in \Lambda_1$ to be a density point of $\Lambda_1$. Recall $S$ and $\ell_0$ in Section 8.2. With $\delta_0 = \min\{\hat{\delta}/4, e^{-2\lambda_0 - 4\varepsilon}/(4\ell_0)\}$, let $E \subset S$ be a set with $\mu(E) > 0$ and

$$E \subset B(x_0, 1/2\delta_0\ell^{-1}_0).$$

Then for $x \in E$ we have $E \subset \Phi_x(\mathbb{R}^k(\delta_0))$. Let $\mathcal{P}$ be any measurable partition of $(X, \mu)$ adapted to $(\{\Phi_x\}, \delta_0)$, refining $\{S, M \smallsetminus S\}$ and $\{E, M \smallsetminus E\}$, and satisfying $H_\mu(\mathcal{P}) < \infty$.

Take $\eta_* = (\eta \vee \mathcal{P})^+$ and for $0 \le i \le r$ take $\eta_i = \eta_* \vee \xi^i$. Note $\eta_0$ is the point partition. For notational convention, we also write $\eta_{r+1} := \eta_*$. Note that if no exponent is zero, it follows from [Led] that $\eta_* = \eta_{r+1} = \eta_r$. In the presence of zero exponents we may have $\eta_r \ne \eta_{r+1}$. Recall for almost every $x$ that

(1) $\Phi_x(\eta_{r+1}(x)) \subset S^{cu}_{\delta_0,x} \subset \widetilde{W}^{r+1}_{x,2\delta_0}(0);$

(2) $\eta_i(x) \subset V^i_{\mathrm{loc},x,\varepsilon}$ for $1 \le i \le r$.

As in [LY1, LY2, Lemmas 3.3.1–3.3.2, 11.1.2–11.1.3] we have

**Lemma 11.2.** *For each $1 \le i \le r$, almost every $x$, and every $y \in \Lambda_0 \cap \eta_{i+1}(x)$, writing $\hat{y} = \Phi_x^{-1}(y)$,*

    (a) $\Phi_x \widetilde{W}_{x,2\delta_0}^i(\hat{y}) \cap \eta_{i+1}(x) = \eta_i(y)$;

    (b) $f^{-1}(\eta_i(y)) = \eta_i(f^{-1}(y)) \cap f^{-1}(\eta_{i+1}(x))$.

*Proof.* For (a), consider $z \in \Phi_x(\widetilde{W}_{x,2\delta_0}^i(\hat{y})) \cap \eta_{i+1}(x)$. We show $z \in \xi^i(y)$. As $\mathcal{P}$ refines $\{S, M \smallsetminus S\}$ and $z \in \mathcal{P}^+(y)$, it follows that $f^m(z) \in S$ for all $m \le 0$ such that $f^m(y) \in S$. For all $m \le 0$ we have

$$\Phi_{f^m(x)}^{-1}(f^m(z)) \in \widetilde{W}_{f^m(x),2\delta_0}(\Phi_{f^m(x)}^{-1}(f^m(y)))$$

For $m \le 0$ such that $f^m(y) \in S$ it follows from Claim 11.1(g) and the choice of $\delta_0$ that $f^m(z) \in V_{\mathrm{loc},f^m(y),\varepsilon}^i$. By the characterization of $\xi^i(y)$ in Remark 8.2, it follows that $z \in \xi^i(y)$ whence $\Phi_x \widetilde{W}_{x,2\delta_0}^i(\hat{y}) \cap \eta_{i+1}(x) \subset \xi^i(y) \cap \eta_{i+1}(x) = \eta_i(y)$. The reverse inclusion follows as $\eta_i(y) \subset V_{\mathrm{loc},y,\varepsilon}^i \cap \eta_*(x)$ whence

$$\Phi_x^{-1}(\eta_i(y)) \subset \Phi_x^{-1}(V_{\mathrm{loc},y,\varepsilon}^i) \cap S_{\delta_0,x}^{cu} \subset \widetilde{W}_{x,2\delta_0}^i(\Phi_x^{-1}(y))$$

For (b), first note that we have

$$f^{-1}(\eta_i(y)) \subset \eta_i(f^{-1}(y)) \cap f^{-1}(\eta_{i+1}(x))$$

because $f^{-1}(\eta_i(y)) \subset \eta_i(f^{-1}(y))$ and $\eta_i(y) \subset \eta_{i+1}(x)$. For the reverse inequality, from part (a)

$$
\begin{aligned}
\Phi_x^{-1}&\left(f(\eta_i(f^{-1}(y))) \cap (\eta_{i+1}(x))\right)\\
&= \tilde{f}_x\left(\widetilde{W}_{f^{-1}(x),2\delta}^i(\Phi_x^{-1}(f^{-1}(y)) \cap \eta_{i+1}(f^{-1}(x))\right) \cap \Phi_x^{-1}(\eta_{i+1}(x)) \cap \mathbb{R}^k(2\delta)\\
&= \left(\widetilde{W}_{x,2\delta}^i(\Phi_x^{-1}(y))\right) \cap \Phi_x^{-1}(\eta_{i+1}(x))\\
&= \Phi_x^{-1}(\eta_i(y)).
\end{aligned}
$$

$\square$

Note that for $1 \le i \le r$ we have $h_\mu(f,\eta_i) = h_\mu(f, \xi^i \vee \eta) = h_i(\eta)$ and, from Lemma 11.2(b), for $1 \le i \le r$,

$$h_\mu(f,\eta_i) = H_\mu(f^{-1}\eta_* \mid \eta_i) \le H_\mu(f^{-1}\eta_* \mid \eta_*) = h_\mu(f,\eta_*) < \infty.$$

11.3. **Transverse metrics on $\eta_k/\eta_{k-1}$ and transverse dimension.** Recall $x_0 \in S$ and $E \subset B(x_0, 1/2\delta_0\ell_0^{-1})$ defined in the previous section. Note that for each $1 \le i \le r$ we have that $\Phi_x^{-1}(\eta_i(x)) \subset W_{x,\delta_0}^i(0)$ and $\Phi_x^{-1}(\eta_*(x)) \subset S_{\delta_0,x}^{cu}$.

Consider $x \in \Lambda_0$. For $1 \le i \le r + 1$ let $V^i = \bigoplus_{j \ge i} \mathbb{R}^j$. Given $y, z \in \eta_*(x)$ note that $\widetilde{W}_{x,2\delta_0}^{i-1}(\Phi_x^{-1}(y)) \cap V^i$ and $\widetilde{W}_{x,2\delta_0}^{i-1}(\Phi_x^{-1}(z)) \cap V^i$ are singletons. Given $x \in \Lambda_0$ let

$$d_x^i(y,z) = \|\widetilde{W}_{x,2\delta_0}^{i-1}(\Phi_x^{-1}(y)) \cap V^i - \widetilde{W}_{x,2\delta_0}^{i-1}(\Phi_x^{-1}(z)) \cap V^i\|.$$

We have that $d_x^i$ defines a metric on $\eta_{i-1}$-equivalence classes in $\eta_*(x)$. As in [LY1, LY2, Lemmas 2.3.2, 8.3.2]

**Claim 11.3.** *for $z \in \eta_i(x)$*

$$d_{f(x)}^i(f(x), f(z)) \le e^{\lambda_i + 3\varepsilon} d_x^i(x, z).$$

We give an alternative metric that is independent of the choice of $x$. For each $1 \le i \le r + 1$ let $T^i \subset B(x_0, 1/2\delta_0\ell_0^{-1})$ be a $(\dim(\oplus_{j \ge i} E_i(x)))$-dimensional embedded disc that is uniformly transverse to each $V_{\mathrm{loc},y,\varepsilon}^{i-1}$ for $y \in E$ and such that for $1 \le i \le r$, $V_{\mathrm{loc},y,\varepsilon}^i \cap T^i$ is an embedded $m_i$-dimensional submanifold for each $y \in E$. For $x \in E$ and $y, z \in \eta_*(x)$ define $d^{T_i,i}(y,z)$ to be

$$d^{T_i,i}(y,z) := d^{T_i}(V_{\mathrm{loc},y,\varepsilon}^{i-1} \cap T^i, V_{\mathrm{loc},z,\varepsilon}^{i-1} \cap T^i)$$

where $d^{T_i}$ is the metric on $T_i$ obtained by a bi-Lipschitz identification of $T_i$ with a subset of $\mathbb{R}^{\dim(\oplus_{j \ge i} E_i)}$. For $x \in E$, $d^{T_i,i}$ defines a metric on $\eta_{i-1}$-equivalence classes in $\eta_*(x)$.

For $x \in E$, let $\widetilde{T}_x^i = \Phi_x^{-1}(T_i)$. We have that $\widetilde{T}_x^i$ and $T^i$ are bi-Lipshitz equivalent with Lipschitz constant uniform over $x \in E$. Moreover, defining

$$d^{\widetilde{T}_i,i}(y,z) := \left\| \left( \widetilde{W}_{x,2\delta_0}^{i-1}(\Phi_x^{-1}(y)) \cap \widetilde{T}_x^i \right) - \left( \widetilde{W}_{x,2\delta_0}^{i-1}(\Phi_x^{-1}(z)) \cap \widetilde{T}_x^i \right) \right\|$$

it follows from the main result of [Bro] that, restricted to $\eta_i(x) \subset \Phi_x(\widetilde{W}_{x,\delta_0}^i(0))$, the metrics $d_x^i$ and $d^{\widetilde{T}_i,i}$ are uniformly Lipschitz equivalent. Using that the charts $\{\Phi_x : x \in E\}$ are uniformly Lipschitz embeddings, it follows

**Claim 11.4.** *There is a $N > 1$ such that for all $x \in E$ and $y, z \in \eta_i(x)$,*

$$\frac{1}{N} d^{T_i,i}(y,z) \le d_x^i(y,z) \le N d^{T_i,i}(y,z).$$

For arbitrary $x \in \Lambda_0$, let $n(x) := \inf\{n \ge 0 : f^{-n}(x) \in E\}$. Note that $n(x)$ is constant on elements of $\eta_*$. For $y, z \in \eta_*(x)$, define

$$d^i(y,z) := d^{T_i,i}(f^{-n(x)}(y), f^{-n(x)}(z)).$$

Observe that, as $\eta_{i-1}(f(x)) \subset f\eta_{i-1}(x)$, $d^i$ defines a metric on $\eta_{i-1}$-equivalence classes in $\eta_*(x)$ for almost all $x$. In particular, $d^i$ induces a metric on the quotient $\eta_i(x)/\eta_{i-1}(x)$.

For $1 \le i \le r+1$, we write

$$B^{T,i}(x,\rho) := \{y \in \eta_i : d^i(y,x) < \rho\}.$$

Note that if $y \in B^{T,i}(x,\rho)$ then $\eta_{i-1}(y) \in B^{T,i}(x,\rho)$.

For $1 \le i \le r+1$ let

$$\tilde{\gamma}_i(x) = \liminf_{\rho \to 0} \frac{\log \mu_x^{\eta_i}(B^{T,i}(x,\rho))}{\log \rho}. \tag{24}$$

**Claim 11.5.** *For $1 \le i \le r$, $\tilde{\gamma}_i(x) \le m_i$ and $\tilde{\gamma}_{r+1}(x) \le \dim E^0 + \dim E^s$.*

*Proof.* For $1 \le i \le r$, note that $\eta_i(x) \subset V_{\mathrm{loc},x,\varepsilon}^i$. For $x \in E$, we may identify $\eta_i(x)/\eta_{i-1}$ with a subset of $V_{\mathrm{loc},x,\varepsilon}^i \cap T^i$. Fixing bi-Lipschitz charts on $V_{\mathrm{loc},x,\varepsilon}^i \cap T^i$, it follows for almost every $x' \in \eta_i(x)$ ([LY1], Lemma 4.1.4]) that $\tilde{\gamma}_i(x') \le m_i$. The result for $i = r+1$ holds identifying $T^{r+1}$ with a subset of $\mathbb{R}^{\dim E^0 + \dim E^s}$. $\square$

The following claim summarizes the discussion in [LY2, (11.3)–(11.4)].[2]

**Claim 11.6.** *For $1 \le i \le r$, suppose $c \ge 0$ is such that $\tilde{\gamma}_i(x) \ge c$ for almost every $x \in E$. Then*

$$c \le \underline{\dim}^i(\mu|\eta^+) - \underline{\dim}^{i-1}(\mu|\eta^+).$$

*Proof.* Note that for $x \in E$ we have $\xi^i(x) \in \Phi_x(W_x^r, 1) \subset \Phi_x(1)$. Applying the backwards dynamics, for some $m \ge 0$, $f^{-m}(x) \in S$ and $\Phi_{f^{-m}(x)}^{-1}(f^{-m}(\xi^i(x))) \in \mathbb{R}^k(\hat{\delta}/2)$.

---

Then the globalized backwards dynamics (23) coincides with the local backwards dynamics on $\Phi_x^{-1}(f^{-m}(\xi^i(x)))$. Using the local Lipschitz property of holonomies established in [Bro] (for the globalized past dynamics on $\tilde{f}_x^{-m}\Phi_x^{-1}(\xi^i(x))$ and pushing forward under $\tilde{f}^m$) for each $x \in E$ there is a bi-Lipschitz identification of $\xi^i(x)$ with a subset of $\mathbb{R}^{\dim \oplus_{j \leq i} E^j}$ such that for $x' \in \xi^i(x)$, $\xi^{i-1}(x')$ is contained in a horizontal slice $\mathbb{R}^{\dim \oplus_{j < i} E^j} \times \{t'\}$ and if $x'' \notin \xi^{i-1}(x')$, then $\xi^{i-1}(x'')$ and $\xi^{i-1}(x')$ are contained in distinct horizontal slices $\mathbb{R}^{\dim \oplus_{j < i} E^j} \times \{t'\}$ and $\mathbb{R}^{\dim \oplus_{j < i} E^j} \times \{t''\}$.

Under this identification, we may push forward the measures $\mu_x^{\xi^i \vee \eta^+}$ and $\mu_y^{\eta_i}, y \in \xi^i(x)$. Note that $\eta_i$ refines $\xi^i \vee \eta^+(x)$. In particular,

$$\mu_x^{\xi^i \vee \eta^+} = \int \mu_y^{\eta_i} \, d\mu_x^{\xi^i \vee \eta^+}(y).$$

Note also that by definition, for $y \in S$, $\xi^{i-1}(y) \cap \eta_i(y) = \eta_{i-1}(y)$. As all identifications discussed are bi-Lipschitz, the claim follows from [LY2, Lemma 11.3.2] and [LY2, Lemma 11.3.1]. $\qquad \square$

Note that if $\hat{\eta}$ is a measurable partition of $(M, \mu)$ with $\hat{\eta}^+ \geq \eta^+$ then we may similarly define $\hat{\eta}_* = (\hat{\eta} \vee \mathcal{P})^+$ and $\hat{\eta}_i = \hat{\eta}_* \vee \xi^i$. We can then define for $1 \leq i \leq r+1$

$$\hat{\gamma}_i(x) = \liminf_{\rho \to 0} \frac{\log \mu_x^{\hat{\eta}_i}(B^{T,i}(x, \rho))}{\log \rho}. \tag{25}$$

From [LY2, Lemma 11.3.2], it follows that

$$\operatorname{ess\,inf} \hat{\gamma}_i(x) \leq \operatorname{ess\,inf} \tilde{\gamma}_i(x).$$

Claim 11.6 combined with Proposition 11.7 below and inequality (I) then imply Proposition 7.6.

11.4. **Key Proposition.** The following proposition is the analogue of [LY1, LY2, Propositions 5.1, 11.2].

**Proposition 11.7.** *For the family of partitions*

$$\eta_1 \geq \eta_2 \geq \cdots \geq \eta_r \geq \eta_{r+1}$$

*as in Section 11.2, for every $1 \leq i \leq r+1$ with $\lambda_i \geq 0$ and a.e. $x$ we have*

$$(\lambda_i + 3\varepsilon)\tilde{\gamma}_i(x) \geq (1 - \varepsilon)\left(h_\mu(f, \eta_i) - h_\mu(f, \eta_{i-1}) - 2\varepsilon\right)$$

*where $\tilde{\gamma}_i(x)$ is the transverse dimension (25) associated to $\eta_i/\eta_{i-1}$.*

The 3 inequalities (II), (III), and (IV) now follow from Proposition 11.7, Claims 11.6 and 11.5, and the arbitrariness of $\varepsilon > 0$. Indeed, from Claim 8.6, for $1 \leq i \leq r$,

$$h_\mu(f, \eta_i) = h_\mu(f, \xi^i \vee \eta) = h_i(\eta)$$

are independent of the partition $\mathcal{P}$. Then (III) follows from Proposition 11.7 and Claim 11.5. (II) follows from Claim 11.6 with

$$c = \frac{(1 - \varepsilon)\left(h_\mu(f, \eta_i) - h_\mu(f, \eta_{i-1}) - 2\varepsilon\right)}{(\lambda_i + 3\varepsilon)}$$

and Proposition 11.7 applied to $x \in E$.

For (IV), first note that $h_r(\eta) = h_\mu(f, \eta \vee \xi^u) = h_\mu(f, \eta \vee \xi^r \vee \mathcal{P}) \leq h_\mu(f \mid \eta)$ where the first equality holds by Proposition 9.1 and the inequality is from definition. Moreover, given any $M < h_\mu(f \mid \eta)$ we may assume $\mathcal{P}$ in Section 11.2 is chosen so that

$$h_\mu(f, \eta \vee \mathcal{P}) := h_\mu(f, \eta_{r+1}) > M. \tag{26}$$

If $\lambda_{r+1} < 0$ then $\eta_r = \eta_{r+1}$ and

$$M \leq h_\mu(f, \eta_{r+1}) = h_\mu(f, \eta \vee \xi^r \vee \mathcal{P}) = h_\mu(f, \eta \vee \xi^u) \leq h_\mu(f \mid \eta)$$

and the result follows. If $\lambda_{r+1} = 0$ then Then Proposition 11.7 gives

$$3\varepsilon \dim(E^0 \oplus E^s) \geq (1 - \varepsilon)(M - h_r(\eta) - 2\varepsilon).$$

As $\varepsilon < \varepsilon_0$ is arbitrary, we have $M \leq h_r(\eta)$ and (IV) follows.

Note that if $\lambda_1 < 0$ then $\eta_*$ is the point partition and it follows that $h_\mu(f) = 0$. Similarly, if $\lambda_1 = 0$ then $\eta_* = \eta_1$ and $\eta_0$ is the point partition. Proposition 11.7 applies to this setting and shows again that $h_\mu(f) = 0$.

11.5. **Proof of Proposition 11.7.** The proof is identical to [LY1, Proposition 5.1]; we include it for completeness. For $1 \leq i \leq r + 1$ with $\lambda_i \geq 0$, define the following measurable functions $g, g_\delta, g_* : M \to \mathbb{R}$

$$g(x) := \mu_x^{\eta_{i-1}}(f^{-1}\eta_i(x)) = \mu_x^{\eta_{i-1}}(f^{-1}\eta_{i-1}(x))$$

$$g_\delta(x) := \frac{1}{\mu_x^{\eta_i}(B^{i,T}(x, \delta))} \int_{B^{i,T}(x,\delta)} \mu_z^{\eta_{i-1}}((f^{-1}\eta_i)(x)) \, d\mu_x^{\eta_i}(z)$$

$$= \frac{\mu_x^{\eta_i}\left[(B^{i,T}(x, \delta)) \cap (f^{-1}\eta_i)(x)\right]}{\mu_x^{\eta_i}(B^{i,T}(x, \delta))}$$

$$g_*(x) := \inf_{\delta > 0} g_\delta(x)$$

The equality in the definition of $g(x)$ follows from Lemma 11.2.

Recall that for $1 \leq i \leq r + 1$ and almost every $x$, the metric $d^i$ identifies the quotient $(\eta_i(x)/\eta_{i-1})$ with a subset of $\mathbb{R}^{\dim(\oplus_{j \geq i} E_i)}$. Pushing forward the measures $\mu_x^{\eta_i}$ under this identification define $h \colon \eta_i(x) \to \mathbb{R}$ to be $h(z) = g\!\upharpoonright_{\eta_i(x)}$ and define the partition $\zeta$ of $\eta_i(x)$ to be the partition $\zeta = f^{-1}\eta_i\!\upharpoonright_{\eta_i(x)}$. Applying [LY1, Lemma 4.1.3] to $h$ and $\zeta$, we have $g_\delta \to g$ $\mu_x^{\eta_i}$-a.e. and that

$$\int g_*(z) \, d\mu_x^{\eta_i}(z) \leq H_{\mu_x^{\eta_i}}(\zeta) + \log c + 1$$

where $c$ is a constant depending only on the dimension of $\mathbb{R}^{\dim(\oplus_{j \geq i} E_i)}$ coming from the Besicovitch Covering Lemma (see [LY1, (4.1)].) It follows that $g_\delta \to g$ $\mu$-a.e. and that

$$\int g_*(z) \, d\mu(z) \leq \int H_{\mu_x^{\eta_i}}(\zeta) d\mu(x) + \log c + 1 \leq H(f^{-1}\eta_i \mid \eta_i) + \log c + 1 < \infty.$$

Take $E' \subset E$ with $0 < \mu(E') < \frac{\varepsilon}{4 \log N}$. Let $r_0 \leq 0 < r_1 < \ldots$ denote the distinct times when $f^{r_i}(x) \in E'$. For $n$ sufficiently large and $0 \leq k < n$ for $r_j \leq k \leq r_{j+1}$ write

$$a(x; n, k) := B^{T,i}\left(f^k(x), e^{-(\lambda_i + 3\varepsilon)(n - r_j)} N^{2j}\right) \subset \eta_i(f^k(x)).$$

**Claim 11.8.** $a(x; n, k) \cap f^{-1}\eta_i(f^k(x)) \subset f^{-1}(a(x; n, k + 1))$.

Consider $\mu_x^{\eta_i}(a(x; n, 0))$. With $p(n) := \lfloor (1 - \varepsilon)n \rfloor$ we write

$$\mu_x^{\eta_i}(a(x; n, 0)) = \mu_{f^p(x)}^{\eta_i}(a(x; n, p)) \prod_{k=0}^{p-1} \frac{\mu_{f^k(x)}^{\eta_i}(a(x; n, k))}{\mu_{f^{k+1}(x)}^{\eta_i}(a(x; n, k + 1))}$$

$$\leq \prod_{k=0}^{p-1} \frac{\mu_{f^k(x)}^{\eta_i}(a(x; n, k))}{\mu_{f^{k+1}(x)}^{\eta_i}(a(x; n, k + 1))}.$$

Renormalizing and applying Claim 11.8 we have

$$
\frac{\mu_{f^k(x)}^{\eta_i}(a(x;n,k))}{\mu_{f^{k+1}(x)}^{\eta_i}(a(x;n,k+1))} = \mu_{f^k(x)}^{\eta_i}(a(x;n,k)) \frac{\mu_{f^k(x)}^{\eta_i}((f^{-1}\eta_i)(f^k(x)))}{\mu_{f^k(x)}^{\eta_i}(f^{-1}(a(x;n,k+1)))}
$$

$$
\leq \frac{\mu_{f^k(x)}^{\eta_i}(a(x;n,k))}{\mu_{f^k(x)}^{\eta_i}[(f^{-1}\eta_i)(f^k(x)) \cap (a(x;n,k))]} \mu_{f^k(x)}^{\eta_i}((f^{-1}\eta_i)(f^k(x)))
$$

$$
= g_{\delta(x;n,k)}(f^k(x))e^{-I(f^k)}
$$

where $\delta(x;n,k) := e^{-(\lambda_i+3\varepsilon)(n-r_{j_k})}N^{2j_k}$ for $r_{j_k} \leq k \leq r_{j_{k+1}}$ and $I(y) := -\log \mu_y^{\eta_i}((f^{-1}\eta_i)(y))$. Then

$$
(\lambda_i+3\varepsilon)\liminf_{r\to 0} \frac{\log\left(\mu_x^{\eta_i}B^{T,i}(x,r)\right)}{\log r}
$$

$$
= (\lambda_i+3\varepsilon)\liminf_{n\to\infty} \frac{\log\left(\mu_x^{\eta_i}B^{T,i}(x,e^{-(\lambda_i+3\varepsilon)n})\right)}{\log e^{-(\lambda_i+3\varepsilon)n}}
$$

$$
= (\lambda_i+3\varepsilon)\liminf_{n\to\infty} \frac{\log\left(\mu_x^{\eta_i}a(x;n,0)\right)}{\log e^{-(\lambda_i+3\varepsilon)n}}
$$

$$
\geq \liminf_{n\to\infty} \frac{1}{n}\sum_{k=0}^{p(n)}\log g_{\delta(x;n,k)}(f^k(x)) + \lim_{n\to\infty}\frac{1}{n}\sum_{k=0}^{p(n)}I(f^k(x)). \quad (27)
$$

For almost every $x$, $\lim_{n\to\infty}\frac{1}{n}\sum_{k=0}^{p(n)}I(f^k(x))$ exists and

$$
\lim_{n\to\infty}\frac{1}{n}\sum_{k=0}^{p(n)}I(f^k(x)) = (1-\varepsilon)H_\mu(\eta_i \mid f\eta_i) = (1-\varepsilon)h_\mu(f,\eta_i). \quad (28)
$$

Now, let

$$
A_{\delta'} := \{x : -\log g_\delta(x) \leq -\log g(x) + \varepsilon \text{ for all } \delta \leq \delta'\}.
$$

Since $\int -\log g_* < \infty$, taking $\delta'$ sufficiently small we can ensure

$$
\int_{M \smallsetminus A_{\delta'}} -\log g_* < \varepsilon.
$$

By the choice of $\mu(E')$, having taken $n$ sufficiently large we have $J_n = \#\{0 \leq i \leq nf^i(x) \in E'\} \leq \varepsilon n M(E')$. Then, for large enough $n$

$$
\delta(x;n,k) := e^{-(\lambda_i+3\varepsilon)(n-r_{j_k})}N^{2j_k} \leq e^{-\varepsilon(n-p)+2\log NJ_n} \leq e^{(-\varepsilon^2+\varepsilon^2/2)n} \leq \delta'
$$

for all $0 \leq k \leq p(n) = \lfloor(1-\varepsilon)n\rfloor$.

Then for such $n$,

$$
\sum_{k=0}^{p(n)} -\log g_{\delta(x;n,k)}(f^k(x))
$$

$$
\leq \sum_{\{0\leq k\leq p(n):f^k(x)\in A_{\delta'}\}} -\log g_{\delta(x;n,k)}(f^k(x)) + \sum_{\{0\leq k\leq p(n):f^k(x)\notin A_{\delta'}\}} -\log g_{\delta(x;n,k)}(f^k(x))
$$

$$
\leq \sum_{\{0\leq k\leq p(n):f^k(x)\in A_{\delta'}\}} \left(-\log g(f^k(x))+\varepsilon\right) + \sum_{\{0\leq k\leq p(n):f^k(x)\notin A_{\delta'}\}} -\log g_*(f^k(x))
$$

whence

$$\limsup_{n\to\infty} \frac{1}{n} \sum_{k=0}^{p(n)} -\log g_{\delta(x;n,k)}(f^k(x)) \le (1-\varepsilon)\left[\int_{A'_\delta}(-\log g + \varepsilon) + \int_{M\smallsetminus A_{\delta'}} -\log g_*\right]$$

$$\le (1-\varepsilon)\left[\int -\log g\ d\mu + 2\varepsilon\right]$$

$$= (1-\varepsilon)\left[H_\mu(f\in\eta_{i-1}|\eta_{i-1}) + 2\varepsilon\right]$$

$$= (1-\varepsilon)\left[h_\mu(f,\eta_{i-1}) + 2\varepsilon\right]$$

Combining the above with (27) and (28), the proposition follows. $\qquad\square$

## 12. Proof of Theorem 7.2 and Proposition 7.3

Let $\mathcal{F}$ be an $f$-invariant, $C^{1+\beta}$-tame, measurable foliation. Let $\eta$ be a measurable partition of $(M,\mu)$ and let $\xi^{\mathcal{F}^u}$ be as in Remark 8.2. We claim

**Claim 12.1.** $h_\mu(f\mid\mathcal{F}\vee\eta) = h_\mu(f,\eta\vee\xi^{\mathcal{F}^u}) = H_\mu(\eta^+\vee\xi^{\mathcal{F}^u}\mid f(\eta^+\vee\xi^{\mathcal{F}^u}))$.

Note that Proposition 7.3 follows.

*Proof.* Consider $\xi$ a measurable partition subordinate to $\mathcal{F}$. Let $\varepsilon < \varepsilon_0$ be sufficiently small, let $\{\Phi_x\}$ be a family of $\varepsilon$-Lyapunov charts, and let $\{\Psi_x\}$ a family of charts adapted to $\mathcal{F}$ built from $\Phi_x$ as in Proposition 5.4. Let $\xi^u$ and $\xi^{\mathcal{F}^u}$ be as constructed Remark 8.3 and let $\rho$ and $S$ be in the construction. Refining $\xi$, we may assume that $\xi(x) \subset V_{\mathrm{loc},y,\varepsilon}^{\mathcal{F}}$ for $y\in S$.

Given $0 < \delta < e^{-2\lambda_0 - 4\varepsilon}$, let $\mathcal{P}$ be a measurable partition with $H_\mu(\mathcal{P}) < \infty$ adapted to $\{\Phi,\delta\}$. We may assume for any $b > 0$ that $\xi$ and $\mathcal{P}$ are chosen so that

$$h_\mu(f,\mathcal{P}\vee\eta\vee\xi) \ge h_\mu(f\mid\mathcal{F}\vee\eta) - b.$$

From Proposition 8.4, we have

$$h_\mu(f,\mathcal{P}\vee\eta\vee\xi) \le h_\mu(f\mid\mathcal{P}\vee\eta\vee\xi) = h_\mu(f,\mathcal{P}\vee\eta\vee\xi\vee\xi^u).$$

Take $\zeta = (\eta\vee\mathcal{P}\vee\xi)^+\vee\xi^u$.

Note that if $x\in S$ and $y\in\zeta(x)\subset(\xi\cap\xi^u)(x)$ then from the choice of $\xi$ we have

$$y\in V^u(x,\rho)\cap V_{\mathrm{loc},x,\varepsilon}^{\mathcal{F}} = V^{\mathcal{F},u}(x,\rho).$$

By the characterization of $\xi^{\mathcal{F}^u}$ in Remark 8.3, it follows that $\zeta\ge\xi^{\mathcal{F}^u}$.

By definition we have $h_\mu(f\mid\mathcal{F}\vee\eta) \ge h_\mu(f\mid\mathcal{F}^u\vee\eta) \ge h_\mu(f,\xi^{\mathcal{F}^u}\vee\eta)$. On the other hand, with $\zeta$ as above, we have

$$h_\mu(f\mid\mathcal{F}\vee\eta) - b \le h_\mu(f,\mathcal{P}\vee\xi\vee\eta) \le h_\mu(f,\mathcal{P}\vee\xi\vee\eta\vee\xi^u) = h_\mu(f,\zeta).$$

As $\zeta\vee\xi^{\mathcal{F}^u} = \zeta$ and $\xi^u\vee\xi^{\mathcal{F}^u} = \xi^{\mathcal{F}^u}$ we also have

$$h_\mu(f,\zeta) = h_\mu(f,\zeta\vee\xi^{\mathcal{F}^u}) = h_\mu(f,\eta\vee\mathcal{P}\vee\xi\vee\xi^{\mathcal{F}^u}) \le h_\mu(f,\xi^{\mathcal{F}^u}\vee\eta)$$

which combined with the above completes the proof of the claim. $\qquad\square$

We now turn to proof of Theorems 7.2. Let $\xi^{\mathcal{F}^u}$ be an increasing generator subordinate to $\mathcal{F}^u$ as in Remark 8.3. From Claim 12.1, we have that $h_\mu(f\mid\mathcal{F}) = h_\mu(f,\xi^{\mathcal{F}^u})$. Let $\eta = \xi^{\mathcal{F}^u}$. We first claim that for $1\le i\le r$ that

$$\gamma^i(\mu\mid\eta) \le m_i(\mathcal{F}).$$

Indeed, in the notation of Sections 8.2 and 11.3 we may having taken the set $S$ to be such that the functions $\check{h}_x^{\mathcal{F}}$ in Proposition 5.2 vary continuously. Let $\mathcal{F}_{\mathrm{loc},x} = \Phi_x(\mathrm{graph}(\tilde{h}_x^{\mathcal{F}}))$.

We may further assume that the transversals $T^i$ in Section 11.3 were taken to be in general position with $\mathcal{F}_{\mathrm{loc},x}$ for every $x \in E$. As before, let $\eta_* = \xi^{\mathcal{F}^u} \vee \mathcal{P}$ where $\mathcal{P}$ is as in Section 11.2 and also satisfying (26). Then for $x \in E$, we have $\eta_i(x) \subset V^i_{\mathrm{loc},x,\varepsilon} \cap \mathcal{F}_{\mathrm{loc},x} \cap T^i$ and the metric $d^{T_i,i}$ identifies $\eta_i(x)/\eta_{i-1}$ with a subset of a $m_i(\mathcal{F})$-dimensional submanifold of $\mathbb{R}^{\dim(\oplus_{j\geq i} E_i)}$. It follows that $\tilde{\gamma}_i(x) \leq m_i(\mathcal{F})$. Inequality (I), Proposition 11.7, and Claim 11.5 of Section 11 then show that $\gamma^i(\mu | \xi^{\mathcal{F}^u}) \leq m_i(\mathcal{F})$ whence the inequality in Theorem 7.2 then follows from Theorem 7.7.

To complete the proof of Theorems 7.2 we claim that if $h_\mu(f, \xi^{\mathcal{F}^u}) = \sum_{1 \leq i \leq r} \lambda_i m_i(\mathcal{F})$, exactly as in [Led] or [LY1], (6.1), it follows from Jensen's inequality that $\mu_x^{\xi^{\mathcal{F},u}}$ is absolutely continuous along $\mathcal{F}^u$ for almost every $x$. We sketch the details here.

Recall that for $x \in S$ we have $\xi^{\mathcal{F}^u}(x) \subset V^{\mathcal{F},u}_{\mathrm{loc},x,\varepsilon}$. Define $n(x) = \min\{n \geq 0 : f^n(x) \in S\}$ and

$$\overline{\xi}^{\mathcal{F}^u}(x) = f^{-n(x)}(\xi^{\mathcal{F}^u}(f^{n(x)}(x))).$$

We then have

(1) $\overline{\xi}^{\mathcal{F}^u}$ is a partition of $(M, \mu)$ subordinate to $\mathcal{F}^u$;
(2) $\overline{\xi}^{\mathcal{F}^u}(x) \subset V^{\mathcal{F},u}_{\mathrm{loc},x,\varepsilon}$ for almost every $x$;
(3) $f^{-1}\overline{\xi}^{\mathcal{F}^u} \geq \overline{\xi}^{\mathcal{F}^u}$;
(4) $h_\mu(f \mid \mathcal{F}) = h_\mu(f, \overline{\xi}^{\mathcal{F}^u})$.

Replacing $\xi^{\mathcal{F}^u}$ with $\overline{\xi}^{\mathcal{F}^u}$ we assume $\xi^{\mathcal{F}^u}$ satisfies (2) for the remainder. Note that this ensures $D_z f$ is defined for all $z \in \xi^{\mathcal{F}^u}(x)$ and almost every $x$.

For each $C \in \xi^{\mathcal{F}^u}$ and $x \in C$ restrict the ambient Riemannian metric to $V^{\mathcal{F},u}_{\mathrm{loc},x,\varepsilon}$ and consider the induced Riemannian volume $m_x$ to $C$. Note that $m_x = m_y$ if $y \in C$ and since $\xi^{\mathcal{F}^u}$ is subordinated to $\mathcal{F}^u$, $m_x$ is a positive measure for almost every $x$.

For every $x$ define $\Delta_x : \xi^{\mathcal{F}^u}(x) \to [0, \infty)$ by

$$\Delta_x(y) = \lim_{n\to\infty} \frac{\prod_{i=1}^n J^{\mathcal{F},u}(f^{-i}(x))}{\prod_{i=1}^n J^{\mathcal{F},u}(f^{-i}(y))}$$

where for $z \in \xi^{\mathcal{F}^u}(x)$ we define

$$J^{\mathcal{F},u}(z) = \left| \det D_z f \!\restriction_{T_z V^{\mathcal{F},u}_{\mathrm{loc},x,\varepsilon}} \right|$$

and the determinant is computed against the Riemannian metric on $U_0$. Note that as $\xi^{\mathcal{F},u}(x) \subset V^{\mathcal{F},u}_{\mathrm{loc},x,\varepsilon}$ for a.e. $x$, for all $y \in \xi^{\mathcal{F}^u}(x)$ we have $f^{-n}(y) \in U_0$ and $J^{\mathcal{F},u}f^{-n}(y)$ is defined for all $n \geq 0$.

Arguing in $\varepsilon$-charts $\{\Phi_x\}$, we have that as in [LY1], (6.1)

**Claim 12.2.**

(a) for almost every $x$, $\Delta_x$ is uniformly $\beta$-Hölder on $\xi^{\mathcal{F}^u}(x)$;
(b) defining $L(x) := \int_{\xi^{\mathcal{F}^u}(x)} \Delta_x(y) \, dm_x(y)$, for almost every $x$
$$0 < L(x) < \infty;$$
(c) $L(f(x)) = J^{\mathcal{F}^u}(x) \int_{(f^{-1}\xi^{\mathcal{F}^u})(x)} \Delta_x(y) \, d\mu_x(y) \leq J^{\mathcal{F}^u}(x) L(x);$
(d) defining $\nu_x$ on $\xi^{\mathcal{F}^u}(x)$ by
$$d\nu_x(y) = \frac{\Delta_x(y)}{L(x)} dm_x$$

we have that $\nu_x = \nu_y$ if $y \in \xi^{\mathcal{F}^u}(x)$, $\{\nu_x\}$ is a measurable family of probability measures, and

$$\int -\log \nu_x(f^{-1}\xi^{\mathcal{F}^u})(x) \, d\mu(x) = \int \log J^{\mathcal{F}^u}(x) \, d\mu(x).$$

Let $\nu$ be the probability measure defined on $M$ by

$$\nu(A) = \int \nu_x(A) \, d\mu(x)$$

for Borel $A$. As the difference between the Riemannian metric on $U_0$ and the metric in the charts given in Section 3.1 is tempered by a slowly increasing function, we have

$$\int \log \left( J^{\mathcal{F}^u}(x) \right) \, d\mu(x) = \sum m_i(\mathcal{F})\lambda_i.$$

Exactly as in [LY1, Lemma 6.1.3] (see also [Led]), Jensen's inequality gives

**Claim 12.3.** *If* $\int \log J^{\mathcal{F}^u}(x) = H_\mu(f^{-1}(\xi^{\mathcal{F}^u}) \mid \xi^{\mathcal{F}^u})$ *then* $\nu$ *and* $\mu$ *coincide on the* $\sigma$-*algebra generated by* $f^{-1}(\xi^{\mathcal{F}^u})$.

As $\xi^{\mathcal{F}^u}$ generates, it follows that the measures $\nu$ and $\mu$ coincide. Theorem 7.2 follows.

University of Chicago, Chicago, IL 60637, USA
*E-mail address*: awb@uchicago.edu

# SMOOTH ERGODIC THEORY OF $\mathbb{Z}^d$-ACTIONS PART 3: PRODUCT STRUCTURE OF ENTROPY

AARON BROWN, FEDERICO RODRIGUEZ HERTZ, AND ZHIREN WANG

ABSTRACT. For a smooth action of $\mathbb{Z}^d$ preserving a Borel probability measure, we show that entropy satisfies a certain "product structure" along coarse unstable manifolds. Moreover, given two smooth $\mathbb{Z}^d$-actions—one of which is a measurable factor of the other—we show that all coarse Lyapunov exponents contributing to the entropy of the factor system are coarse Lyapunov exponents of the total system and derive an Abramov–Rohlin formula for entropy subordinated to coarse unstable manifolds.

## 13. STATEMENT OF RESULTS

As in Part 1, take $M$ to be a $C^\infty$ manifold equipped with a Borel probability measure $\mu$. Let $\alpha\colon \mathbb{Z}^d \times M \to M$ be an action by measure-preserving, measurable transformations. We moreover assume $(M, \mu)$ and $\alpha$ satisfy the standing hypotheses of Section 3.1. We further assume for simplicity that $\mu$ is ergodic.

### 13.1. **Product structure and subadditivity of entropy.** Our first main result of is the following "product structure of entropy" formula. Recall $\hat{\mathcal{L}}$ denotes the coarse Lyapunov exponents of $\alpha$ with respect to $\mu$ and for $\chi \in \hat{\mathcal{L}}$, $\mathscr{W}^\chi$ is the corresponding foliation by coarse Lyapunov manifolds.

**Theorem 13.1.** *Let $\mathcal{F}$ be an $\alpha$-invariant, $C^{1+\beta}$-tame, measurable foliation and let $\eta$ be an $\alpha$-invariant measurable partition. Then for $n \in \mathbb{Z}^d$,*

$$h_\mu(\alpha(n) \mid \mathcal{F} \vee \eta) = \sum_{\{\chi \in \hat{\mathcal{L}}: \chi(n) > 0\}} h_\mu(\alpha(n) \mid \mathcal{F} \vee \mathscr{W}^\chi \vee \eta).$$

In particular, we have

**Corollary 13.2** (Product structure of entropy)**.**

$$h_\mu(\alpha(n)) = \sum_{\{\chi \in \hat{\mathcal{L}}: \chi(n) > 0\}} h_\mu(\alpha(n) \mid \mathscr{W}^\chi). \tag{29}$$

Note that if $f\colon M \to M$ and $g\colon N \to N$ are diffeomorphisms preserving $\mu_1$ and $\mu_2$, respectively, then there is a natural $\mathbb{Z}^2$-action on $M \times N$ preserving $\mu_1 \times \mu_2$. In this case, (29) follows immediately from the classical Ledrappier–Young entropy formula. Our result (29) suggest that, at least at the level of entropy, an ergodic $\alpha$-invariant measure behaves like a product measure along coarse Lyapunov manifolds. It would be of interest to know if the unstable conditional measures are necessarily products of conditional measures along coarse Lyapunov manifolds. In homogeneous settings considered in [EL] and [EK1, EK2], similar product structures of entropy are established by first establishing a product structure of the measure along coarse manifolds.

As a direct corollary of Theorems 13.1 and 7.7, we recover the subadditivity of entropy of $\mathbb{Z}^d$-actions first obtained in [Hu] for commuting $C^2$ diffeomorphisms of compact manifolds.

**Theorem 13.3** (c.f. [Hu], Theorem B). *Let $\mathcal{F}$ be an $\alpha$-invariant, $C^{1+\beta}$-tame, measurable foliation and let $\eta$ be an $\alpha$-invariant measurable partition. Then for all $n, m \in \mathbb{Z}^d$*

(1) $h_\mu(\alpha(n+m) \mid \mathcal{F} \vee \eta) \le h_\mu(\alpha(n) \mid \mathcal{F} \vee \eta) + h_\mu(\alpha(m) \mid \mathcal{F} \vee \eta);$

(2) *moreover, if $\mathcal{F} \vee \mathscr{W}_n^u = \mathcal{F} \vee \mathscr{W}_m^u$ then equality holds.*

We also obtain the following exact dimensionality formula for measures invariant under $\mathbb{Z}^k$-actions. Let $\mathcal{F}$ be an $\alpha$-invariant, $C^{1+\beta}$-tame, measurable foliation and let $\eta$ be an $\alpha$-invariant measurable partition. For $\chi \in \mathcal{L}$ let $d^{\mathcal{F},\chi,\eta}(\mu)$ be the almost surely constant value of the pointwise dimension of $\mu$ along $\mathcal{F} \vee \mathscr{W}^\chi \vee \eta$ and for $n \in \mathbb{Z}^d$ let $d_n^{u,\mathcal{F},\eta}(\mu)$ be the almost-surely constant value of the pointwise dimension of $\mu$ along $\mathcal{F} \vee \mathscr{W}_n^u \vee \eta$. From Theorems 13.1 and 7.7 we obtain

**Corollary 13.4** (Product structure of unstable dimension). *For any $n \in \mathbb{Z}^d$,*
$$d_n^{u,\mathcal{F},\eta}(\mu) = \sum_{\chi(n)>0} d^{\mathcal{F},\chi,\eta}(\mu).$$

13.2. **Measurable factors and coarse Abramov–Rohlin formula.** Consider a second action $\hat\alpha$ of $\mathbb{Z}^d$ on $(N, \nu)$ satisfying the standing hypotheses of Section 3.1. We say that $\hat\alpha$ is a *measurable factor of* $\alpha$ if there is a measurable map $\psi \colon M \to N$ with $\psi_* \mu = \nu$ and $\psi \circ \alpha(n) = \hat\alpha(n) \circ \psi$ for all $n \in \mathbb{Z}^d$. Let $\mathcal{A}^\psi$ denote the $\alpha$-invariant partition of $(M, \mu)$ into level sets of $\psi$. We assume $\mu$ and thus $\nu$ are ergodic. To distinguish data associated to each action, let $\hat{\mathcal{L}}^{\hat\alpha}(\nu)$ and $\hat{\mathcal{L}}^\alpha(\mu)$ denote, respectively, the coarse Lyapunov exponents for the actions $\hat\alpha$ and $\alpha$ on $(N, \nu)$ and $(M, \mu)$.

Consider a coarse Lyapunov exponent $\hat\chi \in \hat{\mathcal{L}}^{\hat\alpha}(\nu)$ of $\hat\alpha$ and suppose that
$$h_\nu(\hat\alpha(n) \mid \mathscr{W}^{\hat\chi}) > 0 \tag{30}$$
for some $n \in \mathbb{Z}^d$ with $\hat\chi(n) > 0$. Let $E = \ker(\hat\chi) \subset \mathbb{R}^d$ be the *Lyapunov hyperplane* determined by $\hat\chi$. It follows from Corollary 13.2 and (30) that for any open cone $C \subset \mathbb{R}^d$ containing $E$, the function
$$n \mapsto h_\nu(\hat\alpha(n))$$
is not a linear function on $C \cap \mathbb{Z}^d$. By the classical Abramov-Rohlin formula (19), for every $n \in \mathbb{Z}^d$ we have that
$$h_\nu(\hat\alpha(n)) = h_\mu(\alpha(n)) - h_\mu(\alpha(n) \mid \mathcal{A}^\psi).$$

If no Lyapunov exponent of $\alpha$ were proportional to $\hat\chi$ then, taking any open cone $C' \subset C \subset \mathbb{R}^d$ containing $E$ and disjoint from the kernels of all non-zero Lyapunov exponents in $\mathcal{L}^\alpha(\mu)$ it follows from Theorems 13.1 and 7.7 that both $h_\mu(\alpha(n))$ and $h_\mu(\alpha(n) \mid \mathcal{A}^\psi)$ coincide with linear functions on $C' \cap \mathbb{Z}^d$ contradicting the choice of $C$ above.

It thus follows that every coarse Lyapunov exponent $\hat\chi \in \hat{\mathcal{L}}^{\hat\alpha}(\nu)$ that contributes entropy to $\hat\alpha$ is proportional to a Lyapunov exponent of $\alpha$. We say $\chi \in \hat{\mathcal{L}}^{\hat\alpha}(\nu)$ is *essential* if
$$h_\nu(\hat\alpha(n) \mid \mathscr{W}^\chi) > 0$$
for some (and hence all) $n \in \mathbb{Z}^d$ with $\chi(n) > 0$. Let $\hat{\mathcal{L}}_{\mathrm{ess}}^{\hat\alpha}(\nu)$ denote the essential coarse Lyapunov exponents of the actions of $\hat\alpha$ on $(N, \nu)$. As remarked above, all essential exponents $\hat\chi$ of $\hat\alpha$ are proportional to Lyapunov exponents of $\alpha$. We show that they are, in fact, positively proportional.

**Theorem 13.5.** *We have $\hat{\mathcal{L}}^{\hat{\alpha}}_{\mathrm{ess}}(\nu) \subset \hat{\mathcal{L}}^{\alpha}(\mu)$.*

Analogous statements to Theorem 13.5 are established (for all coarse Lyapunov exponents) in [KRH, Section 6.2] and [KK, Lemma 2.3] under the assumption that the factor map $\psi$ is Hölder continuous using the exponential contraction along stable manifolds. Our more general statement in Theorem 13.5 follows using only entropy considerations and Theorem 7.7.

Given $\chi \in \hat{\mathcal{L}}^{\alpha}(\mu)$ let $\hat{\chi} \in \hat{\mathcal{L}}^{\hat{\alpha}}(\nu)$ denote the equivalence class of exponents positively proportional to those in $\chi$ if such a class exists; if no such class exists let $\hat{\chi}$ be the 0 linear functional. If $\hat{\chi} = 0$, let $\mathscr{W}^{\hat{\chi}}$ denote the point partition on $(N, \nu)$.

Recall the classical Abramov–Rohlin formula (19). We establish an analogous formula for entropy subordinate to coarse Lyapunov foliations under a measurable factor map $\psi$ between smooth $\mathbb{Z}^d$-actions $\alpha$ and $\hat{\alpha}$.

**Theorem 13.6** (Coarse Abramov–Rohlin formula). *Let $\chi \in \hat{\mathcal{L}}^{\alpha}(\mu)$. Then for $n \in \mathbb{Z}^d$ with $\chi(n) > 0$ we have*

$$h_{\mu}(\alpha(n) \mid \mathscr{W}^{\chi}) = h_{\nu}(\hat{\alpha}(n) \mid \mathscr{W}^{\hat{\chi}}) + h_{\mu}(\alpha(n) \mid \mathscr{W}^{\chi} \vee \mathcal{A}^{\psi}).$$

We note that in Proposition 3.1 and Corollary 3.4 of [EL], an analogous result is derived in the context of joinings of homogeneous actions in which case the factor maps are smooth.

## 14. Preliminaries

Recall we take $M$ to be a $C^{\infty}$ manifold equipped with a Borel probability $\mu$ and $\alpha\colon \mathbb{Z}^d \times M \to M$ an action satisfying the hypotheses of Section 3.1 with $\mu$ ergodic. $\mathcal{L}$ denotes the Lyapunov exponents of $\alpha$ with respect to $\mu$.

### 14.1. Lyapunov hyperplanes, Weyl chambers, and complete classes of exponents.
By extending each $\lambda_i$ to a linear function $\lambda_i\colon \mathbb{R}^d \to \mathbb{R}$, the *Lyapunov hyperplane* associated to a non-zero $\lambda_i \in \mathcal{L}$ is the kernel in $\mathbb{R}^d$ of $\lambda_i$. Note that if $\lambda_i$ and $\lambda_j$ are coarsely equivalent then $\lambda_i$ and $\lambda_j$ induce the same Lyapunov hyperplane.

Recall that relative to a fixed enumeration of $\lambda_i \in \mathcal{L}$, for each $n \in \mathbb{Z}^d$ we fix a permutation as in Section 4.2 so that $\lambda_{\sigma(n)(1)}(n) \geq \lambda_{\sigma(n)(2)}(n) \geq \cdots \geq \lambda_{\sigma(n)(\ell)}(n)$. We say $n \in \mathbb{Z}^d$ is *generic* if the above inequalities are all strict and $n$ is not contained in a Lyapunov hyperplane. A (open) *Weyl chamber* is a connected component of the complement of all Lyapunov hyperplanes in $\mathbb{R}^d$. We say a non-zero $\chi \in \hat{\mathcal{L}}$ is in the *wall* of a Weyl chamber $W$ if the boundary of $W$ contains an open subset of the kernel of $\chi$. A (open) *subchamber* of a Weyl chamber $W$ is a maximal collection of generic $n \in W$ on which the permutation $n \mapsto \sigma(n)$ is constant. As Weyl chambers are open, every Weyl chamber contains a spanning set of generic $n \in \mathbb{Z}^d$. Similarly, subchambers of Weyl chambers contain spanning sets of $n \in \mathbb{Z}^d$.

Given a subset $\mathcal{I} \subset \mathcal{L}^{\alpha}(\mu)$, let $C(\mathcal{I})$ denote the *positive cone* of $\mathcal{I}$:

$$C(\mathcal{I}) := \{n \in \mathbb{Z}^d : \lambda(n) > 0 \text{ for all } \lambda \in \mathcal{I}\}.$$

More generally, if $\mathcal{F}$ is an $\alpha$-invariant, $C^{1+\beta}$-tame, measurable foliation, the *positive cone* of $\mathcal{F}$ is

$$C(\mathcal{F}) := \{n \in \mathbb{Z}^d : \mathcal{F} \text{ is expanding for } \alpha(n)\}$$

where $\mathcal{F}$ is *expanding* for $\alpha(n)$ if $\mathcal{F} \doteq \mathcal{F} \vee \mathscr{W}^u_n$. Note that if $0 \neq \chi \in \hat{\mathcal{L}}$ then $C(\chi)$ is an open half-space called the *Lyapunov half-space* associated to $\chi$.

As non-empty positive cones contain Weyl chambers, we have the following.

**Claim 14.1.** *For any $\mathcal{I}$ either $C(\mathcal{I})$ is empty or contains a spanning set of generic $n \in \mathbb{Z}^d$. The same is true for $C(\mathcal{F})$ for any $\alpha$-invariant, $C^{1+\beta}$-tame, measurable foliation.*

Let $\mathfrak{I}$ denote the set of all subsets of $\mathcal{L}$ with the following property: $\mathcal{I} \in \mathfrak{I}$ if for any $\lambda \in \mathcal{L}$ such that $\lambda(n) > 0$ for all $n \in C(\mathcal{I})$ we have $\lambda \in \mathcal{I}$.

We call such $\mathcal{I} \in \mathfrak{I}$ a *complete class of exponents*. Note that the subsets $\mathcal{I} \in \mathfrak{I}$ are saturated by coarse equivalence classes of Lyapunov exponents. From Proposition 4.5, given $\mathcal{I} \in \mathfrak{I}$ there exists a unique, $\alpha$-invariant, $C^{1+\beta}$-tame, measurable foliation $\mathscr{W}^{\mathcal{I}}$ with $\mathscr{W}^{\mathcal{I}}(x)$ tangent to $\bigoplus_{\lambda \in \mathcal{I}} E_\lambda(x)$ at almost every $x$. Moreover, for $\mathcal{I} \in \mathfrak{I}$, we have $\mathscr{W}^{\mathcal{I}} = \bigvee \mathscr{W}^u_{n_j}$ for some finite set of $n_j \in C(\mathcal{I})$.

## 14.2. **Adjacency and sufficient subsets of $\mathbb{Z}^d$.**

**Definition 14.2.** We say two non-zero coarse Lyapunov exponents $\chi, \chi' \in \hat{\mathcal{L}}$ are *adjacent* if there exist generic $n_1$ and $n_2$ in $\mathbb{Z}^d$ such that

(1) $\sigma(n_1) \neq \sigma(n_2)$ and the permutations differ by pre-composition by disjoint adjacent transpositions: there are $1 \leq i_1 \leq i_1 + 2 \leq i_2 \leq \cdots \leq i_{\ell-1} \leq i_{\ell-1} + 2 \leq i_\ell \leq p - 1$ such that

$$\sigma(n_1) = \sigma(n_2) \circ \sigma_{i_1, i_1+1} \circ \sigma_{i_2, i_2+1} \circ \ldots \sigma_{i_\ell, i_\ell+1};$$

(2) $\{\lambda_{\sigma(n_1)(i_1)}, \ldots, \lambda_{\sigma(n_1)(i_\ell)}\} \subset \chi$ and $\{\lambda_{\sigma(n_1)(i_1+1)}, \ldots, \lambda_{\sigma(n_1)(i_\ell+1)}\} \subset \chi'$;

(3) $\{\lambda_{\sigma(n_2)(i_1)}, \ldots, \lambda_{\sigma(n_2)(i_\ell)}\} \subset \chi'$ and $\{\lambda_{\sigma(n_2)(i_1+1)}, \ldots, \lambda_{\sigma(n_2)(i_\ell+1)}\} \subset \chi$;

(4) $\mathrm{sgn}(\overline{\chi}(n_1)) = \mathrm{sgn}(\overline{\chi}(n_2))$ for all $\overline{\chi} \in \hat{\mathcal{L}}$.

Note that if (1)–(3) of Definition 14.2 hold then for $1 \leq k \leq \ell$ we have $\lambda_{\sigma(n_1)(i_k+1)} = \lambda_{\sigma(n_2)(i_k)}$ and $\lambda_{\sigma(n_1)(i_k)} = \lambda_{\sigma(n_2)(i_k+1)}$.

**Definition 14.3.** We say two generic $n_1$ and $n_2$ in $\mathbb{Z}^d$ are *adjacent* if either

(a) $\sigma(n_1) = \sigma(n_2)$ and $\mathrm{sgn}(\chi(n_1)) = \mathrm{sgn}(\chi(n_2))$ for all $\chi \in \hat{\mathcal{L}}$; or

(b) there exist non-zero, adjacent $\chi, \chi' \in \hat{\mathcal{L}}$ so that (1)–(3) of Definition 14.2 hold with $n_1$ and $n_2$;

(c) there is a non-zero $\chi$ with
    (a) either $\chi(n_1) < 0 < \chi(n_2)$ or $\chi(n_2) < 0 < \chi(n_1)$,
    (b) $\mathrm{sgn}(\chi'(n_1)) = \mathrm{sgn}(\chi'(n_2))$ for all $\chi' \neq \chi$,
    (c) and $\sigma(n_1)^{-1}(i) = \sigma(n_2)^{-1}(i)$ for every $\lambda_i \neq 0$ with $\lambda_i \notin \chi$.

Note in condition (c) of Definition 14.3 that $n_1$ and $n_2$ are necessarily in disjoint Weyl chambers $W_1$ and $W_2$ and $\mathrm{sgn}(\chi'(n_1)) = \mathrm{sgn}(\chi(n_2))$ for all $\chi' \neq \chi \in \hat{\mathcal{L}}$.

**Definition 14.4.** We say two Weyl chambers $W_1$ and $W_2$ are *adjacent* if there is a $\chi$ in the wall of $W_1$ and $W_2$ and $n_1$ and $n_2$ such that (c) of Definition 14.3 holds.

To emphasize the role of $\chi$ in Definition 14.4 we sometimes say $W_1$ and $W_2$ are *adjacent through $\chi$*.

Given $n \in \mathbb{Z}^d$, we say that $\chi(n) < \chi'(n)$ if

$$\lambda(n) < \lambda'(n)$$

for all $\lambda \in \chi$ and $\lambda' \in \chi'$.

**Definition 14.5.** A subset $S \subset \mathbb{Z}^d$ is *sufficient* in a Weyl chamber $W$ if

(1) every $n \in S$ is generic and $S$ contains a spanning set of $\mathbb{R}^d$ in every subchamber of $W$;

(2) for every $n, m \in S \cap W$ there is a sequence $n = n_0, n_1, \ldots, n_\ell = m$ with $n_i \in S \cap W$ and such that $n_{i-1}$ and $n_i$ are adjacent for $1 \leq i \leq \ell$;

(3) for every non-zero $\chi$ in the wall of $W$, there is a $n \in S \cap W$ such that, for every $\chi' \in \hat{\mathcal{L}}$ with $\chi' \neq \chi$,
  - if $\chi(n) > 0$ then either $\chi'(n) \leq 0$ or $\chi(n) < \chi'(n)$;
  - if $\chi(n) < 0$ then either $\chi'(n) \geq 0$ or $\chi(n) > \chi'(n)$.

**Definition 14.6.** A subset $S \subset \mathbb{Z}^d$ is *sufficient* in a collection $W_1, W_2, \ldots, W_\ell$ of Weyl chambers if

(1) for every $n, m \in S$ there is a sequence $n = n_0, n_1, \ldots, n_\ell = m$ with $n_i \in S$ and and such that $n_{i-1}$ and $n_i$ are adjacent for $1 \leq i \leq \ell$;

(2) $S$ is sufficient in each $W_i$.

We have the following claim.

**Lemma 14.7.** *Let $\mathcal{I} \in \mathfrak{I}$ be a complete family of Lyapunov exponents. Then there exists a finite set $S$ that is sufficient in $C(\mathcal{I})$.*

14.3. **Increasing partitions subordinate to expanding foliation.** Let $\mathcal{F}$ be an $\alpha$-invariant, $C^{1+\beta}$-tame, measurable foliation. Following Remark 8.3 of Section 8.2, for each $n \in C(\mathcal{F})$, there is a measurable partition $\xi_n^{\mathcal{F}}$ of $(M, \mu)$ subordinate to $\mathcal{F}$ and increasing for $\alpha(n)$. By adapting the construction in Remark 8.3, as in [Hu, Section 8] we obtain the following.

**Proposition 14.8.** *Let $n_1, \ldots, n_\ell \in C(\mathcal{F})$. Then there exists a measurable partition $\xi^{\mathcal{F}}$ of $(M, \mu)$ with*

(1) $\xi^{\mathcal{F}}$ *subordinate to $\mathcal{F}$;*
(2) $\alpha(n_i)\xi^{\mathcal{F}} \leq \xi^{\mathcal{F}}$ *for $i = 1, \ldots, \ell$;*
(3) $\xi^{\mathcal{F}}$ *generates for $\alpha(n_i)$; that is $\bigvee_{k=0}^{\infty} \alpha(-kn_i)\xi^{\mathcal{F}}$ is the point partition.*

Note that if $\alpha(n)\xi^{\mathcal{F}} \leq \xi^{\mathcal{F}}$ and $\alpha(m)\xi^{\mathcal{F}} \leq \xi^{\mathcal{F}}$ then

$$\alpha(n + m)\xi^{\mathcal{F}} \leq \xi^{\mathcal{F}}.$$

14.4. **Linearity of entropy on positive cones.** We have the following adaptation of [Hu, Proposition 9.1].

**Lemma 14.9.** *Let $\mathcal{F}$ be an $\alpha$-invariant, $C^{1+\beta}$-tame, measurable foliation. Let $n_1, n_2 \in C(\mathcal{F})$ and let $\eta$ be a measurable partition that is increasing for $\alpha(n_1)$ and $\alpha(n_2)$. Then*

$$h_\mu(\alpha(n_1 + n_2) \mid \mathcal{F} \vee \eta) = h_\mu(\alpha(n_2) \mid \mathcal{F} \vee \eta) + h_\mu(\alpha(n_1) \mid \mathcal{F} \vee \eta).$$

*In particular, if $\eta$ is $\alpha$-invariant then*

$$n \mapsto h_\mu(\alpha(n) \mid \mathcal{F} \vee \eta)$$

*coincides on $C(\mathcal{F})$ with a linear function.*

*Proof.* Let $n_1, n_2 \in C(\mathcal{F})$. Take a measurable partition $\xi^{\mathcal{F}}$ of $(M, \mu)$ with $\xi^{\mathcal{F}}$ subordinate to $\mathcal{F}$ and increasing for $\alpha(n_1)$ and $\alpha(n_2)$ as in Proposition 14.8. It follows from Claim 12.1 that

$$h_\mu(\alpha(m) \mid \mathcal{F} \vee \eta) = H_\mu(\alpha(-m)(\xi^{\mathcal{F}} \vee \eta) \mid \xi^{\mathcal{F}} \vee \eta)$$

for $m = n_1, n_2$ and $m = n_1 + n_2$. Moreover

$$h_\mu(\alpha(n_1 + n_2) \mid \mathcal{F} \vee \eta) = H_\mu(\alpha(-n_1 - n_2)(\xi^\mathcal{F} \vee \eta) \mid \xi^\mathcal{F} \vee \eta)$$
$$= H_\mu(\alpha(-n_1)(\alpha(-n_2)(\xi^\mathcal{F} \vee \eta)) \vee (\alpha(-n_2)(\xi^\mathcal{F} \vee \eta)) \mid \xi^\mathcal{F} \vee \eta)$$
$$= H_\mu(\alpha(-n_2)(\xi^\mathcal{F} \vee \eta) \mid \xi^\mathcal{F} \vee \eta) + H_\mu(\alpha(-n_1)(\alpha(-n_2)(\xi^\mathcal{F} \vee \eta)) \mid \alpha(-n_2)(\xi^\mathcal{F} \vee \eta))$$
$$= h_\mu(\alpha(n_2) \mid \mathcal{F} \vee \eta) + h_\mu(\alpha(n_1) \mid \mathcal{F} \vee \eta). \qquad \square$$

## 15. Key proposition and proof of Theorem 13.1

### 15.1. Conditional dimensions and dependence on coarse conditional measures.
Consider an $\alpha$-invariant, $C^{1+\beta}$-tame, measurable foliation $\mathcal{F}$ and a complete set $\mathcal{I} \in \mathfrak{I}$ of Lyapunov exponents. According to Proposition 4.5, for $n \in C(\mathcal{I})$ we obtain a filtration by $\alpha$-invariant, $C^{1+\beta}$-tame, measurable foliations

$$\{x\} \subset \mathcal{F}_n^{1,\mathcal{I}}(x) \subset \mathcal{F}_n^{2,\mathcal{I}}(x) \subset \cdots \subset \mathcal{F}_n^{u(n),\mathcal{I}}(x) := \mathcal{F}^\mathcal{I}(x) \qquad (31)$$

where $\mathcal{F}^\mathcal{I} = \mathcal{F} \vee \mathcal{W}^\mathcal{I}$ and $\mathcal{F}_n^{j,\mathcal{I}} := \mathcal{W}_n^{u,j} \vee \mathcal{F}^\mathcal{I}$.

For our main proposition, fix $\mathcal{I} \in \mathfrak{I}$, let $S$ be a sufficient set in $C(\mathcal{I})$, and let $\eta$ be a measurable partition such that $\alpha(n)\eta \leq \eta$ for all $n \in S$. Consider $n \in S \cap C(\mathcal{I})$. As $n$ is generic, we have $\lambda_i(n) \neq \lambda_j(n)$ for all $\lambda_i \neq \lambda_j \in \mathcal{L}$. For each $1 \leq j \leq \ell$ let $\xi_n^{j,\mathcal{I},\mathcal{F}}$ be a measurable partition of $(M, \mu)$ subordinate to $\mathcal{F}_n^{j,\mathcal{I}}$ as in Remark 8.3. We write

$$\dim_n^j(\mu \mid \mathcal{F} \vee \eta \mid \mathcal{I}) := \lim \frac{\log(\mu_x^{\xi_n^{j,\mathcal{I},\mathcal{F}} \vee \eta}(B(x,r)))}{\log(r)}.$$

The limit exists by Proposition 7.4 as $\xi_n^{j,\mathcal{I},\mathcal{F}} \vee \eta$ is increasing for $\alpha(n)$. With $\dim_n^0(\mu \mid \mathcal{F} \vee \eta \mid \mathcal{I}) := 0$, for $n \in S \cap C(\mathcal{I})$ and $\lambda_i$ such that $\lambda_i(n) > 0$ we write

$$\gamma_n(\lambda_i \mid \mathcal{F} \vee \eta \mid \mathcal{I}) := \dim_n^{\sigma(n)^{-1}(i)}(\mu \mid \mathcal{F} \vee \eta \mid \mathcal{I}) - \dim_n^{\sigma(n)^{-1}(i)-1}(\mu \mid \mathcal{F} \vee \eta \mid \mathcal{I}).$$

From Proposition 7.6 we have the following observation.

**Claim 15.1.** *If $\mathcal{I}' \subset \mathfrak{I}$ satisfies $\mathcal{I}' \subset \mathcal{I}$ and if $\hat{\eta}$ is a measurable partition with $\eta \leq \hat{\eta}$ and $\alpha(n)\hat{\eta} \leq \hat{\eta}$ for all $n \in S$ then for $n \in S$ and all $\lambda_i$ with $\lambda_i(n) > 0$*

$$\gamma_n(\lambda_i \mid \mathcal{F} \vee \hat{\eta} \mid \mathcal{I}') \leq \gamma_n(\lambda_i \mid \mathcal{F} \vee \eta \mid \mathcal{I}).$$

*Proof.* Take $\xi_n^{\mathcal{I},\mathcal{F}}$ and $\xi_n^{\mathcal{I}',\mathcal{F}}$ with $\xi_n^{\mathcal{I},\mathcal{F}} \leq \xi_n^{\mathcal{I}',\mathcal{F}}$ to be measurable partitions of $(M, \mu)$ subordinated to $\mathcal{W}^{\mathcal{I} \vee \mathcal{F}}$ and $\mathcal{W}^{\mathcal{I}' \vee \mathcal{F}}$, respectively, with $\alpha(n)\xi_n^{\mathcal{I},\mathcal{F}} \leq \xi_n^{\mathcal{I},\mathcal{F}}$ and $\alpha(n)\xi^{\mathcal{I}',\mathcal{F}} \leq \xi^{\mathcal{I}',\mathcal{F}}$ as in Remark 8.3. Then

$$\hat{\eta} \vee \xi^{\mathcal{I}',\mathcal{F}} \geq \eta \vee \xi^{\mathcal{I},\mathcal{F}}$$

whence by Proposition 7.6 we have

$$\gamma_n(\lambda_i \mid \mathcal{F} \vee \hat{\eta} \mid \mathcal{I}') \leq \gamma_n(\lambda_i \mid \mathcal{F} \vee \eta \mid \mathcal{I}). \quad \square$$

We show that often the inequality in Claim 15.1 is an equality as the element $n$ and set $\mathcal{I}$ varies.

**Proposition 15.2.** *Fix $\mathcal{I} \in \mathfrak{I}$. Let $S$ be a sufficient set in $C(\mathcal{I})$ and suppose $\eta$ is measurable partition such that $\alpha(n)\eta \leq \eta$ for all $n \in S$. Then*

*(a) for every Weyl chamber $W \subset C(\mathcal{I})$ and $n_1, n_2 \in W \cap S$ we have*

$$\gamma_{n_1}(\lambda_i \mid \mathcal{F} \vee \eta \mid \mathcal{I}) = \gamma_{n_2}(\lambda_i \mid \mathcal{F} \vee \eta \mid \mathcal{I})$$

*for all $\lambda_i$ with $\lambda_i(n_1) > 0$;*

(b) *for every pair of adjacent Weyl chambers $W_1, W_2 \subset C(\mathcal{I})$ and $n_i \in W_i \cap S$ we have*

$$\gamma_{n_1}(\lambda_i \mid \mathcal{F} \vee \eta \mid \mathcal{I}) = \gamma_{n_2}(\lambda_i \mid \mathcal{F} \vee \eta \mid \mathcal{I})$$

*for all $\lambda_i$ with $\lambda_i(n_1) > 0$ and $\lambda_i(n_2) > 0$;*

(c) *for every Weyl chamber $W \subset C(\mathcal{I})$, $\chi'$ in the wall of $W$, and $n \in W \cap S$ we have*

$$\gamma_n(\lambda_i \mid \mathcal{F} \vee \eta \mid \mathcal{I}) = \gamma_n(\lambda_i \mid \mathcal{F} \vee \eta \mid \mathcal{I} \smallsetminus \chi')$$

*for all $\lambda_i \notin \chi'$ with $\lambda_i(n) > 0$.*

From Proposition 15.2, for $\lambda_i \in \mathcal{I}$ we may write

$$\gamma(\lambda_i \mid \mathcal{F} \vee \eta \mid \mathcal{I}) := \gamma_n(\lambda_i \mid \mathcal{F} \vee \eta \mid \mathcal{I})$$

where $n$ is any choice of $n \in S \cap C(\mathcal{I})$.

Note that in the case that $\eta$ is $\alpha$-invariant we may take the set $S$ in Proposition 15.2 to be the collection of all generic elements in which case $\gamma_n(\lambda_i \mid \mathcal{F} \vee \eta \mid \mathcal{I})$ is defined for every generic $n \in C(\mathcal{I})$. In this case we obtain that for $\lambda_i \in \mathcal{I}$, the number $\gamma(\lambda_i \mid \mathcal{F} \vee \eta \mid \mathcal{I})$ depends only on the conditional measures along the coarse manifolds $\mathscr{W}^\chi$ where $\chi$ is the coarse Lyapunov exponent containing $\lambda_i$.

**Corollary 15.3.** *Let $\eta$ be $\alpha$-invariant. Given $\mathcal{I} \in \mathfrak{I}$ with $C(\mathcal{I}) \neq \varnothing$, any coarse exponent $\chi \subset \mathcal{I}$, and $\lambda_i \in \chi$ we have*

$$\gamma(\lambda_i \mid \mathcal{F} \vee \eta \mid \mathcal{I}) = \gamma(\lambda_i \mid \mathcal{F} \vee \eta \mid \chi).$$

### 15.2. Proof of Theorem 13.1. Theorem 13.1 follows immediately from the above.

*Proof of Theorem 13.1.* First consider a generic $n \in \mathbb{Z}^d$. Let $\mathcal{U}(n) := \{\lambda_i \in \mathcal{L} : \lambda_i(n) > 0\}$ be the collection of positive Lyapunov exponents of $n$. We have $\mathcal{U}(n) \in \mathfrak{I}$. From Theorem 7.7, with $f = \alpha(n)$ we have that

$$h_\mu(\alpha(n) \mid \mathcal{F} \vee \eta) = \sum_{\lambda_i(n) > 0} \gamma_n(\lambda_i \mid \mathcal{F} \vee \eta \mid \mathcal{U}(n)) \lambda_i(n).$$

On the other hand, for any coarse Lyapunov exponent $\chi \in \hat{\mathcal{L}}$ with $\chi(n) > 0$, we have (again from Theorem 7.7) that

$$h_\mu(\alpha(n) \mid \mathcal{F} \vee \mathscr{W}^\chi \vee \eta) = \sum_{\lambda_i \in \chi} \gamma_n(\lambda_i \mid \mathcal{F} \vee \eta \mid \chi) \lambda_i(n).$$

The result then follows for all generic $n$ by Corollary 15.3.

Now, consider a non-generic $n \in \mathbb{Z}^d$. From Claim 12.1 we have

$$h_\mu(\alpha(n) \mid \mathcal{F} \vee \eta) = h_\mu(\alpha(n) \mid \mathcal{F} \vee \mathscr{W}_n^u \vee \eta).$$

Take $\hat{\mathcal{F}} = \mathcal{F} \vee \mathscr{W}_n^u$. Then $\hat{\mathcal{F}}$ is expanding for $\alpha(n)$.

Since $C(\hat{\mathcal{F}}) \neq \varnothing$, it follows that $C(\hat{\mathcal{F}})$ contains a spanning set of generic points. From Lemma 14.9, the functions

$$n \mapsto h_\mu(\alpha(n) \mid \hat{\mathcal{F}} \vee \eta)$$

and

$$n \mapsto \sum_{\{\chi \in \hat{\mathcal{L}} : \chi(n) > 0\}} h_\mu(\alpha(n) \mid \hat{\mathcal{F}} \vee \mathscr{W}^\chi \vee \eta)$$

extend from $C(\hat{\mathcal{F}})$ to linear functions on $\mathbb{R}^d$ and coincide on a spanning set. It follows that they agree on $C(\hat{\mathcal{F}})$. As $\hat{\mathcal{F}} \vee \mathscr{W}^\chi = \mathcal{F} \vee \mathscr{W}^\chi$ for all $\chi \in \hat{\mathcal{L}}$ with $\chi(n) > 0$, the result follows for $n$. $\qquad\square$

15.3. **Proof of Proposition 15.2 and Corollary 15.3.** We establish the proposition and corollary, completing the proof of Theorem 13.1.

*Proof of Proposition 15.2.* Fix $\mathcal{I} \in \mathfrak{I}$. Fix a Weyl chamber $W \subset C(\mathcal{I})$ and consider $n_1, n_2 \in W \cap S$. We prove (a) of Proposition 15.2 by induction on the index of the intermediate foliations in the filtration (31) corresponding to $n_1$ and $n_2$. Note for $\lambda_i \in \mathcal{L}$, $\lambda_i(n_1) > 0$ if and only if $\lambda_i(n_2) > 0$.

It follows immediately from the geometric definition of the transverse dimensions that if $\lambda_i(n_1) > 0$ then

$$\gamma_{n_1}(\lambda_i \mid \mathcal{F} \vee \eta \mid \mathcal{I}) = \gamma_{n_2}(\lambda_i \mid \mathcal{F} \vee \eta \mid \mathcal{I})$$

whenever $\sigma(n_1) = \sigma(n_2)$. Indeed in this case the intermediate foliations in (31) for $n_1$ and $n_2$ coincide.

If $\sigma(n_1) \neq \sigma(n_2)$ then, as any $n_1, n_2 \in W \cap S$ may be joined by a string of adjacent $n_i \in W \cap S$, it is therefore enough to verify the proposition in the case that $n_1$ and $n_2$ are adjacent. As no exponent changes sign in $W$, we may thus assume $n_1$ and $n_2$ satisfy (b) of Definition 14.3; that is, we assume there are adjacent $\chi, \chi'$ so that (1)–(3) of Definition 14.2 hold for $n_1$ and $n_2$. Let $i_1, \ldots, i_\ell$ be as in (1)–(3) of Definition 14.2.

Consider any $1 \leq j$. Assume by induction that we have shown

$$\gamma_{n_1}(\lambda_{\sigma(n_1)(m)} \mid \mathcal{F} \vee \eta \mid \mathcal{I}) = \gamma_{n_2}(\lambda_{\sigma(n_2)(m)} \mid \mathcal{F} \vee \eta \mid \mathcal{I})$$

for all $m < j$. Observe

**Claim 15.4.** *If neither $j \neq i_k$ nor $j \neq i_{k+1}$ then*

$$\lambda_{n_1}(\lambda_{\sigma(n_1)(j)} \mid \mathcal{F} \vee \eta \mid \mathcal{I}) = \lambda_{n_2}(\lambda_{\sigma(n_2)(j)} \mid \mathcal{F} \vee \eta \mid \mathcal{I}).$$

Indeed, if neither $j \neq i_k$ nor $j \neq i_{k+1}$ then $\mathcal{F}_{n_1}^{j,\mathcal{I}} = \mathcal{F}_{n_2}^{j,\mathcal{I}}$ and $\mathcal{F}_{n_1}^{j-1,\mathcal{I}} = \mathcal{F}_{n_2}^{j-1,\mathcal{I}}$ whence equality of transverse dimensions follows from definition.

Consider now the case $j = i_k$. Fix $i \neq i'$ with $\sigma(n_1)(i_k) = i$ and $\sigma(n_2)(i_k) = i'$. Then $\sigma(n_1)(i_k + 1) = i'$ and $\sigma(n_2)(i_k + 1) = i$. We show simultaneously the equalities

$$\begin{aligned}
\gamma_{n_1}(\lambda_i \mid \mathcal{F} \vee \eta \mid \mathcal{I}) &= \gamma_{n_2}(\lambda_i \mid \mathcal{F} \vee \eta \mid \mathcal{I}), \\
\gamma_{n_1}(\lambda_{i'} \mid \mathcal{F} \vee \eta \mid \mathcal{I}) &= \gamma_{n_2}(\lambda_{i'} \mid \mathcal{F} \vee \eta \mid \mathcal{I}).
\end{aligned} \tag{32}$$

Part (a) of the proposition then follows from induction on $j$ and Claim 15.4.

To establish (32), take

$$\hat{\mathcal{F}} = \mathcal{F} \vee \mathscr{W}_{n_1}^{u, i_{k+1}} = \mathcal{F} \vee \mathscr{W}_{n_2}^{u, i_{k+1}}$$

and

$$\widetilde{\mathcal{F}} = \mathcal{F} \vee \mathscr{W}_{n_1}^{u, i_{k-1}} = \mathcal{F} \vee \mathscr{W}_{n_2}^{u, i_{k-1}}.$$

We have that

$$h_\mu(\cdot \mid \hat{\mathcal{F}} \vee \eta) - h_\mu(\cdot \mid \widetilde{\mathcal{F}} \vee \eta)$$

extends from $S \cap W$ to a linear function $L$ on $\mathbb{R}^d$. Moreover, as $S$ contains a spanning set in each subchamber of $W$, on $S \cap W$ the function $L$ coincides simultaneously with the linear functions given by the expressions

- $n \mapsto \gamma_{n_1}(\lambda_i \mid \mathcal{F} \vee \eta \mid \mathcal{I})\lambda_i(n) + \gamma_{n_1}(\lambda_{i'} \mid \mathcal{F} \vee \eta \mid \mathcal{I})\lambda_{i'}(n)$
- $n \mapsto \gamma_{n_2}(\lambda_{i'} \mid \mathcal{F} \vee \eta \mid \mathcal{I})\lambda_{i'}(n) + \gamma_{n_2}(\lambda_i \mid \mathcal{F} \vee \eta \mid \mathcal{I})\lambda_i(n)$.

As $\lambda_i$ and $\lambda_{i'}$ are linearly independent and as $S \cap W$ spans $\mathbb{R}^d$, (32) and conclusion (a) follow.

For (b) consider adjacent $W_1$ and $W_2$ in $C(\mathcal{I})$. Let $\chi$ be as in Definition 14.4 and take $n_i \in W_i \cap S$ satisfying (c) of Definition 14.3 with $\chi(n_1) > 0 > \chi(n_2)$. It follows that

$\chi \notin \mathcal{I}$. Let

$$\ell := \max\{j : \lambda_{\sigma(n_1)(j)} \in \mathcal{U}(n_1) \smallsetminus \chi\} = \max\{j : \lambda_{\sigma(n_2)(j)}(n_2) > 0 \text{ and } \lambda_{\sigma(n_1)(j)}(n_1) > 0\}.$$

Since $\sigma(n_1)(j) = \sigma(n_2)(j)$ for all $1 \le j \le \ell$ it follows that the filtrations (31) corresponding to $n_1$ and $n_2$ coincide for all $1 \le j \le \ell$ whence the transverse dimensions coincide for all $\lambda_i \notin \chi$ with $\lambda_i(n_1) > 0$; that is, for all $\lambda_i$ with $\lambda_{\sigma(n_2)(j)}(n_2) > 0$ and $\lambda_{\sigma(n_1)(j)}(n_1) > 0$.

For (c), take $\chi'$ in the wall of $W$ and $n \in S \cap W$ with

$$0 < \lambda'(n) < \lambda(n) \text{ for all } \lambda' \in \chi' \text{ and } \lambda \notin \chi' \text{ with } \lambda(n) > 0. \tag{33}$$

Let

$$\ell = \max\{j : \lambda_{\sigma(n)(j)} \in \mathcal{U}(n) \smallsetminus \chi'\}.$$

Then $\{\lambda_{\sigma(n)(j)} : 1 \le j \le \ell\} \cap \mathcal{I} = \mathcal{I} \smallsetminus \chi'$. From (33) we have

$$\mathcal{F}_n^{j,\mathcal{I}}(x) = \mathcal{F}_n^{j,\mathcal{I} \smallsetminus \chi'}(x)$$

for all $1 \le j \le \ell$. The equality of transverse dimensions follows from definition. $\qquad\square$

We finish with the proof of Corollary 15.3. Recall we now assume $\eta$ is $\alpha$-invariant whence for $\mathcal{I} \in \mathfrak{I}$ and $\lambda_i \in \mathcal{I}$ the dimension $\gamma_n(\lambda_j \mid \mathcal{F} \vee \eta \mid \mathcal{I})$ is defined for every generic $n \in C(\mathcal{I})$ and independent of $n \in C(\mathcal{I})$.

*Proof of Corollary 15.3.* Given $\mathcal{I} \in \mathfrak{I}$ and $\lambda_i \in \mathcal{I}$ let $\chi \in \hat{\mathcal{L}}$ be the coarse Lyapunov exponent containing $\lambda_i$. Let $\chi' \subset \mathcal{I}$ be such that $\chi' \ne \chi$ and $\chi'$ is in the wall of some Weyl chamber $W \subset C(\mathcal{I})$. Let $\mathcal{I}_1 = \mathcal{I} \smallsetminus \{\chi'\}$. Taking a generic $n \in W$, from Proposition 15.2(c) it follows for every $\lambda_j \in \mathcal{I}_1$ that

$$\gamma_n(\lambda_j \mid \mathcal{F} \vee \eta \mid \mathcal{I}) = \gamma_n(\lambda_j \mid \mathcal{F} \vee \eta \mid \mathcal{I}_1).$$

Proceeding recursively, we define $\mathcal{I} \supset \mathcal{I}_1 \supset \cdots \supset \mathcal{I}_r := \chi$ with

$$\gamma_n(\lambda_j \mid \mathcal{F} \vee \eta \mid \mathcal{I}) = \gamma_n(\lambda_j \mid \mathcal{F} \vee \eta \mid \mathcal{I}_j)$$

for each $1 \le j \le r$. $\qquad\square$

## 16. Proof of Theorems 13.5 and 13.6

We retain all notation from Section 13.2. In particular $\hat{\alpha}$ is a measurable factor of $\alpha$ induced by $\psi$ and $\mathcal{A}^\psi$ is the $\alpha$-invariant measurable partition on $(M, \mu)$ induced by the factor map $\psi$.

### 16.1. Key Lemma.
Consider any non-zero exponent $\hat{\chi} \in \hat{\mathcal{L}}^{\hat{\alpha}}(\nu)$. Note that we not yet show the equivalence class $\hat{\chi}$ is an element of $\hat{\mathcal{L}}^\alpha(\mu)$. Take $H = C(\hat{\chi})$ to be the Lyapunov half-space associated with $\hat{\chi}$.

Considering all non-zero $\hat{\chi}' \in \mathcal{L}^{\hat{\alpha}}(\nu)$ and $\chi' \in \mathcal{L}^\alpha(\mu)$ as Lyapunov functionals on $\mathbb{R}^d$, consider *joint Weyl chambers* as connected subsets of the complement of all Lyapunov hyperplanes of all non-zero $\hat{\chi}' \in \mathcal{L}^{\hat{\alpha}}(\nu)$ and $\chi' \in \mathcal{L}^\alpha(\mu)$. Then $H$ is saturated by joint Weyl chambers and we may take a finite set $S \subset H \cap \mathbb{Z}^d$ that is sufficient in $H$ (where sufficiency is relative to the collection of joint Weyl chambers).

Take $\hat{\eta}$ a measurable partition of $(N, \nu)$ as in Proposition 14.8 that is subordinate to $\mathscr{W}^{\hat{\chi}}$ with $\hat{\alpha}(n)(\hat{\eta}) \le \hat{\eta}$ for all $n \in S$. Let $\eta = \psi^{-1}(\hat{\eta})$.

The key observation in the proof of Theorems 13.5 and 13.6 is that for $n \in S \subset H$, every coarse exponent $\chi \in \hat{\mathcal{L}}^\alpha(\mu)$ with $\chi \ne \hat{\chi}$ only contributes fiber-entropy to $h_\mu(\alpha(n) \mid$

$\eta$). Let

$$\mathcal{U}_\mu^\alpha(n) := \{\chi' \in \hat{\mathcal{L}}^\alpha(\mu) : \chi'(n) > 0\}.$$

**Lemma 16.1.** *For $n \in S \subset H$ and $\chi \in \hat{\mathcal{L}}^\alpha(\mu)$ with $\chi \neq \hat{\chi}$ and $\chi(n) > 0$ we have*

$$\gamma_n(\lambda_j \mid \eta \mid \mathcal{U}_\mu^\alpha(n)) = \gamma_n(\lambda_j \mid \mathcal{A}^\psi \mid \chi)$$

*for all $\lambda_j \in \chi$.*

*Proof.* First note that for any $n \in S$ with $\chi(n) > 0$, from Claim 15.1 we have

$$\gamma_n(\lambda_j \mid \eta \mid \mathcal{U}_\mu^\alpha(n)) \geq \gamma_n(\lambda_j \mid \eta \mid \chi) \geq \gamma_n(\lambda_j \mid \mathcal{A}^\psi \mid \chi) \tag{34}$$

for all $\lambda_j \in \chi$.

To prove the reverse inequality, note that as $\chi \neq \hat{\chi}$ the Lyapunov hyperplane associated to $\chi$ intersects the interior of $H$. In particular, there are joint Weyl chambers $W_1 \subset H$ and $W_2 \subset H$ that are adjacent through $\chi$ and $n_i \in W_i \cap S \subset H$ so that (c) of Definition 14.3 holds with $\chi(n_1) > 0$ and $\chi(n_2) < 0$.

Let $W \subset C(\chi) \cap H$ be the joint Weyl chamber containing the $n$ in the lemma. As $S$ is sufficient in $C(\chi) \cap H$ there is a sequence of subsequently adjacent joint Weyl chambers $W = W^1, W^1, \ldots, W^{\ell+1} = W_1$ in $C(\chi) \cap H$ and a sequence of coarse exponents $\chi_1', \chi_2', \ldots, \chi_\ell'$ in $\hat{\mathcal{L}}^\alpha(\mu)$ with $\chi_j' \neq \hat{\chi}$ for $1 \leq j \leq \ell$ such that each pair $W^j, W^{j+1}$ is adjacent through $\chi_j'$. Note also that since each $W^j \subset H$, we have $\chi_j' \neq \hat{\chi}$ for $1 \leq j \leq \ell$. Let $\mathcal{I} = \mathcal{U}_\mu^\alpha(n) \smallsetminus \{\chi_1', \chi_2', \ldots, \chi_\ell'\}$. Then also $\mathcal{I} = \mathcal{U}_\mu^\alpha(n_1) \smallsetminus \{\chi_1', \chi_2', \ldots, \chi_\ell'\}$. From Proposition 15.2 we have for $\lambda_j \in \chi$

(1) $\gamma_n(\lambda_j \mid \eta \mid \mathcal{U}_\mu^\alpha(n)) = \gamma_{n_1}(\lambda_j \mid \eta \mid \mathcal{I})$;
(2) $\gamma_{n_1}(\lambda_j \mid \eta \mid \mathcal{U}_\mu^\alpha(n_1)) = \gamma_n(\lambda_j \mid \eta \mid \mathcal{I})$;
(3) $\gamma_n(\lambda_j \mid \eta \mid \mathcal{I}) = \gamma_{n_1}(\lambda_j \mid \eta \mid \mathcal{I})$;
(4) $\gamma_n(\lambda_j \mid \mathcal{A}^\psi \mid \chi) = \gamma_{n_1}(\lambda_j \mid \mathcal{A}^\psi \mid \chi)$.

To prove the lemma, it is thus sufficient to show for $\lambda_i \in \chi$ that

$$\gamma_{n_1}(\lambda_j \mid \eta \mid \mathcal{U}_\mu^\alpha(n_1)) = \gamma_{n_1}(\lambda_j \mid \mathcal{A}^\psi \mid \chi). \tag{35}$$

From Theorem 7.7, for any $n \in S$ we have

$$h_\mu(\alpha(n) \mid \eta) = \sum_{\lambda_i \in \mathcal{U}_\mu^\alpha(n)} \lambda_i(n) \gamma_n(\lambda_i \mid \eta \mid \mathcal{U}_\mu^\alpha(n)).$$

Moreover, from the Abramov–Rohlin formula (18), Theorem 7.7, and Theorem 13.1 (applied to the trivial foliation and $\alpha$-invariant partition $\mathcal{A}^\psi$) we have for any $n \in \mathbb{Z}^d$ that

$$h_\mu(\alpha(n) \mid \eta) = h_\nu(\hat{\alpha}(n) \mid \hat{\eta}) + h_\mu(\alpha(n) \mid \mathcal{A}^\psi)$$

$$= h_\nu(\hat{\alpha}(n) \mid \hat{\eta}) + \sum_{\chi \in \mathcal{U}_\mu^\alpha(n)} h_\mu(\alpha(n) \mid \mathscr{W}^\chi \vee \mathcal{A}^\psi).$$

$$= h_\nu(\hat{\alpha}(n) \mid \hat{\eta}) + \sum_{\chi' \in \mathcal{U}_\mu^\alpha(n)} \sum_{\lambda_i \in \chi'} \lambda_i(n) \gamma_n(\lambda_i \mid \mathcal{A}^\psi \mid \chi').$$

whence

$$h_\nu(\hat{\alpha}(n) \mid \hat{\eta}) = \sum_{\lambda_i(n) \in \mathcal{U}_\mu^\alpha(n)} \lambda_i(n) \gamma_n(\lambda_i \mid \eta \mid \mathcal{U}(n)) - \sum_{\chi' \in \mathcal{U}_\mu^\alpha(n)} \sum_{\lambda_i \in \chi'} \lambda_i(n) \gamma_n(\lambda_i \mid \mathcal{A}^\psi \mid \chi').$$

It follows from Proposition 14.8 that on $S \subset H$, that the map

$$n \mapsto h_\nu(\hat{\alpha}(n) \mid \hat{\eta})$$

coincides with a linear function $L\colon \mathbb{R}^d \to \mathbb{R}$. From Propositions 14.8 and 15.2 it follows that for all $n \in W_1 \cap S$ and $m \in W_2 \cap S$

$$L(n) = \sum_{\lambda_i \in \mathcal{U}_\mu^\alpha(n_1)} \lambda_i(n)\gamma_{n_1}(\lambda_i \mid \eta \mid \mathcal{U}_\mu^\alpha(n_1)) - \sum_{\chi' \in \mathcal{U}_\mu^\alpha(n_1)} h_\mu(\alpha(n) \mid \mathscr{W}^{\chi'} \vee \mathcal{A}^\psi) \quad (36)$$

$$L(m) = \sum_{\lambda_i \in \mathcal{I}^*} \lambda_i(m)\gamma_{n_2}(\lambda_i \mid \eta \mid \mathcal{U}_\mu^\alpha(n_2)) - \sum_{\chi' \in \mathcal{I}^*} h_\mu(\alpha(m) \mid \mathscr{W}^{\chi'} \vee \mathcal{A}^\psi). \quad (37)$$

where $\mathcal{I}^* = (\mathcal{U}_\mu^\alpha(n_1) \smallsetminus \{-\chi\}) \cup \{-\chi\}$ or $\mathcal{I}^* = \mathcal{U}_\mu^\alpha(n_1) \smallsetminus \chi$ depending, respectively, on whether or not $-\chi$ is a coarse exponent in $\hat{\mathcal{L}}^\alpha(\mu)$.

From Proposition 15.2, for $\chi' \neq \chi$ and $\chi' \neq -\chi$ and $\lambda_i \in \chi'$ we have equalities

$$\gamma_{n_1}(\lambda_i \mid \eta \mid \mathcal{U}_\mu^\alpha(n_1)) = \gamma_{n_2}(\lambda_i \mid \eta \mid \mathcal{U}_\mu^\alpha(n_2))$$

$$\gamma_{n_1}(\lambda_i \mid \mathcal{A}^\psi \mid \chi') = \gamma_{n_2}(\lambda_i \mid \mathcal{A}^\psi \mid \chi').$$

Let $L_1$ be the function

$$L_1(n) = \sum_{\lambda_i \in \chi} \lambda_i(n)\gamma_{n_1}(\lambda_i \mid \eta \mid \mathcal{U}_\mu^\alpha(n_1)) - \sum_{\lambda_i(n) \in \chi} \lambda_i(n)\gamma_{n_1}(\lambda_i \mid \mathcal{A}^\psi \mid \chi)$$

and (if $-\chi$ is a coarse Lyapunov exponent) let

$$L_2(n) = \sum_{\lambda_i \in -\chi} \lambda_i(n)\gamma_{n_2}(\lambda_i \mid \eta \mid \mathcal{U}_\mu^\alpha(n_2)) - \sum_{\lambda_i(n) \in -\chi} \lambda_i(n)\gamma_{n_2}(\lambda_i \mid \mathcal{A}^\psi \mid -\chi).$$

Comparing righthand sides of (36) and (37) and canceling common linear terms it follows that either

$$L_1 = 0 \qquad \text{or} \qquad L_1 = L_2.$$

From (34) we have $L_1(n) \geq 0$ and $L_2(n) \leq 0$ for $n \in W_1$ whence either case above implies

$$\sum_{\lambda_i \in \overline{\chi}} \lambda_i(n_1) \left( \gamma_{n_1}(\lambda_i \mid \eta \mid \mathcal{U}_\mu^\alpha(n_1)) - \gamma_{n_1}(\lambda_i \mid \mathcal{A}^\psi \mid \chi) \right) = 0.$$

(35) then follows from (34). □

From Lemma 16.1 we obtain Theorem 13.5.

*Proof of Theorem 13.5.* We retain all notations from above. In particular, we take $\hat{\chi} \in \hat{\mathcal{L}}^{\hat{\alpha}}(\nu)$ with

$$h_\nu(\hat{\alpha}(n) \mid \hat{\chi}) > 0$$

for some $n$ with $\hat{\chi}(n) > 0$ and fix $\hat{\eta}$ as above. Suppose that every coarse exponent $\chi \in \hat{\mathcal{L}}^\alpha(\mu)$ is distinct from $\hat{\chi}$. Then, by the Abramov-Rohlin formula (19), Theorem 7.7, and Lemma 16.1 we obtain a contradiction as for any $n \in S \subset H$ we have

$$h_\nu(\hat{\alpha}(n) \mid \mathscr{W}^{\hat{\chi}}) = h_\nu(\hat{\alpha}(n) \mid \hat{\eta})$$
$$= h_\mu(\alpha(n) \mid \eta) - h_\mu(\alpha(n) \mid \mathcal{A}^\psi)$$
$$= \sum_{\lambda_i \in \mathcal{U}_\mu^\alpha(n)} \lambda_i(n)\gamma_n(\lambda_i \mid \eta \mid \mathcal{U}_\mu^\alpha(n)) - \sum_{\chi \in \hat{\mathcal{L}}^\alpha(\mu): \chi(n) > 0} \sum_{\lambda_i \in \chi} \lambda_i(n)\gamma_n(\lambda_i \mid \mathcal{A}^\psi \mid \chi)$$
$$= 0. \qquad \qquad \square$$

Having established Theorem 13.5, given any coarse exponent $\hat{\chi} \in \hat{\mathcal{L}}^{\hat{\alpha}}(\nu)$ contributing entropy to the factor system $\hat{\alpha}$ it follows that $\hat{\chi} \in \hat{\mathcal{L}}^\alpha(\mu)$. That is, if

$$\nu(\hat{\alpha}(n) \mid \mathscr{W}^\chi) > 0$$

for some $n$, then $\hat{\chi}$ is also a coarse exponent for the action of $\alpha$ on $(M, \mu)$. With $\hat{\eta}$, $\eta$, and $S \subset H = C(\hat{\chi})$ as in Lemma 16.1, we obtain the following.

**Corollary 16.2.** *For $n \in S$ we have*
$$h_\mu(\alpha(n) \mid \mathscr{W}^{\hat{\chi}} \vee \eta) = h_\nu(\hat{\alpha}(n), \hat{\eta}) + h_\mu(\alpha(n) \mid \mathcal{A}^\psi \vee \mathscr{W}^{\hat{\chi}})$$

*Proof.* Note from Corollary 15.3 that for $n \in S$ and $\chi' \subset \mathcal{U}_\mu^\alpha(n)$, for $\lambda_i \in \chi'$
$$\gamma_n(\lambda_i \mid \mathcal{A}^\psi \mid \chi') = \gamma_n(\lambda_i \mid \mathcal{A}^\psi \mid \mathcal{U}_\mu^\alpha(n))$$
and is defined independent of $n$. Let $\gamma(\lambda_i \mid \mathcal{A}^\psi)$ denote this constant.

From Lemma 16.1 and Theorem 7.7 we have
$$h_\mu(\alpha(n) \mid \eta) = \sum_{\lambda_i(n) \in \mathcal{U}_\mu^\alpha(n)} \lambda_i(n) \gamma_n(\lambda_i \mid \eta \mid \mathcal{U}_\mu^\alpha(n))$$
$$= \sum_{\lambda_i(n) \in \mathcal{U}_\mu^\alpha(n)} \lambda_i(n) \gamma(\lambda_i \mid \mathcal{A}^\psi) + \sum_{\lambda_i(n) \in \hat{\chi}} \lambda_i(n) \left( \gamma_n(\lambda_i \mid \eta \mid \mathcal{U}_\mu^\alpha(n)) - \gamma(\lambda_i \mid \mathcal{A}^\psi) \right)$$
$$= h_\mu(\alpha(n) \mid \mathcal{A}^\psi) + \sum_{\lambda_i(n) \in \hat{\chi}} \lambda_i(n) \left( \gamma_n(\lambda_i \mid \eta \mid \mathcal{U}_\mu^\alpha(n)) - \gamma(\lambda_i \mid \mathcal{A}^\psi) \right).$$

whence
$$h_\nu(\hat{\alpha}(n), \hat{\eta}) = h_\mu(\alpha(n) \mid \eta) - h_\mu(\alpha(n) \mid \mathcal{A}^\psi)$$
$$= \sum_{\lambda_i(n) \in \hat{\chi}} \lambda_i(n) \left( \gamma_n(\lambda_i \mid \eta \mid \mathcal{U}_\mu^\alpha(n)) - \gamma(\lambda_i \mid \mathcal{A}^\psi) \right).$$

From Proposition 15.2 (taking a sequence of points in $S$ meeting every Weyl chamber of $\hat{\mathcal{L}}^\alpha(\mu)$ in $H = C(\hat{\chi})$) we have for $\lambda_i \in \hat{\chi}$ that
$$\gamma_n(\lambda_i \mid \eta \mid \mathcal{U}_\mu^\alpha(n)) = \gamma_n(\lambda_i \mid \eta \mid \hat{\chi}).$$

Thus
$$h_\nu(\hat{\alpha}(n), \hat{\eta}) = \sum_{\lambda_i(n) \in \hat{\chi}} \lambda_i(n) \left( \gamma_n(\lambda_i \mid \eta \mid \hat{\chi}) - \gamma(\lambda_i \mid \mathcal{A}^\psi) \right)$$
$$= h_\mu(\alpha(n) \mid \mathscr{W}^{\hat{\chi}} \vee \eta) - h_\mu(\alpha(n) \mid \mathcal{A}^\psi \vee \mathscr{W}^{\hat{\chi}}). \qquad \square$$

**16.2. Proof of Theorem 13.6.** We proceed with the proof of the theorem.

*Proof of Theorem 13.6.* Consider any fixed, non-zero $\chi \in \hat{\mathcal{L}}^\alpha(\mu)$. Fix $n \in \mathbb{Z}^d$ with $\chi(n) > 0$ and let $W$ be the Weyl chamber of $\hat{\mathcal{L}}^\alpha(\mu)$ containing $n$. Fix a finite set $S$ that is sufficient in $W$.

Choose an enumeration of all coarse exponents $\chi_j \subset \mathcal{U}_\mu^\alpha(n)$ with $\chi = \chi_0$. For each $\chi_j \in \hat{\mathcal{L}}^\alpha(\mu)$ let $\hat{\chi}_j = \chi_j$ if $\chi_j$ is a coarse Lyapunov exponent in $\hat{\mathcal{L}}^{\hat{\alpha}}(\nu)$; otherwise, let $\hat{\chi}_j$ denote the 0 functional. For each $\chi_j \subset \mathcal{U}_\mu^\alpha(n)$ with $\hat{\chi}_j = \chi_j$ let $\hat{\eta}_j$ be a measurable partition of $(N, \nu)$ as in Proposition 14.8 that is subordinate to $\mathscr{W}^{\hat{\chi}_j}$ with $\hat{\alpha}(m)(\hat{\eta}_j) \leq \hat{\eta}_j$ for all $m \in W \cap S$. If $\hat{\chi}_j$ is not a coarse exponent in $\hat{\mathcal{L}}^{\hat{\alpha}}(\nu)$ take $\hat{\eta}_j$ to be the point partition. Let $\eta_j = \psi^{-1}(\hat{\eta}_j)$.

From Theorems 13.1 and 13.5 we have for $m \in W \cap S$
$$(1) \quad h_\nu(\hat{\alpha}(m)) = \sum_{\chi_j \in \mathcal{U}_\mu^\alpha(m)} h_\nu(\hat{\alpha}(m), \hat{\eta}_j);$$
$$(2) \quad h_\mu(\alpha(m) \mid \mathcal{A}^\psi) = \sum_{\chi_j \in \mathcal{U}_\mu^\alpha(m)} h_\mu(\alpha(m) \mid \mathcal{A}^\psi \vee \mathscr{W}^{\chi_j});$$

(3) $h_\mu(\alpha(m)) = \sum_{\chi_j \in \mathcal{U}_\mu^\alpha(m)} h_\mu(\alpha(m) \mid \mathscr{W}^{\chi_j})$.

Combined with Corollary 16.2, it follows that for $m \in W \cap S$,

$$\sum_{\chi_j \in \mathcal{U}_\mu^\alpha(m)} h_\mu(\alpha(m) \mid \mathscr{W}^{\chi_j} \vee \eta_j) = h_\nu(\hat{\alpha}(m)) + h_\mu(\alpha(m) \mid \mathcal{A}^\psi)$$

$$= \sum_{\chi_j \in \mathcal{U}_\mu^\alpha(m)} h_\mu(\alpha(m) \mid \mathscr{W}^{\chi_j}). \tag{38}$$

From Claim 8.6, for each $j$ we have

$$h_\mu(\alpha(m) \mid \mathscr{W}^{\chi_j} \vee \eta_j) \le h_\mu(\alpha(m) \mid \mathscr{W}^{\chi_j})$$

for every $m \in W \cap S$. Combined with (38) it follows for each $j$ and $m \in W \cap S$ that

$$h_\mu(\alpha(m) \mid \mathscr{W}^{\chi_j} \vee \eta_j) = h_\mu(\alpha(m) \mid \mathscr{W}^{\chi_j})$$

whence the result follows for $m \in W \cap S$ from Corollary 16.2. The result then follows for all $m \in W$ by linearity.  $\square$

University of Chicago, Chicago, IL 60637, USA
*E-mail address*: awb@uchicago.edu

Pennsylvania State University, State College, PA 16802, USA
*E-mail address*: hertz@math.psu.edu

Pennsylvania State University, State College, PA 16802, USA
*E-mail address*: zhirenw@psu.edu